\documentclass[11pt]{article}
\usepackage{a4wide}
\usepackage{fourier}
\usepackage[utf8]{inputenc}
\usepackage{microtype}
\usepackage{amsmath}
\usepackage{amsthm}
\usepackage{amssymb}
\usepackage{graphicx}
\usepackage{mathtools}
\usepackage{wrapfig}

\usepackage[dvipsnames,table,xcdraw]{xcolor}
\usepackage[colorlinks=true,
            citecolor=JungleGreen,
            urlcolor=MidnightBlue,
            breaklinks=true]{hyperref}   
\usepackage[shortlabels]{enumitem}
\usepackage{natbib}

\theoremstyle{plain}
\newtheorem{theorem}{Theorem}[section]       
\newtheorem{proposition}[theorem]{Proposition}
\newtheorem{definition}[theorem]{Definition}
\newtheorem{corollary}[theorem]{Corollary}
\newtheorem{lemma}[theorem]{Lemma}

\newtheorem{example}[theorem]{Example}

\usepackage{xspace}
\newcommand*{\eg}{e.g.\@\xspace}
\newcommand*{\ie}{i.e.\@\xspace}
\newcommand*{\ia}{i.a.\@\xspace}
\newcommand*{\cf}{cf.\@\xspace}
\newcommand*{\p}{p.\@\xspace}
\newcommand*{\resp}{resp.\@\xspace}
\newcommand*{\Reals}{\mathbb{R}}
\newcommand*{\Naturals}{\mathbb{N}}
\newcommand*{\Mba}{M}
\newcommand*{\Mca}{M_\sigma}
\newcommand*{\mba}{m}
\newcommand*{\mdba}{m^\circ}
\newcommand*{\co}{\overline{\mathrm{co}} \ }
\newcommand*{\Dhull}{\operatorname{D}_{\sigma}}
\newcommand*{\PreDhull}{\operatorname{D}}
\newcommand*{\Pba}{\Delta}
\newcommand*{\Pca}{\Delta_{\sigma}}
\newcommand*{\Q}{Q}
\newcommand*{\uP}{\overline{\mathbb{E}}}
\newcommand*{\lP}{\underline{\mathbb{E}}}
\newcommand*{\setofprecisegambles}{\mathcal{S}}
\newcommand*{\aP}{\nu}
\newcommand*{\pE}{E}
\newcommand*{\setoflinearprevisions}{\Delta}
\newcommand*{\fixedaP}{\psi}
\newcommand*{\mbageneral}{m}
\newcommand*{\mdbageneral}{m^\circ}
\newcommand*{\SpaceOfGambles}{\operatorname{B}(\Omega)}
\newcommand*{\SpaceOfLinearFunctionalsOnGambles}{\operatorname{ba}(\Omega)}
\newcommand*{\lin}{\operatorname{lin}}

\title{Systems of Precision: Coherent Probabilities on Pre-Dynkin-Systems and Coherent Previsions on\\ Linear Subspaces}
\author{Rabanus Derr\\ 
        {\small University of T\"{u}bingen} \\
        \texttt{\scriptsize rabanus.derr@uni-tuebingen.de }\\ 
\and 
      Robert C. Williamson\\
        {\small University of T\"{u}bingen} \\ {\small and Tübingen AI Centre} \\
        \texttt{\scriptsize bob.williamson@uni-tuebingen.de }
}

\begin{document}

\maketitle

\begin{abstract}
    In literature on imprecise probability little attention is paid
    to the fact that imprecise probabilities are precise on a set of events.
    We call these sets \emph{systems of precision}.
    We show that, under mild assumptions, the system of precision of a
    lower and upper probability form a so-called
    (pre-)Dynkin-system.
    Interestingly, there are several settings,
    ranging from machine learning on partial data over frequential probability theory
    to quantum probability theory and decision making under uncertainty,
    in which a priori the probabilities are only desired to be precise on a specific underlying set system.
    Here, (pre-)Dynkin-systems have been adopted as systems of precision, too.
    We show that, under extendability conditions, those pre-Dynkin-systems equipped with probabilities can be embedded into algebras of sets.
    Surprisingly, the extendability conditions elaborated in a strand of work in quantum probability are equivalent to 
    coherence from the imprecise probability literature. On this basis, we spell out a lattice duality which relates
    systems of precision to credal sets of probabilities.
    We conclude the presentation with a generalization of the framework to expectation-type counterparts of imprecise probabilities. The analogue of pre-Dynkin-systems turn out to be (sets of) linear subspaces in the space of bounded, real-valued functions. We introduce partial expectations, natural generalizations of probabilities defined on pre-Dynkin-systems. Again, coherence and extendability are equivalent. A related, but more general lattice duality preserves the relation between systems of precision and credal sets of probabilities.
\end{abstract}

{
\footnotesize\raggedleft\emph{When posing problems in probability calculus, \\
                 it should be required to indicate for which\\
                 events the probabilities are assumed to exist.}\\
                 ---  Andre{\u{\i}} Kolmogorov \citeyearpar[page 52]{kolmogorov1929}\\
}

\section{Introduction}
\label{introduction}
Scholarship in imprecise probability largely focuses on the imprecision of probabilities. However, imprecise probability models often
lead to \emph{precise} probabilistic statements on certain events or gambles, \ie bounded, real-valued functions.
In this work, we follow a hitherto not taken route investigating the \emph{system of precision}, \ie the set structure on which an imprecise probability is precise.\footnote{We elaborate the exact definition of imprecise probabilities and expectation used here in Section~\ref{From Porbability measures to Dynkin systesms} and Section~\ref{sec:A More General Perspective - The Set of Gambles with Precise Expectation}.}
It turns out that (pre-)Dynkin-systems\footnote{As shown in Appendix~\ref{names of dynkin-systems},
(pre-)Dynkin-systems appear under plenty of names.} describe the set of events with precise
probabilities (\cf\S~\ref{From Porbability measures to Dynkin systesms}).
This event structure is a neglected object in the literature on imprecise probability.
In particular, it constitutes a parametrized choice somewhat ``orthogonal'' to the standard.
Roughly stated, existing approaches to imprecise probability generalize the probability measure $\mu_\sigma$ in a classical probability space $(\Omega, \mathcal{F}_\sigma, \mu_\sigma)$.\footnote{
Following Kolmorogov's classical setup, $\Omega$ is the base set, $\mathcal{F}_\sigma$ a $\sigma$-algebra and $\mu_\sigma$ a countably additive probability on $\mathcal{F}_\sigma$. Approaches to imprecise probability often do not even presuppose an underlying measure space (\eg \citep{walley1991statistical}). However, they are often linked to finitely additive
measure spaces $(\Omega, \mathcal{F}, \mu)$, where $\mu$ is a finitely additive probability and $\mathcal{F}$ is an algebra of sets (sometimes called a field).}
We start by generalizing $\mathcal{F}_\sigma$ from a $\sigma$-algebra to a pre-Dynkin-system.

This suggestion is practically motivated: what do the following scenarios have in common?
\begin{description}
    \item [(a)] A machine learning algorithm
has access to a restricted subset of attributes. It cannot jointly query all attributes simultaneously.
This is called ``learning on partial, aggregated information'' \citep{eban2014discrete}.
The reasons might be manifold: for privacy preservation, ``not-missing-at-random'' features,
restricted data base access for acceleration or multi-measurement data sets.
    \item [(b)] Quantum physical
quantities, \eg location and impulse, are (statistically) incompatible \citep{gudder1979stochastic}.
    \item [(c)] A preference ordering
on a set of acts gives rise to precise beliefs on a set of events, whereas this belief
is not necessarily precise for intersections of such events \citep{Epstein2001Subjective, zhang2002subjective}.
\end{description}
In all of these scenarios, there does not exist a precise probability over all
attributes and events. Or, there is no such precise probability accessible.
Two attributes might each on their own exhibit a precise probabilistic description, 
while a joint precise probabilistic description does not exist.
On a more fundamental level, no intersectability is provided. A precise probabilistic description of
two events does \emph{not} imply that the
intersection of those events possesses a precise probability.
The set system for the description of the events
with precise probabilities which independently turned up in the various, previously mentioned fields
of research is, again, the (pre-)Dynkin-system.

\begin{wrapfigure}[7]{r}{0.2\textwidth}
    \includegraphics[width=0.2\textwidth]{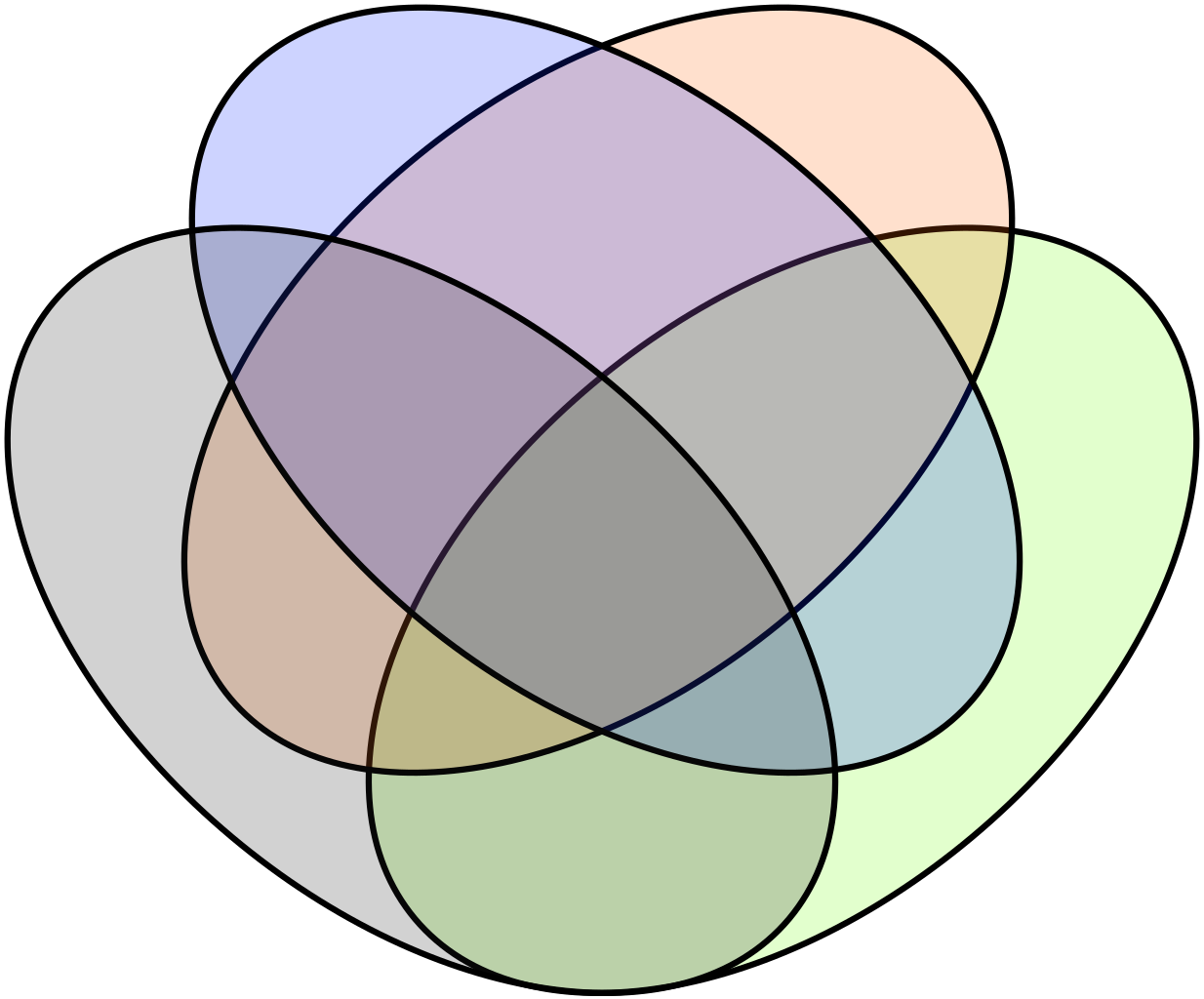}
\end{wrapfigure}
The question of intersectability (or ``intersectionality'') is of considerable interest in the social sciences
where it is used as a label to describe the problem of the \emph{joint} effect of various individual attributes
on social outcomes \citep{cole2009intersectionality, shields2008, weldon2008}. That this notion of intersectionality has 
\emph{something} to do with set systems is clear already from the fact that the Venn diagram pictured on the right\footnote{By RupertMillard, CC BY-SA 3.0.}
is used as an illustration both for the Wikipedia articles on  Hypergraphs \citep{hypergraph2023}
(another name for a set system \citep{berge1989}) and Intersectionality \citep{intersectionality2023}.
Needless to say, the concept as used in the social sciences is rich, complex, and somewhat vague, which is
not necessarily held to be a weakness:
``at least part of its success has been attributed to its vagueness'' \citep[page 260]{hancock2013}.
Our interest is in under what circumstances precise probabilities can be ascribed to events; we speculate that
such formal results may well contribute to a deeper empirical understanding of social intersectionality, without
resorting to fuzzy logic \citep{hancock2007multiplication} with 
its renowned lack of operational definition \citep{cooke2004anatomy}.

By rethinking the domain of probability measures one might wonder about the origins of Kolmogorov's $\sigma$-algebra as \emph{the} set system for events which possess probabilities.
This links back to the old problem of \emph{measurability} \citep[page 1-5]{elstrodt2018mass}. The measurability problem is the mathematical problem to assign a uniform measure to all subsets of a continuum. Giuseppe Vitali showed 1905 that this problem is not solvable for countably additive measures \citep[page 5]{elstrodt2018mass} (from \citep{vitali1905sul}). Hence, more restricted set systems such as the $\sigma$-algebra arose. \citet{isaacs2022non} reconsidered this century-old discussion to argue for rationality of imprecise probabilities. 
We take their argument even further. Inside the borders of mathematical measurability, the set of events which ought to be assigned probabilities is a modelling choice. Measurability is a modelling tool. We show that it is naturally parametrized by the set of (pre-)Dynkin-systems.

All of the preceding considerations bring us to the main question of this paper:
\textbf{What is the system of precision and how does it
relate to an imprecise probability on ``all'' events?}
We approach this question from three perspectives. 
\begin{enumerate}
    \item First, we show that, under mild assumptions,
a pair of lower and upper probabilities assign precise probabilities, \ie lower and upper probability coincide, to
events which form a pre-Dynkin-system or even a Dynkin-system.
    
    \item Second, we define probabilities
on pre-Dynkin-systems in accordance with the literature on 
quantum probability, in particular \citep{gudder1969quantum}.
We argue that probabilities on pre-Dynkin-systems, as well as their inner and outer extension,
exhibit few desirable properties, \eg subadditivity cannot be guaranteed.
Hence, extendability,
the ability to extend a probability from a pre-Dynkin-system
to a larger set structure, turns out to be crucial, as it implies
coherence of the probability defined on the pre-Dynkin-system. This 
observation links together the research from probabilities
defined on weak set structures \citep{gudder1969quantum, zhang2002subjective, schurz2008finitistic}
to imprecise probabilities \citep{walley1991statistical, augustin2014introduction}.
Furthermore, extendability guarantees the existence of a nicely behaving, so-called coherent extension.
We finally show that the inner and outer extension of a probability defined on a pre-Dynkin-system
is always more pessimistic than its corresponding lower and upper coherent extension.

    \item Last, we develop a duality theory between pre-Dynkin-systems
on a predefined base measure space and their respective credal sets of probabilities. The credal sets consist
of all probabilities which coincide with the pre-defined measure on a pre-Dynkin-system.
A so-called Galois connection links together the
containment structure on the set of set systems with
the containment structure on the set of credal sets.
\end{enumerate}
We conclude our presentation with a generalization to expectation-type counterparts of imprecise probabilities in Section~\ref{sec:A More General Perspective - The Set of Gambles with Precise Expectation}. These are often called previsions, \eg in \citep{walley1991statistical}. Our main question thus generalizes to: \textbf{What is the system of precision and how does it
relate to an imprecise expectation on ``all'' gambles?} In this case, by ``system of precision'' we mean the set of gambles with on which a lower and upper expectation coincide.
\begin{enumerate}
    \item First, we propose a generalization of a finitely additive probability defined on a pre-Dynkin-system. More concretely, we define \emph{partial expectations} which correspond to expectation functionals which are only defined on a set of linear subspaces of the space of all gambles. However, on those linear subspaces they behave like ``classical'' (finitely additive) expectations.
    \item Second, we show that under some properties, imprecise expectations are precise on a linear subspace of the linear space of gambles. (\cf Section~\ref{From Porbability measures to Dynkin systesms})
    \item Third, we present a natural generalization of extendability for partial expectations, which again turns out to be equivalent to coherence of the partial expectation.
    \item Last, analogous to the lattice duality\footnote{A lattice is a poset with pairwise existing minimum and maximum. The duality is expressed via an antitone lattice isomorphism.} described in Section~\ref{sec:the credal set and its relation to Dynkin-system structure}, we present a lattice duality for linear subsets of the space of gambles and credal sets which define coherent lower and upper previsions.
\end{enumerate}

In summary, our work makes contributions
in-between the research field of imprecise probabilities, probabilities defined on general set structures, and partially defined expectation functionals.
Part of this work has been presented on the International Symposium on Imprecise Probabilities: Theories and Applications under the title ``The Set Structure of Precision'' \citep{derr2023set}. The following version is a more complete and exhaustive presentation of this conference version. We included all omitted proofs of the conference version. We elaborated the content of Section~\ref{sec:the credal set and its relation to Dynkin-system structure}. We added an entire section about the generalization to expectation-type counterparts of probabilities on pre-Dynkin-systems (Section~\ref{sec:A More General Perspective - The Set of Gambles with Precise Expectation}). We presented relations to different research areas in more detail. We shortly discussed countably additive probabilities on Dynkin-systems in the appendix, as we put emphasis on finitely additive probabilities in the main text.
Before we begin the structural
investigation of pre-Dynkin-systems, we first introduce the used notation
and fix the mathematical framework.

\subsection{Notation and Technical Details}
\label{notation and technical details}
As we deal with a lot of sets, sets of sets, and rarely even sets of sets of sets in this paper,
we agree on the following notation: sets are written with capital latin or greek letter, \eg $A$ or $\Omega$.
Sets of sets are denoted $\mathcal{A}$. Sets of sets of sets obtain the notation $\mathfrak{A}$.
As usual, $\Reals$ is reserved for the set of real numbers, $\Naturals$ for the 
natural numbers.
The power set of a set $A$ is written as $2^A$.

In the course of this work, we require the notions of $\sigma$-algebras and algebras (of sets).
An algebra is a subset of $2^\Omega$ which contains the empty set and is closed under complement and finite union.
A $\sigma$-algebra is an algebra which is closed under countable union \citep[Definition 1.1]{williams1991probability}.\footnote{
Our notion of an algebra should not be confused with the notion of an algebra over a field.}
(Probability) measures are denoted by lowercase greek letters, \eg $\mu$, $\nu$ and $\psi$, except for $\sigma$.
Generally, we use ``$\sigma$'' to emphasize the countable nature of a mathematical object. This becomes clear when we define Dynkin-systems (Definition~\ref{def:dynkin-system}).
Other functions are denoted by lowercase latin letters, \eg $f$ and $g$.

Regarding the technical setting, we roughly follow the setup of \citep[\S 3.6 and Appendix D]{walley1991statistical}.
For a summary, see Table~\ref{tab:notation}.
\begin{table}[t]
\centering
    \begin{tabular}{r|l}
        $\Omega, 2^\Omega$ & Base set and its power set\\
        $[n]$ & Set $\{ 1, \ldots , n   \}$\\
        $\mathcal{D}$ & Pre-Dynkin-system on $\Omega$ (Definition~\ref{def:dynkin-system})\\
        $\mathcal{D}_\sigma $ & Dynkin-system on $\Omega$ (Definition~\ref{def:dynkin-system})\\
        $\PreDhull(\mathcal{A})$ & Pre-Dynkin hull of a set system $\mathcal{A} \subseteq 2^\Omega$ (Definition~\ref{def:dynkin-system})\\

        $\mu$ & Finitely additive probability defined on $\mathcal{D}$ (Definition~\ref{def:probameasure on Dynkin-system})\\
        $\mu_*$, $\mu^*$ & Inner respectively outer extension (Proposition~\ref{prop:innerouterextension})\\
        $\underline{\mu}_\mathcal{D}$, $\overline{\mu}_\mathcal{D}$ & Lower respectively upper coherent extension (Corollary~\ref{corollary:Coherent Extension of Probability})\\
        $\Mba(\mu, \mathcal{D})$ & Credal set of $\mu$ on $\mathcal{D}$ (Corollary~\ref{corollary:Coherent Extension of Probability})\\
        $\nu$ & Finitely additive probability defined on $2^\Omega$\\
        $\psi$ & Fixed, finitely additive probability defined on $2^\Omega$\\
        $\chi_A$ & Indicator function of the set $A \subset \Omega$\\
        $\Pba$ & Set of finitely additive probability measures on $2^\Omega$, set of linear previsions\\
        $\mba\colon 2^{2^\Omega} \rightarrow 2^{\Pba}$ & Credal set function (Definition~\ref{def:Credal Set Function})\\
        $\mdba\colon 2^{\Pba} \rightarrow 2^{2^\Omega}$ & Dual credal set function (Definition~\ref{def:Dual Credal Set Function})\\
        $\co$ & Convex, Weak$^\star$ Closure\\

        $\SpaceOfGambles$ & Set of real-valued, bounded functions on $\Omega$\\
        $\SpaceOfLinearFunctionalsOnGambles$ & Set of bounded, signed, finitely additive measures on $2^\Omega$\\
        $\pE$ & Partial Expectation (Definition~\ref{def:partial expectation})\\
        $\operatorname{S}(\Omega, \mathcal{A})$ & Linear space of simple gambles on the set system $\mathcal{A}$\\
        $\operatorname{B}(\Omega, \mathcal{F}_\sigma)$ & Linear space of bounded, $\mathcal{F}_\sigma$-measurable functions\\
        $\lP$ & Coherent lower prevision (Definition~\ref{def:coherent prevision})\\
        $\uP$ & Coherent upper prevision (Definition~\ref{def:coherent prevision})\\
        $\nu$ & Linear prevision defined on $\SpaceOfGambles$ (equivalent to $\nu$ above)\\
        $\psi$ & Fixed, linear prevision defined on $\SpaceOfGambles$ (equivalent to $\psi$ above)\\
        $\mba\colon 2^{\SpaceOfGambles} \rightarrow 2^{\Pba}$ & Generalized credal set function (Definition~\ref{def:generalized credal set function})\\
        $\mdba\colon 2^{\Pba} \rightarrow 2^{\SpaceOfGambles}$ & Generalized dual credal set function (Definition~\ref{def:generalized dual credal set function})
    \end{tabular}
    \caption{Summary of important, used notations.}
    \label{tab:notation}
\end{table}
Let $\Omega$ be an arbitrary set. In several examples $\Omega = [n]$, where $[n]$ denotes the set $\{ 1, \ldots , n  \}$.
The set $\SpaceOfGambles$ is defined as the set of all real-valued, bounded functions on $\Omega$. We call those functions \emph{gambles}. For instance, $\chi_A$, the indicator function of $A \subseteq \Omega$, is in $\SpaceOfGambles$.
The supremum norm, $\| f \|_{\sup} \coloneqq \sup_{\omega \in \Omega} |f(\omega)|$ makes $\SpaceOfGambles$ a topological linear vector space \citep{hildebrandt1934bounded}.
With $\SpaceOfLinearFunctionalsOnGambles$ we denote the set of all bounded, signed, finitely additive measures
on $2^\Omega$. In fact, $\SpaceOfLinearFunctionalsOnGambles$ is the topological dual space of $\SpaceOfGambles$.
So in particular, every continuous linear functional $\phi \in \SpaceOfGambles^*$ can be identified with
a bounded, signed, finitely additive measure \citep{hildebrandt1934bounded}.\footnote{
As the linear functional is defined on a normed space, continuity and boundedness are equivalent.}
For this reason, we use, with minor abuse of notation, the same notation for bounded, signed, finitely additive measures and for continuous linear functionals in the dual space of $\SpaceOfGambles$, \ie we write $\nu(f) = \int f d\nu$.
The dual space $\SpaceOfLinearFunctionalsOnGambles$ gets equipped with the weak$^\star$ topology,
\ie the weakest topology which makes all
evaluation functionals of the form $f^* \in \SpaceOfLinearFunctionalsOnGambles^*$ such that $f^*(\nu) \coloneqq \int f d\nu$ for some $f \in \SpaceOfGambles$ continuous.
With $\Pba \subseteq \SpaceOfLinearFunctionalsOnGambles$ we denote the convex, weak$^\star$-closed subset of finitely additive probability measures.
The set $\Pba$ plays a major role in Walley's theory of previsions, as the measures in $\Pba$ are in one-to-one correspondence to his linear previsions \citep[Theorem 3.2.2]{walley1991statistical}. The operator $\co$ is the convex, weak$^\star$ closure on the space $\SpaceOfLinearFunctionalsOnGambles$.
We further introduce the following two notations: let $\mathcal{F} \subseteq 2^\Omega$ be an algebra. Then $\operatorname{S}(\Omega, \mathcal{F}) \subseteq \SpaceOfGambles$ denotes the linear subspace of simple functions on $\mathcal{F}$, \ie scaled and added indicator functions of a finite number of disjoint sets (\cf \citep[Definition 4.2.12]{rao1983theory}). 
Let $\mathcal{F}_\sigma \subseteq 2^\Omega$ be a $\sigma$-algebra. Then $\operatorname{B}(\Omega, \mathcal{F}_\sigma) \subseteq \SpaceOfGambles$ denotes the linear subspace of all bounded, real-valued, $\mathcal{F}_\sigma$-measurable functions. Equipped with these notions and tools we are ready for a first preliminary question.

\section{What Is a (Pre-)Dynkin-System?}
\label{what is a dynkin-system}
In this work, the main objects under consideration are pre-Dynkin-systems and Dynkin-systems.
A (pre-)Dynkin-system is a set system on $\Omega$.
It contains the empty set, is closed under complement and (countable) disjoint union.
More formally:
\begin{definition}[(Pre-)Dynkin-system]
\label{def:dynkin-system}
    We say $\mathcal{D} \subseteq 2^\Omega$ is a \emph{pre-Dynkin-system} on some set $\Omega$ if and only
    if all of the following conditions hold:
    \begin{enumerate}[(a)]
        \item $\emptyset \in \mathcal{D}$,
        \item $D \in \mathcal{D}$ implies $D^c \coloneqq \Omega \setminus D \in \mathcal{D}$
        \item $C,D \in \mathcal{D}$ with $C \cap D = \emptyset$ implies $C\cup D \in \mathcal{D}$.
    \end{enumerate}
    We call $\mathcal{D}_\sigma \subseteq 2^\Omega$ a \emph{Dynkin-system} if and only if the conditions (a), (b) and
    \begin{enumerate}[(c')]
        \item let $\{ D_i\}_{i \in \Naturals} \subseteq \mathcal{D}_\sigma$, if for all $i,j \in \Naturals$ with $i \neq j$ it holds $D_i \cap D_j = \emptyset$ then $\bigcup_{i \in \Naturals} D_i \in \mathcal{D}_\sigma$,
    \end{enumerate}
    are fulfilled.
\end{definition}
Observe that every Dynkin-system is a pre-Dynkin-system. We will denote pre-Dynkin-systems by the use of $\mathcal{D}$, in contrast to
$\mathcal{D}_\sigma$ for Dynkin-systems. This should not be confused with $\PreDhull(\mathcal{A})$
for $\mathcal{A} \subseteq 2^\Omega$, which is the intersection 
of all pre-Dynkin-systems which contain $\mathcal{A}$, \ie the smallest pre-Dynkin-system
containing $\mathcal{A}$.\footnote{
For $\mathcal{A} = \emptyset$ we define $\PreDhull(\mathcal{A}) = \{ \emptyset, \Omega\}$.} In other words, $\PreDhull(\mathcal{A})$ is the
\emph{pre-Dynkin-hull} generated by $\mathcal{A}$.
The following short lemma will be helpful in later proofs.
\begin{lemma}[Closedness under Set Difference]
\label{lemma:closedness under set difference}
    Let $\mathcal{D} \subseteq 2^\Omega$ be a pre-Dynkin-system, if $A,B \in \mathcal{D}$ and $A \subseteq B$, then $B \setminus A \in \mathcal{D}$.
\end{lemma}
\begin{proof}
    Let $A,B \in \mathcal{D}$ and $A \subseteq B$. Then $(B \setminus A)^c = B^c \cup A$. Since $\mathcal{D}$ is closed under complement and disjoint union $B^c \cup A \in \mathcal{D}$. But then again the complement, $B \setminus A$, is in $\mathcal{D}$.
\end{proof}

In classical probability theory, Dynkin-systems appear
as a technical object required for the measure-theoretic link between
cumulative distribution functions and probability measures
(\cf\citep[Proof of Lemma 1.6]{williams1991probability}). In particular,
every $\sigma$-algebra, the well-known domain of probability measures, is
a Dynkin-system. Thus, all statements within this work are generalizations
of classical probability theoretical results. We give a short example of a pre-Dynkin-system, which is not an algebra in the following.
This example gets reused to illustrate forthcoming statements.
{
\color{darkgray}
\begin{example}
\label{runningexample:pre-Dynkin-system}
The smallest pre-Dynkin-system which is not an algebra can be defined on
$\Omega_4 \coloneqq \{ 1,2,3,4\}$. It is given by $\mathcal{D}_4 \coloneqq \{ \emptyset, 12,34,13,24, \Omega_4\}$,
where we write $12$ as a shorthand for $\{ 1,2\}$.
\end{example}
}
Pre-Dynkin- and Dynkin-systems naturally arise in probability theory.
For instance, the set of all subsets $A \subseteq \Naturals$,
such that the natural density $\mu(A) = \lim_{n \rightarrow \infty}\frac{|A \cap [n]|}{n}$
exists (\cf\citep{schurz2008finitistic}) is a pre-Dynkin-system $\mathcal{D}_\Naturals$,
but not an algebra\footnote{Intriguingly, 
    this was used as an example by \citet{kolmogorov1929} of a measure defined on a restricted set system
    for which it is desired to extend the measure to the power set $2^\Naturals$ (\cf \S~\ref{Extendable Probabilitiies on Pre-Dynkin-Systems}); see the discussion in 
    \citep[pages 11-14]{khrennikov2009interpretations} who observed (page 14) that 
    ``the main problem is \emph{non-uniqueness} of an extension'' and that
    such extended measures are impossible to verify from observed frequencies, because
    the relative frequencies do not converge for events in $2^\Naturals\setminus\mathcal{D}_\Naturals$.
    The non-uniqueness
    is naturally handled in the present paper
    by working with lower and upper previsions (or lower and upper probabilities) (\cf \citep{Frohlich2023StrictlyFI}).}. 
It is sometimes called the density logic \citep{ptak2000concrete} and constitutes
the foundation of von Mises' century-old frequential theory of probability \citep{mises1919grundlagen} (refined and summarized in \citep{mises1964mathematical}).

Another class of Dynkin-systems occurs in so-called marginal scenarios \citep{cuadras2002distributions}.
Marginal scenarios are settings in which marginal probability distributions
for a subset of a set of random variables are given, but not the entire joint distribution.
This restricted ``joint measurability'' of the involved random variables can be expressed via Dynkin-systems \citep[Example 4.2]{gudder1984extension}
\citep{vorob1962consistent}.

Pre-Dynkin-systems are so helpful because they structurally align with finitely additive probability measures.
The same statement holds for Dynkin-systems and countably additive probabilities.
If we know the probability of an event, then we know the probability of
the complement, \ie the event does not happen. 
If we know the probability of several events which are disjoint, then
we know the probability of the union, which is just the sum.
Probabilities following their standard definition go hand in hand with Dynkin-systems. We
see this observation manifested in many following statements.

Remarkably, (pre-)Dynkin-systems appeared under a variety of names (\cf Appendix~\ref{names of dynkin-systems}).
Fundamental to all its regular, independent occurences
in many research areas is the need for a set structure which does not allow for arbitrary
intersections.

\subsection{Compatibility}
\label{sec:compatibility}
(Pre-)Dynkin-systems are not necessarily closed under intersections.
However, when the intersection of two sets (events) is contained in the
(pre-)Dynkin-system, we call the two events \emph{compatible}.
\begin{definition}[Compatibility]
\label{def:compatibility}
    Let $A, B$ be elements in a pre-Dynkin-system $\mathcal{D}$, then $A$ and $B$ are \emph{compatible} if and only if $A \cap B \in \mathcal{D}$.
\end{definition}
This definition follows the definitions given in \eg \citep{gudder1969quantum, gudder1973generalized, gudder1984extension}.\footnote{
It should not be confused with the very similar, and sometimes equivalent, notion of commutativity
in logical structures \citep[Definition 14]{narens2016introduction} (\cf Appendix~\ref{appendix:From Set Systems to Logical Structures and Back}).}
Compatibility in pre-Dynkin-systems is a symmetric relation, but it is not necessarily transitive.
Furthermore, it is complement inherited, \ie if $A,B$ are compatible in a pre-Dynkin-system then so are $A, B^c$
\citep[Lemma 3.6]{gudder1979stochastic}.
Lastly, compatibility, even though expressed as intersectability, \ie ``closed under intersection'',
can be equivalently expressed as unifiability, \ie ``closed under union''.
\begin{lemma}[Cup gives Cap gives Cup]
\label{lemma:cup gives cap gives cup}
Let $\mathcal{D}$ be a pre-Dynkin-system and $A,B \in \mathcal{D}$. Then
\begin{align*}
    A \cap B \in \mathcal{D} \Leftrightarrow A \cup B \in \mathcal{D}
\end{align*}
\end{lemma}
\begin{proof}
    Using Lemma~\ref{lemma:closedness under set difference} for pre-Dynkin-systems we can quickly see that the following two decompositions give the desired equivalence:\\
    For the ``$\Rightarrow$''-direction: $A \cup B = (A \setminus (A \cap B)) \cup B$. The fact $A, A \cap B, B \in \mathcal{D}$ implies $(A \setminus (A \cap B)) \cup B \in \mathcal{D}$.\\
    For the ``$\Leftarrow$''-direction: $A \cap B = A \setminus ((A \cup B) \setminus B)$. The fact $A, A \cup B, B \in \mathcal{D}$ implies $A \setminus ((A \cup B) \setminus B) \in \mathcal{D}$.
    (A related result for Dynkin-systems is given in \citep[5.1]{gudder1969quantum}.)
\end{proof}
{
\color{darkgray}
\begin{example}
\label{runningexample:compatibility}
We reconsider the set $\Omega_4$ and pre-Dynkin-system $\mathcal{D}_4$ from Example~\ref{runningexample:pre-Dynkin-system}.
The elements $12$ and $34$ are intersectable $12 \cap 34 = \emptyset \in \mathcal{D}_4$
and unifiable $12 \cup 34 = \Omega_4 \in \mathcal{D}_4$. The elements $12$ and $13$ are not intersectable
$12 \cap 13 = 1 \notin \mathcal{D}_4$ and not unifiable $12 \cup 13 = 123 \notin \mathcal{D}_4$.
\end{example}
}
The term ``compatibility'' underlines that closedness under intersection gets loaded with further meaning
in the context of probability theory. As we define in the next section, $\mathcal{D}$ is the 
set of ``measurable'' events, \ie events which get assigned a probability. Hence, two events $A,B$ are called compatible if and only if a
precise joint probabilistic description, \ie a precise probability of $A \cap B$, exists.\footnote{For a more thorough discussion
of the nature of compatibility (and its cousin commutativity) we point to the literature on quantum probability, e.g.
\citep[Definition 3.12]{khrennikov2009contextual},
or \citep{rivas2019role}.}

Compatibility is not only a property of elements in a pre-Dynkin-systems. One can take compatibility
as a primary notion, \ie one requires the statements of Lemma~\ref{lemma:cup gives cap gives cup} and \citep[Lemma 3.6]{gudder1979stochastic} to hold.
Then, a set structure which contains the empty set and the entire base set and is equipped with this notion of compatibility is a pre-Dynkin-system \citep[Definition 5.1]{khrennikov2009interpretations}.\footnote{It is called semi-algebra in  \citep[Definition 5.1]{khrennikov2009interpretations}.}

Interestingly, the assumption of arbitrary compatibility is
fundamental to most parts of probability theory.
$\sigma$-algebras, the domain of probability measures,
are exactly those Dynkin-systems in which all events are compatible with all others \citep[Theorem 2.1]{gudder1973generalized}. Algebras are exactly those pre-Dynkin-systems in which all events are compatible with all others.
Surprisingly, it turns out that, as well, all pre-Dynkin-systems can be dissected into such ``blocks'' of full compatibility.
Every pre-Dynkin-system consists of a set of maximal algebras which we call \emph{blocks}.
In particular, maximality here stands for: there is no algebra contained
in $\mathcal{D}$ such that some $\mathcal{A}_i$ is a strict sub-algebra of this algebra.\footnote{
Similar and related results can be found in \citep{katrivnak1973certain, vsipovs1978subalgebras,brabec1979compatibility, vallander2016structure}.}
\begin{theorem}[Pre-Dynkin-Systems Are Made Out of Algebras]
\label{thm:pre-dynkin-systems are made out of algebras}
    Let $\mathcal{D}$ be a pre-Dynkin-system on $\Omega$. Then there is a unique
    family of maximal algebras $\{\mathcal{A}_i \}_{i \in I}$
    such that $\mathcal{D} = \bigcup_{i \in I} \mathcal{A}_i$. We
    call these algebras the \emph{blocks}
    of $\mathcal{D}$.
\end{theorem}
\begin{proof}
    For the proof we require the definition of a compatible subset of $\mathcal{D}$. A subset $\mathcal{A} \subseteq \mathcal{D}$ is compatible, if all elements are completely compatible, \ie any finite intersection of elements in $\mathcal{A}$ is contained in $\mathcal{D}$. This is indeed a stronger requirement than pairwise compatibility (\cf Definition~\ref{def:compatibility}). Certainly, every subset $\mathcal{A} \subseteq \mathcal{D}$ is compatible if and only if every finite subset of $\mathcal{A}$ is compatible. Hence, compatibility is a property of so-called finite character \citep[Definition 3.46]{schechter1997handbook}.
    Then, Tuckey's lemma (\eg \citep[Theorem 6.20.AC5]{schechter1997handbook}) guarantees that any compatible subset of $\mathcal{D}$ is contained in a maximal compatible subset. Since every element $D \in \mathcal{D}$ is in at least one compatible subset, \eg
    $\{\emptyset, D, D^c, \Omega\} \subseteq \mathcal{D}$, the (unique) set of maximal compatible subsets $\{\mathcal{A}_i \}_{i \in I}$ covers the entire pre-Dynkin-system $\mathcal{D}$. It remains to show that the maximal compatible subsets are algebras.
    Consider a maximal compatible subset $\mathcal{A}_i$. First, $\emptyset \in \mathcal{A}_i$ as $\emptyset$ is compatible to all sets in $2^\Omega$. Second, $\mathcal{A}_i$ is closed under finite intersection, otherwise there would exist a finite combination of elements $A_1, \ldots, A_n \subseteq \mathcal{A}_i$ such that $A_\cap \coloneqq \bigcap_{j = 1}^n A_j \in \mathcal{D}$, but  $A_\cap \notin \mathcal{A}_i$. Then, one can easily see that $\mathcal{A}_i \cup \{ A_\cap\}$ would be a compatible subset which strictly contained $\mathcal{A}_i$. This is impossible, since $\mathcal{A}_i$ is maximal. Finally, $\mathcal{A}_i$ is closed under complement. Consider $A \in \mathcal{A}_i$, we show that $\mathcal{A}_i \cup \{ A^c\}$ is again a compatible subset. Let $A_1,\ldots , A_n \subseteq \mathcal{A}_i$ be an arbitrary finite collection of subsets, then $A \cap \bigcap_{j =1}^n A_n \in \mathcal{D}$, hence $A^c \cap \bigcap_{j =1}^n A_n \in \mathcal{D}$ \citep[Lemma 3.6]{gudder1979stochastic}. By maximality of $\mathcal{A}_i$ we then know $A^c \in \mathcal{A}_i$.
\end{proof}
{
\color{darkgray}
\begin{example}
\label{runningexample:predynkinsystem into blocks}
The pre-Dynkin-system $\mathcal{D}_4$ of Example~\ref{runningexample:pre-Dynkin-system} consists of the algebras $\{ \emptyset, 12,34, \Omega_4\}$ and $\{ \emptyset, 13, 24, \Omega_4\}$.
\end{example}
}
Theorem~\ref{thm:pre-dynkin-systems are made out of algebras} simplifies several follow-up observations. Instead of pre-Dynkin-systems we can equivalently
consider a set of algebras. However, not every union of algebras is a pre-Dynkin-system.
If these algebras form a compatibility structure, \ie a set of
maximal $\pi$-systems\footnote{Non-empty set systems which are closed under finite intersections are called $\pi$-systems.}, then their union is a pre-Dynkin-system (Definition~\ref{def:compatibility structure} and Theorem~\ref{theorem:Union of Algebras is Pre-Dynkin-system if Compatibility Structure} in Appendix).
Analogous results for Dynkin-systems and $\sigma$-algebras exist and are given in Appendix~\ref{appendix:Dynkin-Systems}.
In summary, pre-Dynkin-systems are set structures which do not allow for arbitrary intersections, but can be
split into maximal intersectable subsets, their \emph{blocks}.

\subsection{Probabilities on Pre-Dynkin-Systems}
We require a notion of probability on a pre-Dynkin-system to elaborate the relationship
of imprecise probability and the system of precision in the following.
Probabilities are classically defined on $\sigma$-algebras.
We generalize this definition as \eg stated in \citep[page 18f]{williams1991probability}
to pre-Dynkin-systems.
\begin{definition}[Probability Measure on a Pre-Dynkin-System]
\label{def:probameasure on Dynkin-system}
    Let $\mathcal{D}$ be a pre-Dynkin-system. We call a function
    $\mu \colon \mathcal{D} \rightarrow [0,1]$  a \emph{countably additive probability measure on $\mathcal{D}$} if and only if it fulfills the following two conditions:
    \begin{enumerate}[(a)]
        \item \emph{Normalization}: $\mu(\emptyset) = 0$ and $\mu(\Omega) = 1$.
        \item $\sigma$\emph{-Additivity}: let $I \subseteq \Naturals$ and $\{ A_i\}_{i \in I}$ such that $A_i \in \mathcal{D}$ for all $i \in I$ and $A_i \cap A_j = \emptyset$ for $i \neq j$, $i,j \in I$. Then
        $\mu(\bigcup_{i \in I} A_i) = \sum_{i \in I} \mu(A_i)$.
    \end{enumerate}
    If condition (b) holds at least for finite $I$, we say that $\mu$ is a \emph{finitely additive probability measure}.
\end{definition}
For the sake of readability, we use ``probability'' and ``probability measure'' exchangeably.
Probabilities on pre-Dynkin-systems are monotone, \ie for $A, B \in \mathcal{D}$, if $A \subseteq B$, then $\mu(A) \le \mu(B)$.\footnote{This can be seen when applying Lemma~\ref{lemma:closedness under set difference} and Definition~\ref{def:probameasure on Dynkin-system}.}
But, in contrast to a probability defined on a $\sigma$-algebra,
a probability on a pre-Dynkin-system is not necessarily \emph{modular}, \ie for $A, B \in \mathcal{D}$, $\mu(A) + \mu(B) = \mu(A \cup B) + \mu(A \cap B)$ \citep[page 16]{denneberg1994non}.\footnote{
It is, however, possible to define modular probabilities on pre-Dynkin-systems. This leads to
a fixed parametrization of probability functions already on simple examples \citep[page 125]{navara1998considering}.}
It is that sophisticated interplay of set structure and probability function which leads us
through this paper. In particular, why should we consider pre-Dynkin-systems?

\section{Imprecise Probabilities Are Precise on a Pre-Dynkin-System}
\label{From Porbability measures to Dynkin systesms}
As we now demonstrate, pre-Dynkin-systems are, under mild assumptions, the systems of precision.
To make this formal, we solely require a normed, conjugate pair of lower and upper probability
which fulfill super\@\xspace(\resp sub)-additivity and possibly a continuity assumption.
\begin{theorem}[Imprecise Probability Induces a (Pre-)Dynkin-System]
\label{thm:lower and upper probablity define Dynkin-system}
    Let $\ell \colon 2^\Omega \rightarrow [0,1]$ and $u \colon 2^\Omega \rightarrow [0,1]$ be 
    two set functions, for which all the following properties hold:
    \begin{enumerate}[(a)]
        \item Normalization: $u(\emptyset) = \ell(\emptyset) = 0$.
        \item Conjugacy: $u(A) = 1 - \ell(A^c)$ for $A, A^c \in 2^\Omega$.
        \item Subadditivity of $u$: for $A,B \in 2^\Omega$ such that $A \cap B = \emptyset$
        then $u(A \cup B) \le u(A) + u(B)$.
        \item Superadditivity of $\ell$: for $A,B \in 2^\Omega$ such that $A \cap B = \emptyset$
        then $\ell(A \cup B) \ge \ell(A) + \ell(B)$.
    \end{enumerate}
     Then $u$ and $\ell$ define a finitely additive probability measure
    $\mu \coloneqq u|_{\mathcal{D}} = \ell|_{\mathcal{D}}$ on a pre-Dynkin-system
    $\mathcal{D} \subseteq 2^\Omega$.
    If either $u$ fulfills
    \begin{enumerate}[(e)]
        \item Continuity from below: for $A_n \in 2^\Omega$ with $A_n \subseteq A_{n+1}$
        such that 
        $\bigcup_{n= 1}^\infty A_n = A \in 2^\Omega$, then\\
        $\lim_{n \rightarrow \infty} u(A_n) = u(A)$,
    \end{enumerate}
    or $\ell$ fulfills
    \begin{enumerate}[(e')]
        \item Continuity from above: for $A_n \in 2^\Omega$ with $A_{n+1} \subseteq A_{n}$
        such that 
        $\bigcap_{n= 1}^\infty A_n = A \in 2^\Omega$, then\\
        $\lim_{n \rightarrow \infty} \ell(A_n) = \ell(A)$,
    \end{enumerate}
    then $u$ and $\ell$ define a countably additive probability measure
    $\mu_\sigma \coloneqq u|_{\mathcal{D}_\sigma} = \ell|_{\mathcal{D}_\sigma}$ on a Dynkin-system
    $\mathcal{D}_\sigma \subseteq 2^\Omega$.
\end{theorem}
\begin{proof}
    We start proving the first part of the theorem.
    Let
    \begin{equation}
    \label{eq:sets with precise probability}
        \mathcal{D} \coloneqq \{A \in 2^\Omega \colon  \ell(A) = u(A)\}.
    \end{equation}
    We show that $\mathcal{D}$ is a pre-Dynkin-system.
    First, $\emptyset \in \mathcal{D}$ by assumption (a).
    Second, let $D \in \mathcal{D}$. Then $u(D^c) = 1 - \ell(D) = 1 - u(D) = \ell(D^c)$ by the conjugacy relation.
    Third, let $\{ A_i\}_{i \in I} \subseteq \mathcal{D}$ for finite $I \subseteq \Naturals$
    such that $A_i \cap A_j = \emptyset$ for all $i \neq j$, then
    \begin{align*}
        \sum_{i \in I} \ell(A_i) &\overset{(d)}{\le} \ell \left(\bigcup_{i \in I} A_i \right)\\
        &\overset{(\star)}{\le} u\left(\bigcup_{i \in I} A_i\right)\\
        &\overset{(c)}{\le} \sum_{i \in I} u(A_i)\\
        &\overset{Eq.\ref{eq:sets with precise probability}}{=} \sum_{i \in I} \ell(A_i).
    \end{align*}
    For $(\star)$, observe that $\ell(A) \le u(A)$ for all $A \in 2^\Omega$, since
    \begin{align*}
        \ell(A) + \ell(A^c) &\le \ell(A \cup A^c) = 1 = u(A \cup A^c) \le u(A) + u(A^c),
    \end{align*}
    and thus,
    \begin{align*}
        \ell(A) + \ell(A^c) &\le u(A) + u(A^c)\\
        \Leftrightarrow \ell(A) + 1 - u(A) &\le u(A) + 1 - \ell(A)\\
        \Leftrightarrow \ell(A) &\le u(A).
    \end{align*}
    Concluding, we define $\mu \coloneqq \ell|_\mathcal{D} = u|_\mathcal{D}$ for which it is trivial to show that
    it is a finitely additive probability on $\mathcal{D}$.
    
    For the second part, we first notice that continuity from below and from above are equivalent
    for conjugate set functions
    on set systems which are closed under complement \citep[Proposition 2.3]{denneberg1994non}.
    Next, we show that subadditivity of $u$
    and continuity from below (of $u$) imply
    $\sigma$-subadditivity of $u$: for $\{ A_i\}_{i \in I} \subseteq 2^\Omega$
    such that $I \subseteq \Naturals$ and $A_i \cap A_j = \emptyset$ for all $i \neq j$ with $i, j \in I$
    then $u\left( \bigcup_{i \in I} A_i \right ) \le \sum_{i \in I}u(A_i)$.
    In case that $I$ is finite, subadditivity of $u$ is provided by assumption. For infinite $I$
    we can construct an increasing sequence of sets, namely
    $B_j = \bigcup_{1 \le i \le j} A_i$, so that $B_j \subseteq B_{j+1}$.
    Furthermore, $\bigcup_{j = 1}^\infty B_j = \bigcup_{i \in I} A_i$. Thus,
    \begin{align*}
        u\left( \bigcup_{i \in I} A_i \right ) &= u\left( \bigcup_{j = 1}^\infty B_j \right) \\
        &\overset{(e)}{=} \lim_{j \rightarrow \infty} u\left(B_j \right)\\
        &= \lim_{j \rightarrow \infty} u\left(\bigcup_{1 \le i \le j} A_i \right) \\
        &\overset{(d)}{\le} \lim_{j \rightarrow \infty} \sum_{1 \le i \le j} u\left( A_i \right) \\
        &= \sum_{i \in I} u\left( A_i \right).
    \end{align*}
    The same argument holds analogously for superadditivity and continuity from above of $\ell$ which is implied by continuity from below 
    and the conjugacy relationship \citep[Proposition 2.3]{denneberg1994non}.
    In summary, the proof of the first part can then be applied again, now
    without the restriction that $I \subseteq \mathbb{N}$ is finite. Instead it potentially is countable.
\end{proof}
\begin{figure}
    \centering
    \includegraphics[width=\textwidth]{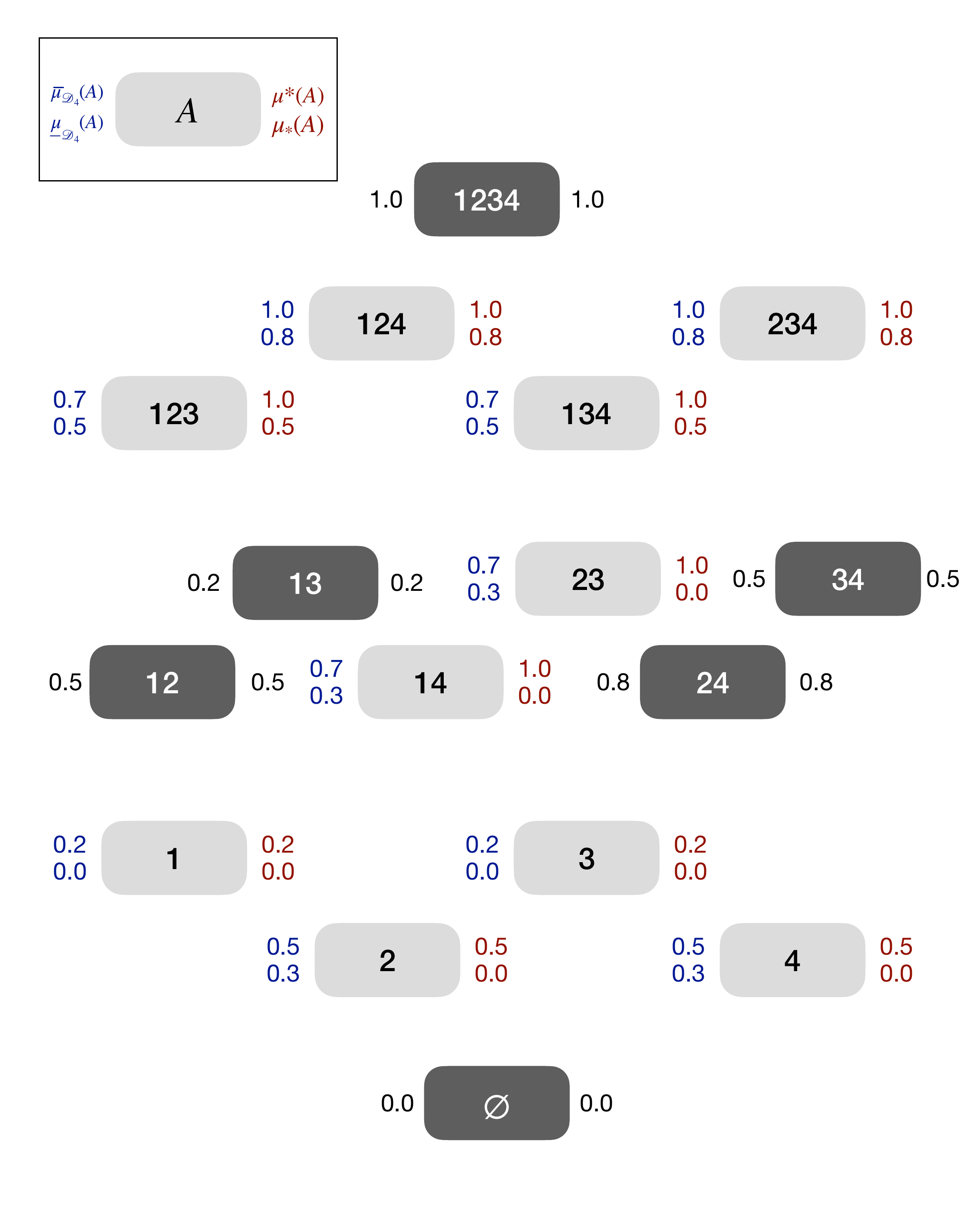}
    \caption{Illustration of the running example. The dark elements are contained in the pre-Dynkin-system $\mathcal{D}$ on $\Omega = \{ 1,2,3,4\}$. The lower and upper coherent extension, respectively the inner and outer extension are denoted at the sides of the elements in the set system as shown in the example in the left upper corner. Elements in $\mathcal{D}$ possess a precise probability.}
    \label{fig:runningexample}
\end{figure}
{
\color{darkgray}
\begin{example}
\label{runningexample:imprecise probability give pre-Dynkin-system}
Remember, $\Omega_4 = \{ 1,2,3,4\}$. We define
$\ell\colon 2^{\Omega_4} \rightarrow [0,1]$ by $\ell(1) = 0$, $\ell(2) = 0.3$, $\ell(3) = 0$, $\ell(4) = 0.3$, $\ell(12) = 0.5$, $\ell(34) = 0.5$, $\ell(13) = 0.2$, $\ell(24) = 0.8$, $\ell(14) = 0.3$, $\ell(23) = 0.3$, $\ell(123) = 0.5$, $\ell(124) = 0.8$, $\ell(134) = 0.5$, $\ell(234) = 0.8$, $\ell(\Omega_4) = 1$
and $u\colon 2^{\Omega_4} \rightarrow [0,1]$ by $u(A) = 1- \ell(A^c)$.
It is easy to show that $\ell$ and $u$ fulfill the assumptions (a), (b), (c) and (d) in
Theorem~\ref{thm:lower and upper probablity define Dynkin-system}.
The imprecise probabilities $u$ and $\ell$ coincide on $\{\emptyset, 12, 34, 13,24,\Omega_4\}$, the pre-Dynkin-system described in Example~\ref{runningexample:pre-Dynkin-system}.
The example is illustrated in Figure~\ref{fig:runningexample}. In this figure $\ell = \underline{\mu}_{\mathcal{D}_4}$ and  $u = \overline{\mu}_{\mathcal{D}_4}$.
\end{example}
}

In summary, imprecise probabilities are, under mild assumptions, precise on a pre-Dynkin-system or even a Dynkin-system. This, importantly, is also the case if the system of precision is strictly larger than the trivial pre-Dynkin-systems $\{ \emptyset, \Omega\}$.
Exemplarily, a pair of conjugate, coherent lower and upper probability (\eg \citep[\S 2.7.4]{walley1991statistical}) fulfills the conditions (a) -- (d). However, in several cases (\eg distorted probability distributions)
imprecise probabilities are just precise on the \emph{system of certainty}, \ie the events which possess $0$ or $1$ probability (Proposition~\ref{prop:events of precise probability for distorted probabilities}).
Concluding, the system of precision is a pre-Dynkin-system $\mathcal{D} \subseteq 2^\Omega$. 
What if we first define precise, finitely additive probabilities
on a pre-Dynkin-system, \ie we fix a system of precision? We can then ask for ``imprecise probabilities'' deduced from this probability which are defined on a larger set structure, \eg
an algebra in which the pre-Dynkin-system is contained.

\section{Extending Probabilities on Pre-Dynkin-Systems}
\label{Extendable Probabilitiies on Pre-Dynkin-Systems}
Precise probabilities on pre-Dynkin-system naturally arise in many, distinct, applied scenarios as we argued in the Introduction (\S~\ref{introduction}).
However, we acknowledge that the definition of probabilities on
pre-Dynkin-systems is mathematically cumbersome.
The possibilities to prove standard theorems is very limited
as the approaches by \citet{gudder1973generalized, gudder1979stochastic, gudder1981generalized} demonstrate.
However, if we consider a probability defined on a pre-Dynkin-system as an imprecise probability on a larger set system with a fixed system of precision,
we possibly obtain a richer, mathematical toolkit to work with. In this case the larger set system preferably is an algebra in which the pre-Dynkin-system is contained.
It remains to clarify how we construct the imprecise probability from the precise probability on the pre-Dynkin-sytem.

\subsection{Inner and Outer Extension}
A simple but, as we show, unsatisfying solution is the use of an inner and outer measure extension.
It does not rely on imposing any conditions on the
probability defined on the pre-Dynkin-system.
We pay for this generality with the few
properties that we can derive
for the obtained extension.
\begin{proposition}[Inner and Outer Extension]\citep[Lemma 2.2]{zhang2002subjective}
\label{prop:innerouterextension}
    Let $\mathcal{D}$ be a pre-Dynkin-system on $\Omega$ and
    $\mu$ a finitely additive probability measure
    on $\mathcal{D}$. The \emph{inner probability measure}
    \begin{align*}
        \mu_*(A) &:= \sup \{ \mu(B): A \supseteq B \in \mathcal{D} \}, \ \forall A \in 2^\Omega,
    \end{align*}
    and \emph{outer probability measure}
    \begin{align*}
        \mu^*(A) &:= \inf \{ \mu(B): A \subseteq B \in \mathcal{D} \}, \ \forall A \in 2^\Omega,
    \end{align*}
    define $\mu_*, \mu^*: 2^\Omega \rightarrow [0,1]$, i.e.
    all of the following conditions are fulfilled:
    \begin{enumerate}[(a)]
        \item Normalization: $\mu^*(\emptyset) = 0$, $\mu_*(\Omega) = 1$.
        \item Conjugacy: $\mu^*(A) = 1 - \mu_*(A^c), \ \forall A \in 2^\Omega$.
        \item Monotonicity: for $A, B \in 2^\Omega$, if $A \subseteq B$ then $\mu^*(A) \le \mu^*(B)$.
    \end{enumerate}
    Furthermore, $\mu_*$ is superadditive, for $A, B \in 2^\Omega$ if $A \cap B =\emptyset$ then $\mu_*(A \cup B) \ge \mu_*(A) + \mu_*(B)$.
    But $\mu^*$ is not generally subadditive.
\end{proposition}
{
\color{darkgray}
\begin{example}
\label{runningexample:inner and outer extension}
For $\mathcal{D}_4$ as in Example~\ref{runningexample:imprecise probability give pre-Dynkin-system} let $\mu \colon \mathcal{D}_4 \rightarrow [0,1]$ be defined as $\mu(\emptyset) = 0, \mu(12) = 0.5, \mu(34) = 0.5, \mu(13) = 0.2, \mu(24) = 0.8, \mu(\Omega) = 1$.
The inner and outer extension of $\mu$ on $\mathcal{D}_4$
is $\mu_*(\emptyset) = 0, \mu_*(1) = \mu_*(2) = \mu_*(3) = \mu_*(4) = 0, \mu_*(12) = 0.5,\mu_*(34) = 0.5,\mu_*(13) = 0.2,\mu_*(24) = 0.8,\mu_*(14) = 0,\mu_*(23) = 0,\mu_*(123) = 0.5, \mu_*(124) = 0.8, \mu_*(134) = 0.5, \mu_*(234) = 0.8, \mu_*(\Omega_4) = 1$ and $\mu^* = 1 - \mu_*$.
The inner and outer extension are not coherent (Definition~\ref{def:coherent probability}).
In particular, the outer extension is not subadditive: $\mu^*(14) = 1 - \mu_*(23) = 1 > 0.2 + 0.5 = (1 - \mu_*(234)) + (1- \mu_*(123)) = \mu^*(1) + \mu^*(4)$.
The example is illustrated in Figure~\ref{fig:runningexample}.
\end{example}
}
In conclusion, the inner and outer extension provides an imprecise probability, which is not necessarily coherent (\cf Definition~\ref{def:coherent probability}) and it
does not fulfill the conditions required for Theorem~\ref{thm:lower and upper probablity define Dynkin-system}
to post-hoc guarantee that the set of precision is a pre-Dynkin-system. We remark that there exist normalized, conjugate, monotone superadditive but not subadditive pairs of probabilities, hence possibly inner and outer probabilities as defined here, whose system of precision is not a pre-Dynkin-system (see Example~\ref{example:almost predynkin-system}). For this reason we now explore another, more powerful extension method.
{
\color{darkgray}
\begin{example}
    \label{example:almost predynkin-system}
    Let $\Omega_3 \coloneqq \{ 1,2,3\}$. The probability pair defined by $\underline{\nu}(\emptyset) = \overline{\nu}(\emptyset) = 0$, $\underline{\nu}(1) = \overline{\nu}(1) = 0.2$,  $\underline{\nu}(2) = \overline{\nu}(2) = 0.2$ and $\underline{\nu}(3) = 0 , \overline{\nu}(3) = 0.6$ is precise on
    $\{ \emptyset, 1, 2, 23, 13, \Omega_3\}$, which is obviously not a pre-Dynkin-system.
\end{example}
}

\subsection{Extendability and Its Equivalence to Coherence}
In the following, we try to entirely embed pre-Dynkin systems equipped with a probability into
larger algebras. Then, we extend the probability defined on
the pre-Dynkin-system in all possible ways to probabilities on the algebra.
It turns out that this embedding is only possible under certain conditions on the probability defined on the pre-Dynkin-system.
We call this condition \emph{extendability}.
For the sake of generality, we focus on the extension of finitely additive
probabilities from pre-Dynkin-systems to algebras here. We treat countably additive probabilities, Dynkin-systems and $\sigma$-algebras
in Appendix~\ref{appendix:Dynkin-Systems and Countably Additive Probability}.
In addition, all results until Subsection~\ref{sec:extensions of probabilities on pre-dynkin-systems}
can be formulated in more general terms for non-structured set systems. For the sake of simplicity,
we remain within the setting of probabilities defined on pre-Dynkin-systems in this work.

Extendability is the property that a probability measure defined on
a pre-Dynkin-system can be extended to a probability measure on an
algebra containing the pre-Dynkin-system. Formally:
\begin{definition}[Extendability]
\label{def:extendability}
    Let $\mathcal{D}$ be a pre-Dynkin-system on $\Omega$.
    We call a finitely additive probability measure $\mu$
    on $\mathcal{D}$ \emph{extendable} to $2^\Omega$
    if and only if there is a finitely additive probability measure $\nu\colon 2^\Omega \rightarrow [0,1]$ such that $\nu|_\mathcal{D} = \mu$.
\end{definition}
We defined extendability with respect to the power set $2^\Omega$. In fact, any relativization to an arbitrary sub-algebra of $2^\Omega$ is equivalent. A finitely additive probability defined on $\mathcal{D}$ is extendable to any sub-algebra of $2^\Omega$ which contains $\mathcal{D}$ if and only if it is extendable to $2^\Omega$ \citep[Theorem 3.4.4]{rao1983theory}.

The definition is non-vacuous \citep{gudder1984extension, de2007extending}.
For instance, a probability measure on a pre-Dynkin-System is not generally
extendable to a measure on the generated algebra (\eg Example 3.1 in \citep{gudder1984extension}).
If a probability is extendable, its extension is in general non-unique.

Extendability of probabilities on (pre-)Dynkin-systems
has already been part of discussions in quantum probability since 1969 \citep{gudder1969quantum}
up to more current times \citep{de2010measures}.
Several necessary and/or sufficient conditions on the structure of $\mathcal{D}$ and/or the values of $\mu$ are known
\citep{gudder1984extension, de2007extending, de2010measures}. We present here a sufficient and necessary condition discovered by \citet{horn1948measures} and restated in \citep[Theorem 3.2.10]{rao1983theory}.\footnote{In fact, Theorem~\ref{theorem:extendability condition - Horn-Tariski} can be stated for a more general definition of probabilities on arbitrary set systems \eg \citep[Theorem 3.2.10]{rao1983theory} \citep[Propostion 2.2]{Simone2006ExtendingCM}. For the sake of simplicity, we restrict this result to pre-Dynkin-systems and probabilities defined on pre-Dynkin-systems.}
\begin{theorem}[Extendability Condition]\citep[Theorem 3.2.10]{rao1983theory}
\label{theorem:extendability condition - Horn-Tariski}
    Let $\mathcal{D}$ be a pre-Dynkin-system on $\Omega$.
    A finitely additive probability measure $\mu$
    on $\mathcal{D}$ is extendable to $2^\Omega$ if and only if
    \begin{align}
    \label{eq:extendability condition}
    \sum_{k = 1}^m \chi_{B_k}(\omega) - \sum_{j = 1}^n \chi_{A_j}(\omega) \ge 0, \ \forall \omega \in \Omega \implies \sum_{k = 1}^m \mu(B_k) - \sum_{j = 1}^n \mu(A_j) \ge 0
    \end{align}
    for all finite families of sets in $\mathcal{D}$: $A_1, \ldots, A_n, B_1, \ldots, B_m \in \mathcal{D}$.
\end{theorem}
{
\color{darkgray}
\begin{example}
\label{runningexample:extendable probability on pre-Dynkin-system}
For $\mathcal{D}_4$ as in Example~\ref{runningexample:imprecise probability give pre-Dynkin-system} let $\mu \colon \mathcal{D}_4 \rightarrow [0,1]$ be defined as $\mu(\emptyset) = 0, \mu(12) = 0.5, \mu(34) = 0.5, \mu(13) = 0.2, \mu(24) = 0.8, \mu(\Omega) = 1$. The probability $\mu$ on $\mathcal{D}_4$ meets the extendability condition.
\end{example}
}
Extendability proves to be more than 
a helpful mathematical property for embedding pre-Dynkin-systems and
their respective probabilities into algebras.
Whether a probability defined on $\mathcal{D}$ can be extended to a
probability on $2^\Omega$ is directly connected to the question whether the probability
measure on $\mathcal{D}$ is coherent in the sense of \citep[page 68, page 84]{walley1991statistical} or not.
Coherence is a minimal consistency requirement for probabilistic descriptions which has been introduced in the fundamental work of \citet{de2017theory} and developed by \citet{walley1991statistical}.
Shortly summarizing, an incoherent imprecise probability
is tantamount to an irrational betting behavior, thus the name. Thus, extendability is, besides its mathematical convenience, a desirable property of probabilities in pre-Dynkin settings.

We adapt here the definition of coherence of previsions in \citep[Definition 2.5.1]{walley1991statistical}
to probabilities.
\begin{definition}[Coherent Probability]
    \label{def:coherent probability}
    Let $\mathcal{A} \subseteq 2^\Omega$ be an arbitrary collection of subsets.
    A set function $\underline{\nu} \colon \mathcal{A} \rightarrow [0,1] $ is a \emph{coherent lower probability} if and only if
    \begin{align*}
        \sup_{\omega \in \Omega} \sum_{i = 1}^j (\chi_{A_i}(\omega) - \underline{\nu}(A_i)) - m (\chi_{A_0}(\omega) - \underline{\nu}(A_0)) \ge 0,
    \end{align*}
    for any non-negative $n,m \in \Naturals$ and any $A_0, A_1, \ldots A_n \in \mathcal{A}$.
    If $\mathcal{A}$ is closed under complement, the conjugate \emph{coherent upper probability} is given by $\overline{\nu}(A) \coloneqq 1 - \underline{\nu}(A^c)$ for all $A \in \mathcal{A}$.
    If furthermore $\overline{\nu}(A) = \underline{\nu}(A)$ for all $A \in \mathcal{A}$, we call $\nu \coloneqq \overline{\nu}$ a \emph{coherent additive probability}.
\end{definition}
At first sight, the Horn-Tarski condition given in Theorem~\ref{theorem:extendability condition - Horn-Tariski} and the coherence condition presented here already appear similar.
This becomes even more apparent in Walley's reformulation of coherence for additive probabilities \citep[Theorem 2.8.7]{walley1991statistical}.
In the following, we show that this superficial similarity is indeed based on a rigorous link. Surprisingly, Walley did not mention Horn and Tarski's work in his foundational book.
\begin{theorem}[Extendability Equals Coherence]
\label{thm:Extendability Gives Coherence}
    Let $\mathcal{D}$ be a pre-Dynkin-system on $\Omega$.
    A finitely additive probability measure $\mu$
    on $\mathcal{D}$ is extendable to $2^\Omega$ if and only if it is a coherent additive probability on $\mathcal{D}$.
\end{theorem}
\begin{proof}
    If $\mu$ is a coherent additive probability on $\mathcal{D}$, then the linear extension theorem \citep[Theorem 3.4.2]{walley1991statistical} applies.
    Hence, a coherent additive probability $\nu \colon 2^\Omega \rightarrow [0,1]$ exists, such that $\nu|_\mathcal{D} = \mu$. In particular, $\nu$ is a finitely additive probability following Definition~\ref{def:probameasure on Dynkin-system} on $2^\Omega$ \citep[Theorem 2.8.9]{walley1991statistical}.

    For the converse direction, we observe that if $\mu$ possess an extension following Definition~\ref{def:extendability}, then such an extension is a finitely additive probability on $2^\Omega$ following Definition~\ref{def:probameasure on Dynkin-system}. Hence, \citet[Theorem 2.8.9]{walley1991statistical} guarantees that the extension is a coherent additive probability (Definition~\ref{def:coherent probability}). Any restriction to a subdomain $\mathcal{D} \subseteq 2^\Omega$ keeps the probability coherent and additive.
\end{proof}
The linear extension theorem in Walley \citep[Theorem 3.4.2]{walley1991statistical} used here is a generalization of de Finetti's fundamental theorem of probability \citep[Theorems 3.10.1 and 3.10.7]{de2017theory}.
De Finetti's theorem is furthermore interesting, as he explicitly states that a coherent additive probability defined on an arbitrary collection of sets can be extended in a precise way (so lower and upper probability coincide) to \emph{some} sets. De Finetti does not characterize this collection. Our Theorem~\ref{thm:lower and upper probablity define Dynkin-system}, however, gives an answer to this question: the collection forms a pre-Dynkin-system.

Theorem~\ref{thm:Extendability Gives Coherence}
provides a missing link between two strands of work: 
on the one hand, probabilities on pre-Dynkin-systems and related weak set structures
have been closely investigated in foundational quantum probability theory \citep{gudder1969quantum, gudder1979stochastic} and decision theory \citep{Epstein2001Subjective,zhang2002subjective}. On the other hand,
coherent probabilities are central to imprecise probability, in particular, the more general formulations of
coherent previsions and risk measures \citep{walley1991statistical, delbaen2002coherent, pelessoni2003imprecise}.

The reader familiar with the literature on imprecise probability might well not be surprised by the equivalence of extendability and coherence. We still think that this link is indeed valuable to be spelled out explicitly here. The concept of extendability and coherence have been developed separately in two communities with different goals in focus. Coherence tries to capture ``rational'' betting behavior \citep{de2017theory, walley1991statistical}. Extendability links to what is sometimes called ``quantum weirdness''.

\subsubsection{Extendability, Compatibility and Contextuality}
Extendability in quantum theory tightly interacts with a series of properties and concepts which pervade discussions about the ``specialness'' of quantum theory in comparison to other classical physical theories: compatibility, contextuality, hidden variables and more. To be concrete, two measurements are \emph{compatible} if, for any initial state\footnote{States are often represented as probability distributions in foundational quantum theory (\cf \citep{gudder1979stochastic}).}, there exists a joint measurement such that a fixed joint distribution for both measurement outcomes exists, whose marginals are the distribution of the single measurement \citep{Busch2012ComparingTD, xu2018NecessaryAS}. If measurements are incompatible, then there are potentially still states such that a joint distribution of measurements exist. Only in the cases that no joint distribution of measurements exists, \ie extendability is not provided, a measuring observer observes \emph{contextual} behavior \citep{xu2018NecessaryAS}. Translated to the language of imprecise probability, contextuality amounts to non-coherence of a probabilistic description. Compatibility, in contrast, is a structural notion. If any finitely additive probability on a pre-Dynkin-system is extendable, then the pre-Dynkin-system, very roughly, resembles compatible measurements. We are indeed not the first to notice intriguing links between imprecise probability and concepts therein to quantum mechanics. Benavoli and collaborators recovered the four postulates of quantum mechanics with desirability as a starting point \citep{Benavoli2016QuantumMT}. Desirability is a very general framework for imprecise probability \citep{walley2000towards}.

\subsubsection{Extendability and Marginal Scenarios}
Not far from the relation between extendability and coherence, \citet{Miranda2018Compatibility} and \citet{casanova2022information} bridged desirability to marginal scenarios. Marginal scenarios can equivalently be expressed in terms of probabilities on pre-Dynkin-systems \citep{vorob1962consistent, kellerer1964verteilungsfunktionen}\citep[Example 4.2]{gudder1984extension}. In rough terms, the marginal problem for marginal scenarios asks whether for a given set of marginal probability distributions (not necessarily disjoint) there exists a joint distribution.\footnote{The attentive reader might have noticed the similarity to the notion of compatibility. For good reasons \citep[\S V.B.2]{Budroni2021KochenSpeckerC}. Marginal scenarios are used to represent multi-measurement settings  and compatibility among the measurements.} This question has been, some while ago, asked for probabilities on finite spaces \citep{vorob1962consistent}, countably additive probabilities \citep{kellerer1964verteilungsfunktionen}, finitely additive probabilities (\cf \citep{maharam1972consistent}) and recently for even more general probability models -- sets of desirable gambles \citep{Miranda2018Compatibility}.
A recurring theme in all those studies is the so-called running intersection property which characterizes all those marginal structures for which a joint probabilistic description can always be guaranteed.
To bring extendability to this picture, one should think of it as a more fine-grained concept: solutions to the marginal problem show under which circumstances \emph{every} marginal distribution of a certain structure is extendable. But there exist marginal problems for which only \emph{specific} instantiations of the marginal distributions allow for extendability. The running intersection property is a property of a structure. Extendability is a property of a structure and a probability on this structure.

\subsection{Coherent Extension}
\label{sec:extensions of probabilities on pre-dynkin-systems}
A probability on a pre-Dynkin-system $\mathcal{D}$, even when extendable, only allows for probabilistic statements
on $\mathcal{D}$ itself.
However, extendability guarantees that a ``nice'' embedding into
a larger system of measurable sets exists. More specifically,
extendability expressed in terms of credal sets
provides a well-known tool for the worst-case extension of a probability from
a pre-Dynkin-system to a larger algebra.

If a finitely additive probability on a pre-Dynkin-system is extendable, then we can obtain lower and upper probabilities of events which are not in the pre-Dynkin-system
but on a larger algebra. We follow the idea of natural extensions, \eg as described by \citep[page 136]{walley1991statistical}.
In particular, \citep[Theorem 3.3.4 (b)]{walley1991statistical} directly applies as long as a probability on a pre-Dynkin-system is extendable.
\begin{corollary}[Coherent Extension of Probability]
\label{corollary:Coherent Extension of Probability}
    Let $\mathcal{D}$ be a pre-Dynkin-system on $\Omega$.
    For a finitely additive probability measure $\mu$ on $\mathcal{D}$ we define the \emph{credal set}
    \begin{align*}
        \Mba(\mu, \mathcal{D}) \coloneqq \left \{ \nu \in \Pba \colon \nu(A) = \mu(A), \ \forall A \in \mathcal{D} \right \}.
    \end{align*}
    If $\mu$ on $\mathcal{D}$ is extendable to $2^\Omega$, then $ \forall A \in 2^\Omega$,
    \begin{align*}
        \underline{\mu}_{\mathcal{D}}(A) \coloneqq \inf_{\nu \in \Mba(\mu, \mathcal{D})}  \nu(A),
        \quad
        \overline{\mu}_{\mathcal{D}}(A) \coloneqq \sup_{\nu \in \Mba(\mu, \mathcal{D})}  \nu(A).
    \end{align*}
    define a coherent lower respectively upper probability on $2^\Omega$.
\end{corollary}
{
\color{darkgray}
\begin{example}
\label{runningexample:coherent extension}
The coherent extension of $\mu$ on $\mathcal{D}_4$ as defined in Example~\ref{runningexample:extendable probability on pre-Dynkin-system}
is $\underline{\mu}_{\mathcal{D}_4} = \ell$ and $\overline{\mu}_{\mathcal{D}_4} = u$
where, $\ell$ and $u$ are defined as in Example~\ref{runningexample:imprecise probability give pre-Dynkin-system} (\cf \citep[page 122]{walley1991statistical}).
Figure~\ref{fig:runningexample} illustrates the coherent extensions.
Even though coherent, $\underline{\mu}_{\mathcal{D}_4}$ is neither supermodular nor submodular:
\begin{align*}
    \underline{\mu}_{\mathcal{D}_4}(12) + \underline{\mu}_{\mathcal{D}_4}(13) = 0.5 + 0.2 > 0.5 + 0.0 = \underline{\mu}_{\mathcal{D}_4}(123) + \underline{\mu}_{\mathcal{D}_4}(1),\\
    \underline{\mu}_{\mathcal{D}_4}(1) + \underline{\mu}_{\mathcal{D}_4}(2) = 0.0 + 0.3 < 0.5 + 0.0 = \underline{\mu}_{\mathcal{D}_4}(12) + \underline{\mu}_{\mathcal{D}_4}(\emptyset).
\end{align*}
This implies that as well $\overline{\mu}_{\mathcal{D}_4}$ is neither supermodular nor submodular \citep[Proposition 2.3]{denneberg1994non}.
\end{example}
}
These lower and upper probabilities allow for at least two interpretations:
We can assume that a precise probability on a pre-Dynkin-systems $\mathcal{D} \subseteq 2^\Omega$
just reveals its values on $\mathcal{D}$, but is actually defined over $2^\Omega$. Then the
lower and upper probability constitute lower and upper bounds of the precise ``hidden probability'' on $2^\Omega$,
which is solely accessible on $\mathcal{D}$. On the other hand, we can even reject the existence of such precise ``hidden probability''.
Then lower and upper probability \emph{are} the inherently imprecise probability of an event in $2^\Omega$ but not in $\mathcal{D}$.\footnote{
As remarked by \citet[page 138]{walley1991statistical}, \citet{de2017theory} surprisingly only considered the first mentioned interpretation.}

The obtained lower and upper probabilities represent the imprecise interdependencies between all events of precise probabilities.
We illustrate this statement:
in the variety of updating methods in imprecise probability we pick the generalized Bayes' rule \citep[\S 6.4]{walley1991statistical} to exemplarily compute the conditional probability of two events
for the coherent extension of a probability from a pre-Dynkin-system. For $A,B \in \mathcal{D}$ such that $\mu(B) > 0$ the generalized Bayes' rule gives \citep[Theorem 6.4.2]{walley1991statistical}:
\begin{align*}
    \overline{\mu}_\mathcal{D}(A|B) &\coloneqq \sup_{\nu \in \Mba(\mu, \mathcal{D})} \frac{\nu(A \cap B)}{\nu(B)} = \frac{\sup_{\nu \in \Mba(\mu, \mathcal{D})} \nu(A \cap B)}{\mu(B)} = \frac{\overline{\mu}_\mathcal{D}(A \cap B)}{\mu(B)}, \qquad \forall A,B \in \mathcal{D}
\end{align*}
We can easily rearrange the above as
$\overline{\mu}_\mathcal{D}(A \cap B) = \overline{\mu}_\mathcal{D}(A|B) \mu(B)$.
In this case the imprecision of the probability of the intersected event is \emph{purely}
controlled by the conditional probability $\overline{\mu}(A|B)$ and \emph{not} by the marginal,
which is precise. So, the imprecision captured by the lower and upper probabilities
locates solely in the interdependency of the events.
We remark that Dempster's rule gives the same conditional probability here \citep{dempster1967upper}.

\subsection{Inner and Outer Extension Is More Pessimistic Than Coherent Extension}
We have presented two extension methods for probabilities defined on pre-Dynkin-systems. We relate the methods in the following. In the case of an extendable probability
we can guarantee the following inequalities to hold.
\begin{theorem}[Extension Theorem -- Finitely Additive Case]
\label{thm:Extension Theorem - Finitely Additive Case}
    Let $\mathcal{D}$ be a pre-Dynkin-system on $\Omega$ and
    $\mu$ a finitely additive probability
    on $\mathcal{D}$ which is extendable to $2^\Omega$. Then
    \begin{align*}
        \mu_*(A) \le \underline{\mu}_\mathcal{D}(A) \le  \overline{\mu}_\mathcal{D}(A) \le \mu^*(A), \qquad \forall A \in 2^\Omega.
    \end{align*}
\end{theorem}
\begin{proof}
    Since $\mathcal{D} \subseteq 2^\Omega$, we easily obtain
    \begin{align*}
        \mu_*(A) &= \sup \{ \mu(B) \colon A \supseteq B \in \mathcal{D}\}\\
        &= \sup \left \{ \underline{\mu}_\mathcal{D}(B) \colon A \supseteq B \in \mathcal{D} \right \}\\
        &\le \sup \left \{ \underline{\mu}_\mathcal{D}(B) \colon A \supseteq B \in 2^\Omega \right \}\\
        &= \underline{\mu}_\mathcal{D}(A),
    \end{align*}
    for all $A \in 2^\Omega$.
    The other inequalities follow by the conjugacy of inner and outer measure, and lower and upper
    coherent extension.
\end{proof}
In words, Theorem~\ref{thm:Extension Theorem - Finitely Additive Case} states that the
inner and outer extension is more ``pessimistic'' than the coherent extension. We use ``pessimistic''
in the sense of giving a looser bound for the probabilities assigned to elements
not in the pre-Dynkin-system $\mathcal{D}$ but in $2^\Omega$.
In Appendix~\ref{appendix:Dynkin-Systems and Countably Additive Probability} we demonstrate
an analogous result for countably additive probabilities on Dynkin-systems.

\section{The Credal Set and its Relation to Pre-Dynkin-System Structure}
\label{sec:the credal set and its relation to Dynkin-system structure}
In the earlier parts of the paper, we derived pre-Dynkin-systems as the system of precision for relatively general imprecise probabilities.
Then, we showed that, under extendability conditions, a precise probability on a pre-Dynkin-system gives rise
to a coherent imprecise probability on an encompassing algebra.
In other words, imprecise probabilities can be ``mapped'' to pre-Dynkin-systems and vice-versa.
We concretize these mappings in the following. This manifestation then reveals
structure in the interplay between the systems of precision, \ie pre-Dynkin-systems,
and coherent imprecise probabilities.
In particular, we argue that
the order structure of pre-Dynkin-systems can be mapped
to the space of finitely additive probabilities. This provides 
a (lattice) duality for coherent imprecise probabilities with precise
probabilities on pre-Dynkin-systems. More concretely, the duality
allows for the interpolation from imprecise probabilities which are precise on ``all'' events
to imprecise probabilities which are precise only on the empty set and the entire set.

In the following discussion, we assume, in addition to the technicalities presented
in Section~\ref{notation and technical details}, that a fixed finitely additive probability
on $2^\Omega$, which we call $\psi$, is given. The finitely additive probability $\psi$ with the algebra
$2^\Omega$ and the base set $\Omega$ constitute our ``base measure space'' analogous
to the choice of a base measure space in the theory
of coherent risk measures \citep{delbaen2002coherent}.
In comparison to the previous sections, we use $\psi$ instead of $\mu$ as ``reference measure'' to emphasize the difference that $\mu$ was defined on a relatively arbitrary pre-Dynkin-systems $\mathcal{D}$ on $\Omega$, while
$\psi$ is defined and fixed on the algebra $2^\Omega$ on $\Omega$.
    


\subsection{Credal Set Function Maps From Pre-Dynkin-Systems to Coherent Probabilities}
Equipped with a reference measure $\psi$ we define the credal set function. The name arises due to its close link to the credal set as defined in Corollary~\ref{corollary:Coherent Extension of Probability}.
\begin{definition}[Credal Set Function]
\label{def:Credal Set Function}
    Let $\Pba$ be the set of all finitely additive probabilities on $2^\Omega$. For a fixed finitely additive probability $\psi \in \Pba$ we call
    \begin{align*}
        \mba \colon 2^{2^\Omega} \rightarrow 2^{\Pba}, \quad
        \mba(\mathcal{A}) \coloneqq \{\nu \in \Pba \colon \nu(A) = \psi(A), \ \forall A \in \mathcal{A} \},
    \end{align*}
    the \emph{credal set function}.
\end{definition}
{
\color{darkgray}
\begin{example}
\label{runningexample:credal set function}
Let $\Omega_4 = \{ 1,2,3,4\}$ as in Example~\ref{runningexample:pre-Dynkin-system}. With abuse of notation we define the probability $\psi\colon 2^{\Omega_4} \rightarrow [0,1]$ via its corresponding point on the simplex $\psi \in \Delta$, $\psi_1 = 0.2, \psi_2 = 0.3, \psi_3 = 0.5, \psi_4 = 0$.
It follows, \eg $\mba(\{ 12,3\}) = \{ \nu \in \Pba \colon \nu_1 + \nu_2 = 0.5, \nu_3 = 0.5\}$. In the subsequent examples we implicitly assume the here defined $\psi$.
\end{example}
}
We stress that although not notated explicitly, the credal set function depends upon the choice of $\psi$.
For a fixed $\psi$ on $2^\Omega$, the credal set function $\mba$ maps a subset of the algebra
$2^\Omega$ to the set of all finitely additive probabilities which coincide with $\psi$ on this subset.
It should be noticed that by definition of $\psi$, $\mba(\mathcal{A}) \neq \emptyset$ for every non-empty
$\mathcal{A} \subseteq 2^\Omega$, because $\psi \in \mba(\mathcal{A})$ for every $\mathcal{A} \subseteq 2^\Omega$.

This defined mapping now simplifies our discussion about how pre-Dynkin-systems and imprecise probabilities correspond.
For instance, one can easily see that
the extreme case $\mathcal{A} = 2^\Omega$ corresponds to $\mba(\mathcal{A}) = \{ \psi\}$ and
$\mathcal{A} = \emptyset$ to $\mba(\mathcal{A}) = \Pba$.
More generally, we observe the following two properties of the credal set function.
\begin{proposition}[Credal Set Function is Invariant to Pre-Dynkin-Hull]
\label{prop:Credal Set Function is Invariant to Pre-Dynkin-Hull}
    Let $\mba$ be the credal set function.
    For any $\mathcal{A} \subseteq 2^\Omega$
    \begin{align*}
        \mba(\mathcal{A}) =  \mba(\PreDhull(\mathcal{A})).
    \end{align*}
\end{proposition}
\begin{proof}
    We need to show that
    \begin{align*}
        \{\nu \in \Pba\colon \nu(A) = \psi(A), \ A \in \mathcal{A}\} = \{\nu \in \Pba\colon \nu(A) = \psi(A), \ A \in \PreDhull(\mathcal{A})\}.
    \end{align*}
    The set inclusion of the right hand side in the left hand side is trivial.
    For the reverse direction, consider an element $\nu \in \Pba$ such that $\nu(A) = \psi(A)$ for $A \in \mathcal{A}$. Let
    \begin{align*}
        \mathcal{H} \coloneqq \{ A \in \PreDhull(\mathcal{A})\colon \nu(A) = \psi(A) \}.
    \end{align*}
    By Theorem~\ref{thm:lower and upper probablity define Dynkin-system} $\mathcal{H}$ is a pre-Dynkin-system.
    Since $\mathcal{A} \subseteq \mathcal{H}$ we know $\PreDhull(\mathcal{A}) \subseteq \mathcal{H}$.
    Hence, $\nu(A) = \mu(A)$ for $A \in \PreDhull(\mathcal{A})$.
    This gives the desired inclusion.
    We remark that for $\mathcal{A} = \emptyset$ the equality still holds, since $\PreDhull(\emptyset) = \{ \emptyset, \Omega\}$.
\end{proof}
{
\color{darkgray}
\begin{example}
\label{runningexample:invariance to pre-dynkin hull}
The credal set $\mba(\{ 12,3\}) = \{ \nu \in \Pba \colon \nu_1 + \nu_2 = 0.5, \nu_3 = 0.5\}$ given in Example~\ref{runningexample:credal set function} nicely illustrates Proposition~\ref{prop:Credal Set Function is Invariant to Pre-Dynkin-Hull}:
\begin{align*}
    \mba(\{ 12,3\}) = \{ \nu \in \Pba \colon \nu_1 + \nu_2 = 0.5, \nu_3 = 0.5, \nu_4 = 0\} = \mba(\{ \emptyset, 12, 3,4,34, 123,124, \Omega_4\}).
\end{align*}
\end{example}
}
\begin{proposition}[Credal Set Function Maps to Weak$^\star$-Closed Convex Sets]
\label{prop:Credal Set Function Maps to Closed Convex Sets}
    Let $\mba$ be the credal set function.
    For every non-empty $\mathcal{A} \subseteq 2^\Omega$, $\mba(\mathcal{A})$ is weak$^\star$-closed convex.
\end{proposition}
\begin{proof}
    The reference probability $\psi$ is by definition coherent. Hence, for all non-empty $\mathcal{A} \subseteq 2^\Omega$, the set $\mba(\mathcal{A})$ is the 
    set of all $\nu \in \Pba$ which dominate $\psi$. This set is, by Theorem 3.6.1 in \citep{walley1991statistical}, weak$^\star$-closed and convex.
\end{proof}
In words, the credal set on some set system coincides with
the credal set on its generated pre-Dynkin-system. And the credal set
of probabilities is always weak$^\star$-closed and convex.
Proposition~\ref{prop:Credal Set Function is Invariant to Pre-Dynkin-Hull} allows us to work with credal sets of arbitrary
set systems instead of the entire pre-Dynkin-system.
Thus, it resembles the well-known $\pi$-$\lambda$-Theorem, which is
fundamental to classical probability theory \citep[Lemma A.1.3]{williams1991probability}.
On the other hand, this result justifies our 
focus on pre-Dynkin-systems instead of arbitrary set systems.
We do not lose generality when considering pre-Dynkin-systems instead of 
non-structured sets of sets.

Proposition~\ref{prop:Credal Set Function Maps to Closed Convex Sets} guarantees that the images
of the credal set function behave ``nicely''. Specifically, these weak$^\star$-closed convex sets correspond to coherent previsions, \ie generalizations of coherent probabilities
as already stated by \citet[Theorem 3.6.1]{walley1991statistical}. We elaborate this observation in Section~\ref{sec:interpolation from algebra to trivial dynkin-system}.
In conclusion, credal set functions map pre-Dynkin-systems to coherent probabilities. What about the reverse mapping?

\subsection{The Dual Credal Set Function}
The following is a natural definition of a dual credal set function. We justify this name by Proposition~\ref{prop:galois connection by (dual) credal set function} below.
\begin{definition}[Dual Credal Set Function]
\label{def:Dual Credal Set Function}
    Let $\Pba$ be the set of all finitely additive probabilities on $2^\Omega$. Fix a finitely additive probability $\psi$ on $2^\Omega$. We call
    \begin{align*}
        \mdba \colon 2^{\Pba} \rightarrow 2^{2^\Omega}, \quad
        \mdba(\Q) \coloneqq \{A \in 2^\Omega \colon \nu(A) = \psi(A), \forall \nu \in \Q \}
    \end{align*}
    the \emph{dual credal set function}.
\end{definition} The dual credal set function also depends upon $\psi$, but we do not notate this explicitly. The dual credal set function maps an arbitrary set of finitely additive probability measures
on $2^\Omega$ to the (largest) set of events on which all contained probabilities coincide.
We remark that each set of finitely additive probability measures can be linked to an imprecise probability.

We suggestively called the antagonist to the credal set function the ``dual credal set function''.
The duality appearing here is a well-known fundamental relationship
between partially ordered sets: a Galois connection. A \emph{Galois connection} is a pair of mappings $f\colon X \rightarrow Y$ and $f^\circ \colon Y\rightarrow X$ 
on partially ordered sets $(Y,\le)$ and $(X, \le)$ which preserves order structure (\cf Corollary~\ref{corollary:rules for dual credal set function}).
More formally, $f$ and $f^\circ$ are a Galois connection if and only if for all $x \in X, y \in Y$, $x \le f(y) \Leftrightarrow y \le f^\circ(x)$ \citep[\S V.8]{birkhoff1940lattice}.
Galois connections, even though they do \emph{not} form an order isomorphism, induce a
lattice duality. We exploit this lattice duality to provide an order-theoretic interpolation from
no compatibility at all to full compatability. For this purpose, we first establish the Galois connection.
\begin{proposition}[Galois Connection by (Dual) Credal Set Function]
\label{prop:galois connection by (dual) credal set function}
    The credal set function $\mba$ and the dual credal set function $\mdba$ form a Galois connection.
\end{proposition}
\begin{proof}
    $\mba$ and $\mdba$
    form a Galois connection if and only if $\mathcal{A} \subseteq \mdba(\Q) \Leftrightarrow \Q \subseteq \mba(\mathcal{A})$ \citep[\S V.8]{birkhoff1940lattice}.
    First, we show the left to right implication. We assume $\mathcal{A} \subseteq \mdba(\Q)$, i.e.
    every $\nu \in \Q$ coincides with $\psi$ on $\mathcal{A}$. Hence,
    \begin{align*}
        \nu \in \Q &\Rightarrow  \nu(A) = \psi(A), \forall A \in \mathcal{A}\\
        &\Rightarrow \nu \in \{\nu' \in \Pba\colon \nu'(A) = \psi(A), \forall A \in \mathcal{A} \} = \mba(\mathcal{A}).
    \end{align*}
    In case of the right to left implication we suppose $\Q \subseteq \mba(\mathcal{A})$. Thus,
    \begin{align*}
        A \in \mathcal{A} &\Rightarrow \nu(A) = \psi(A), \forall \nu \in \Q\\
        &\Rightarrow A \in \{A' \in 2^\Omega\colon \nu(A') = \psi(A'), \forall \nu \in \Q \} = \mdba(\Q).
    \end{align*}
\end{proof}
The mappings involved in the Galois connection are antitone, \ie they reverse the order structure from domain to codomain.
Their pairwise application
is extensive, \ie the image of an object contains the object. In summary, the following rules of calculation hold:
\begin{corollary}[Rules for (Dual) Credal Set Function]
\label{corollary:rules for dual credal set function}
    Let $\mba$ be the credal set function
    and $\mdba$ be the dual credal set function.
    For arbitrary $\mathcal{A}_1, \mathcal{A}_2, \mathcal{A} \subseteq 2^\Omega$ and $\Q_1, \Q_2, \Q \subseteq \Pba$,
    \begin{align*}
        \mathcal{A}_1 \subseteq \mathcal{A}_2 &\Rightarrow \mba(\mathcal{A}_2) \subseteq \mba(\mathcal{A}_1), & \Q_1 \subseteq \Q_2 &\Rightarrow \mdba(\Q_2) \subseteq \mdba(\Q_1),  &\text{(antitone)}\\
        \Q &\subseteq \mba(\mdba(\Q)), & \mathcal{A} &\subseteq \mdba(\mba(\mathcal{A})), &\text{(extensive)}\\
        \mba(\mathcal{A}) &= \mba(\mdba(\mba( \mathcal{A} ))), & \mdba(\Q) &= \mdba(\mba(\mdba( \Q ))),  &\text{(pseudo-inverse).}
    \end{align*}
\end{corollary}
\begin{proof}
    \citep[\S V.7 and V.8]{birkhoff1940lattice}
\end{proof}
Proposition~\ref{prop:galois connection by (dual) credal set function} provides a tool to further investigate the dual credal set function.
The reader might have noticed the similarity of the dual credal set function
and the main question of Section~\ref{From Porbability measures to Dynkin systesms}: given a lower and upper probability, on which
set systems do both coincide? In fact, we obtain an analogous result to Theorem~\ref{thm:lower and upper probablity define Dynkin-system}, again an imprecise probability is mapped to the set of events on which it is precise.
\begin{proposition}[Dual Credal Set Function Maps to Pre-Dynkin-Systems]
\label{prop:Dual Credal Set Function Maps to Pre-Dynkin-Systems}
    Let $\mdba$ be the dual credal set function.
    For all non-empty $\Q \subseteq \Pba$, $\mdba(\Q)$ is a pre-Dynkin-system.
\end{proposition}
\begin{proof}
    We show the statement by establishing the equality $\PreDhull(\mdba(\Q)) = \mdba(\Q)$.
    Trivially, $\PreDhull(\mdba(\Q)) \supseteq \mdba(\Q)$. Furthermore,
    \begin{align*}
        \PreDhull(\mdba(\Q)) &\overset{C.\ref{corollary:rules for dual credal set function}}{\subseteq} \mdba( \mba(\PreDhull(\mdba(\Q))))\\
        &\overset{P.\ref{prop:Credal Set Function is Invariant to Pre-Dynkin-Hull}}{=} \mdba( \mba(\mdba(\Q))) \overset{C.\ref{corollary:rules for dual credal set function}}{=} \mdba(\Q).
    \end{align*}
\end{proof}
\begin{proposition}[Dual Credal Set Function is Invariant to Weak$^\star$-Closed Convex Hull]
\label{prop:Dual Credal Set Function is Invariant to Closed Convex Hull}
    Let $\mdba$ be the dual credal set function.
    For any $\Q \subseteq \Pba$, $\mdba(\Q) = \mdba(\co \Q)$.\footnote{As introduced in Section~\ref{notation and technical details}, $\co$ denotes the convex, weak$^\star$ closure of a set in $\Pba \subseteq \SpaceOfLinearFunctionalsOnGambles$.}
\end{proposition}
\begin{proof}
    By Corollary~\ref{corollary:rules for dual credal set function} we know $\mdba(\Q) = \mdba(\mba(\mdba(\Q)))$.
    Furthermore, $\Q \subseteq \mba(\mdba(\Q))$. Via, Proposition~\ref{prop:Credal Set Function Maps to Closed Convex Sets} we obtain $\co \Q \subseteq \mba(\mdba(\Q))$.
    Hence, the result follows.
\end{proof}
{
\color{darkgray}
\begin{example}
\label{runningexample:dual credal set function maps to pre-dynkin-systems}
Let $\Q = \{ \nu\}$ where we identify the probability $\nu$ on $2^{\Omega_4}$ with an element $\nu \in \Delta_4$. For instance, $\nu_1 = 0, \nu_2 = 0.5, \nu_3 = 0, \nu_4 = 0.5$. It is easy to see that
$\mdba(\Q) = \PreDhull(\{12, 3\})$ as given in Example~\ref{runningexample:invariance to pre-dynkin hull}.
\end{example}
}
In summary, credal set functions map pre-Dynkin-systems to weak$^\star$-closed convex credal sets. Dual credal set
functions map (weak$^\star$-closed convex) credal sets to pre-Dynkin-systems. In addition, the two functions form a Galois connection.
In fact, every Galois connection defines closure operators, \ie extensive, monotone and idempotent maps \citep[Definition 4.5.a]{schechter1997handbook}.
The closure operators are defined as the sequential application of the credal set function and the dual credal set function to subsets of $\Pba$ or $2^\Omega$. In symbols: $
    \mathcal{A} \mapsto \mdba(\mba(\mathcal{A})), 
    \Q \mapsto \mba(\mdba(\Q))$
In particular, these closure operators define bipolar-closed sets.

\subsection{Bipolar-Closed Sets}
\label{sec:bipolar-closed sets}
\begin{figure}
    \centering
    \includegraphics[width=\textwidth]{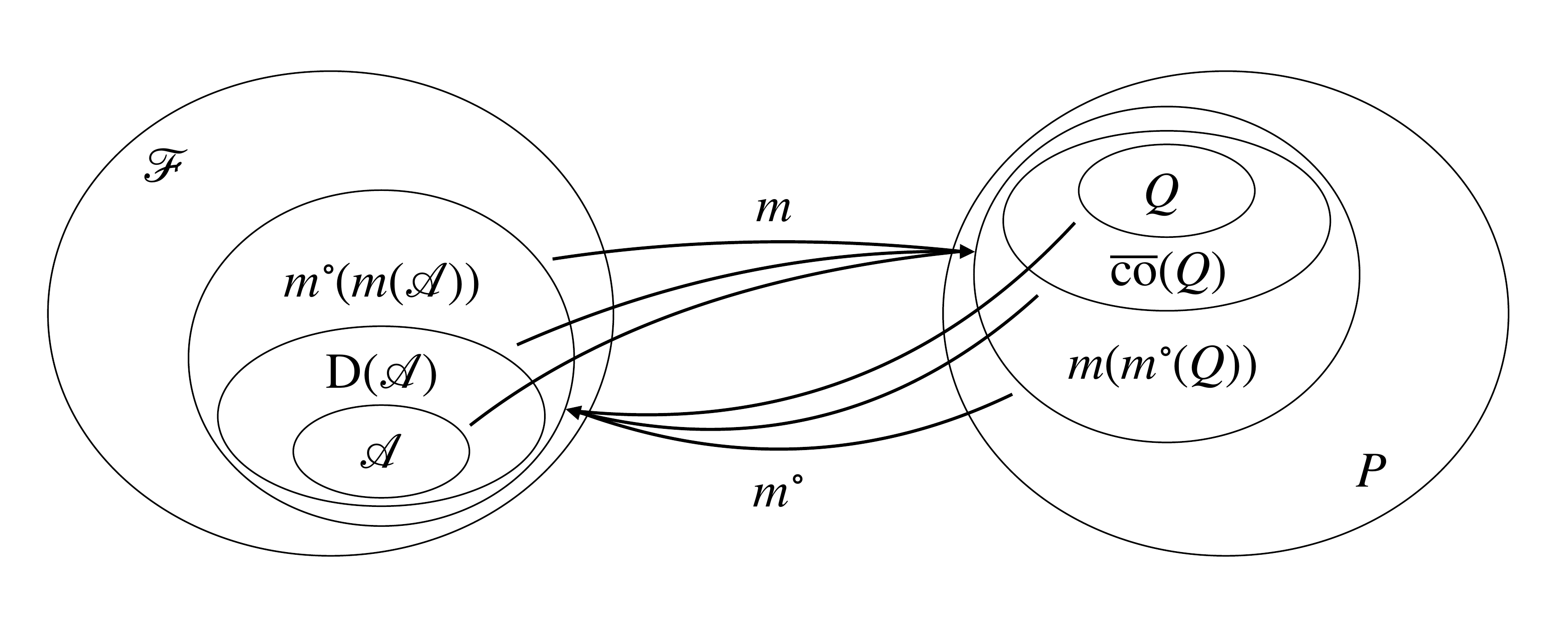}
    \caption{Galois connection between the lattice of pre-Dynkin-systems and the set of credal sets. In the illustrated case, we have $\mdba(\Q) = \mdba(\mba(\mdba(\mathcal{A})))$ respectively
    $\mba(\mathcal{A}) = \mba(\mdba(\mba(\Q)))$.
    The set containment on both sides follows from Proposition~\ref{prop:Credal Set Function is Invariant to Pre-Dynkin-Hull}, Corollary~\ref{corollary:rules for dual credal set function} and Proposition~\ref{prop:Dual Credal Set Function is Invariant to Closed Convex Hull}.}
    \label{fig:galois connection}
\end{figure}

Bipolar-closed sets are sets $\mathcal{A} \subseteq 2^\Omega$ such that 
$\mathcal{A} = \mdba(\mba(\mathcal{A}))$,
respectively $\Q \subseteq \Pba$ such that 
 $\Q = \mba(\mdba(\Q))$.
Most importantly, the bipolar-closed sets
form two antitone isomorphic lattices
ordered by set inclusion \citep[Theorem V.8.20]{birkhoff1940lattice}.
This relationship gives us a lattice duality between set systems and credal sets of probabilities.
See Figure~\ref{fig:galois connection} for an illustration of bipolar-closed sets and the Galois connection.
{
\color{darkgray}
\begin{example}
\label{runningexample:bipolar closed sets}
The pre-Dynkin-system $\PreDhull(\{ 12,3\})$ already discussed in Example~\ref{runningexample:invariance to pre-dynkin hull} and Example~\ref{runningexample:dual credal set function maps to pre-dynkin-systems} is a bipolar-closed set. In contrast, the set $\{ 12,3\}$ cannot be a bipolar-closed set, as it is not a pre-Dynkin-system.
\end{example}
}
More precisely, bipolar-closed sets in the set of finitely additive probability distributions
are weak$^\star$-closed convex (Proposition~\ref{prop:Credal Set Function Maps to Closed Convex Sets}).
These map to
bipolar-closed subsets of $2^\Omega$, which are pre-Dynkin-systems (Proposition~\ref{prop:Dual Credal Set Function Maps to Pre-Dynkin-Systems}).
All of the stated properties of bipolar-closed sets are necessary. But are they sufficient?

\subsubsection{Sufficient Conditions for Bipolar-Closed Sets}
In the search for sufficient conditions for bipolar-closed sets we 
focus on bipolar-closed subsets of $2^\Omega$. 
Bipolar-closed subsets of ${\Pba}$ require further investigation. 
It is already difficult to characterize sufficient conditions for bipolar-closed subsets of $2^\Omega$.
\begin{corollary}
\label{corollary:dynkin-system is inside of bipolar closure}
    Let $\mba$ be the credal set function and $\mdba$ be the dual credal set function.
    For an arbitrary subset $\mathcal{A} \subseteq 2^\Omega$ we have
    \begin{align*}
        \PreDhull(\mathcal{A}) \subseteq \mdba (\mba(\mathcal{A})).
    \end{align*}
\end{corollary}
\begin{proof}
    By Proposition~\ref{prop:Credal Set Function is Invariant to Pre-Dynkin-Hull},
    $\mdba ( \mba(\mathcal{A})) =  \mdba ( \mba(\PreDhull(\mathcal{A})))$. By Corollary~\ref{corollary:rules for dual credal set function}, the statement follows.
\end{proof}
This corollary gives rise to the follow-up question: under which circumstances does $\PreDhull(\mathcal{A}) = \mdba (\mba(\mathcal{A}))$?
As the following theorem demonstrates, this question is closely connected to the sets
of measure zero of the base probability $\psi$ and its problems (\cf\citep{rota2001twelve}).
\begin{proposition}[``Closedness'' under Measure Zero Sets]
\label{proposition:closedness under measure zero sets}
    Let $\mba$ be the credal set function and $\mdba$ be the dual credal set function. Let $\mathcal{A} \subseteq 2^\Omega$.
    For any $A \in \mdba (\mba(\mathcal{A}))$, if $\psi(A) = 0$, then all $B, C \in 2^\Omega$ such that
    $B \subseteq A $ and $C \supseteq A^c$
    are in $\mdba (\mba(\mathcal{A}))$.
\end{proposition}
\begin{proof}
    Let  $A \in \mdba (\mba(\mathcal{A}))$ for arbitrary $\mathcal{A} \subseteq 2^\Omega$
    with $\psi(A) = 0$.
    Consider two sets $B,C \in 2^\Omega$ such that $B \subseteq A $ respectively
    $C \supseteq A^c$. Then $\psi(B) = 0$ and $\psi(C) = 1$.
    Furthermore, for any $\nu \in \mba(\mdba (\mba(\mathcal{A})))$, we have
    $\nu(B) = 0$ (respectively $\nu(C) = 1$). Since $ \mba(\mdba (\mba(\mathcal{A})))= \mba(\mathcal{A})$,
    we obtain $B, C \in \mdba(\mba(\mathcal{A}))$.
\end{proof}
A pre-Dynkin-system $\PreDhull(\mathcal{A})$ can only coincide with $\mdba (\mba(\mathcal{A}))$
if subsets of measure zero sets are included.
Thus, Proposition~\ref{proposition:closedness under measure zero sets} provides a further necessary
condition for bipolar-closed subsets of $2^\Omega$.
Yet, it turns out that the sets of measure zero as well can give a sufficient condition
for bipolar-closed sets at least in a finite setting.
\begin{theorem}[Pre-Dynkin-Hull is Bipolar-Closure in Finite, Discrete Setting]
\label{thm:Pre-Dynkin-Hull is Bipolar-Closure in Finite, Discrete Setting}
    Let $\Omega = [n]$. Fix a finitely additive probability
    $\psi$ on $2^\Omega$ such that $\psi(A)> 0$ for every $A \in 2^\Omega \setminus \{ \emptyset \}$.
    Let $\mathcal{D} \subseteq 2^\Omega$ be a pre-Dynkin-system. Then,
    \begin{align*}
        \mathcal{D} = \mdba (\mba(\mathcal{D})).
    \end{align*}
\end{theorem}
\begin{proof}
    If $\mathcal{D} = 2^{[n]}$ the result follows directly. Thus, from now on we
    assume $\mathcal{D} \subsetneq 2^{[n]}$.
    The set inclusion $\mathcal{D} \subseteq \mdba (\mba(\mathcal{D}))$
    is given by Corollary~\ref{corollary:dynkin-system is inside of bipolar closure}.
    We show $\mathcal{D} \supseteq \mdba (\mba(\mathcal{D}))$
    by proving that for every $B \in 2^{[n]} \setminus \mathcal{D}$ there exists $\nu \in \mba(\mathcal{D})$
    such that $\nu(B) \neq \psi(B)$.
    
    Let us consider an arbitrary $B \in 2^{[n]} \setminus \mathcal{D}$.
    Without loss of generality (Lemma~\ref{lemma:decomposition into dynkin set and atom})
    we can decompose
    \begin{align*}
        B = B_\mathcal{D} \cup A,
    \end{align*}
    where $B_{\mathcal{D}} \in \mathcal{D}$ and $A \notin \mathcal{D}$
    is a weak atom with respect to $\mathcal{D}$ (Definition~\ref{def:weak atom wrt pre-Dynkin-system}).
    Thus, we can leverage Lemma~\ref{lemma:Weak Atoms Are Not in the Bipolar-Closed Sets - Finite, Discrete Setting}:
    there exists $\nu \in \mba(\mathcal{D})$ such that $\nu(A) \neq \psi(A)$.
    Thus,
    \begin{align*}
        \nu(B) = \nu(B_\mathcal{D}) + \nu(A) = \psi(B_\mathcal{D}) + \nu(A) \neq \psi(B_\mathcal{D}) + \mu(A).
    \end{align*}
    It follows that $B \notin \mdba (\mba(\mathcal{D}))$, concluding the proof.
\end{proof}
Whether this theorem can be extended to more general sets $\Omega$ is an open question. Seemingly, proofs along the line of Theorem~\ref{thm:Pre-Dynkin-Hull is Bipolar-Closure in Finite, Discrete Setting} are doomed to fail, since one cannot argue via probabilities on atoms of $\Omega$.

\subsection{Interpolation From Algebra to Trivial Pre-Dynkin-System}
\label{sec:interpolation from algebra to trivial dynkin-system}
In probability theory there is a choice to be made regarding which events should get assigned probabilities \citep[page 52]{kolmogorov1929}. 
This significant choice has (mathematically) been standardized to form a ($\sigma$-)algebra (\cf standard probability space).
But, already Kolmogorov, the ``father'' of modern probability theory, emphasized that this choice is not universal, but should depend on the problem at hand. 
More recently, \citet{khrennikov2016bell} argued that a more appropriate probabilistic modeling should
appeal to weaker domains for probabilities to, for instance, represent physical observations such as quantum phenomena.

In particular, it cannot always be taken for granted
that all events are compatible with all others, as implied by a ($\sigma$-)algebra (\cf Section~\ref{sec:compatibility}).
For instance, von Mises' axiomatization 
of probability inherently reflects potential incompatible events in terms of a pre-Dynkin-system \citep{mises1964mathematical, schurz2008finitistic}.
In other words, there is a choice to be made about the system of precision.
Which sets should be compatible to each other, which should not? How do the choices of the systems of precision
relate to each other?

We neglect, without loss of generality, arbitrary systems of precision
and focus on pre-Dynkin-systems (\cf Proposition~\ref{prop:Credal Set Function is Invariant to Pre-Dynkin-Hull}).
The range of choices is captured by
the system of pre-Dynkin-systems.
\begin{proposition}[Set of Pre-Dynkin-Systems is a Lattice]
\label{prop:Set of Pre-Dynkin-Systems Is a Lattice}
    The set~$\mathfrak{D} \coloneqq \{ \mathcal{D} \subseteq 2^\Omega\colon \mathcal{D} = \PreDhull(\mathcal{D})\}$ is ordered by set inclusion. Furthermore, $(\mathfrak{D}, \subseteq)$ is a lattice with
\begin{align*}
    \bigvee_{i = 1}^n \mathcal{D}_i &= \PreDhull\left(\bigcup_{i  = 1}^n \mathcal{D}_i\right)\\
    \bigwedge_{i = 1}^n \mathcal{D}_i &= \bigcap_{i  = 1}^n \mathcal{D}_i.
\end{align*}
\end{proposition}
\begin{proof}
    On the one hand, it is easy to show that the intersection of pre-Dynkin-systems forms a pre-Dynkin-system again. On other hand, the smallest pre-Dynkin-system which contains a finite set of pre-Dynkin-systems is by definition the pre-Dynkin-system generated by the union over all elements in this finite set of pre-Dynkin-systems.
\end{proof}
{
\color{darkgray}
\begin{example}
\label{runningexample:lattice of pre-dynkin-systems}
Let $\Omega_4$ be as defined in Example~\ref{runningexample:credal set function}. The minimal element in $\mathfrak{D}$ then is $\{ \emptyset, \Omega_4\}$. The maximal element is $2^{\Omega_4}$. For the sake of brevity, we omit all further elements in $\mathfrak{D}$ and remain with the observation that $\mathcal{D}_4$ of Example~\ref{runningexample:pre-Dynkin-system} and $\PreDhull(\{ 12,3\})$ are elements of~$\mathfrak{D}$.
\end{example}
}
The lattice $\mathfrak{D}$ spans a range of choices from $\mathcal{D} = 2^\Omega$, \ie 
complete compatibility and only a single probability distribution
in its credal set, namely $\mba(2^\Omega) = \{ \psi\}$,
to $\mathcal{D} = \{ \emptyset, \Omega\}$, \ie no compatibility and
the entire space of probability distributions constitute its credal set
$\mba(\{ \emptyset, \Omega\}) = \Pba$.
How ``close'' $\mathcal{D}$ is to the algebra $2^\Omega$ determines how ``classical'' the credal set behaves.
In other words, $(\mathfrak{D}, \subseteq)$ parametrizes a family of credal sets. Thus, it parametrizes
coherent probabilities. The knob of compatibility can be turned from
trivially nothing ($\{ \emptyset, \Omega \}$), to everything ($2^\Omega$).
How does the ``amount of compatibility'' of the pre-Dynkin-system map to the
credal sets? Or, \eg given two pre-Dynkin-systems on which a probability is
defined, what is the credal set of the union of these systems?
\begin{proposition}[Lattice of Dynkin-Systems and Credal Sets]
    The credal set function $\mba$ (Definition~\ref{def:Credal Set Function}) together with the lattice $(\mathfrak{D}, \subseteq)$ provides
    a parametrized family of credal sets for which hold ($\forall \mathcal{D}_1, \mathcal{D}_2 \in \mathfrak{D}$):
    \begin{align*}
        \mba(\mathcal{D}_1 \vee \mathcal{D}_2) &= \mba(\mathcal{D}_1 \cup \mathcal{D}_2) = \mba(\mathcal{D}_1) \cap \mba(\mathcal{D}_2)\\
        \mba(\mathcal{D}_1 \wedge \mathcal{D}_2) &=  \mba(\mathcal{D}_1 \cap \mathcal{D}_2) \supseteq \mba(\mathcal{D}_1) \cup \mba(\mathcal{D}_2).
    \end{align*}
\end{proposition}
\begin{proof}
    Concerning the first equality, we observe that for arbitrary $\mathcal{D}_1, \mathcal{D}_2 \in \mathfrak{D}$
    \begin{align*}
        \mba(\mathcal{D}_1 \vee \mathcal{D}_2) = \mba(\PreDhull(\mathcal{D}_1 \vee \mathcal{D}_2)) = \mba(\mathcal{D}_1 \cup \mathcal{D}_2),
    \end{align*}
    by Proposition~\ref{prop:Credal Set Function is Invariant to Pre-Dynkin-Hull}. Consequently,
    \begin{align*}
         \mba(\mathcal{D}_1 \cup \mathcal{D}_2) &= \{\nu \in \Pba\colon \nu(D) = \psi(D), \ \forall D \in \mathcal{D}_1 \cup \mathcal{D}_2\}\\
         &= \{\nu \in \Pba\colon \nu(D) = \psi(D), \ \forall D \in \mathcal{D}_1\} \cap  \{\nu \in \Pba\colon \nu(D) = \psi(D), \ \forall D \in \mathcal{D}_2\} \\
         &= \mba(\mathcal{D}_1) \cap \mba(\mathcal{D}_2).
    \end{align*}
    The second line follows by the definition of infimum on the lattice of pre-Dynkin-systems and simple set containment: $\mba(\mathcal{D}_1) \subseteq \mba(\mathcal{D}_1 \cap \mathcal{D}_2)$ and $\mba(\mathcal{D}_2) \subseteq \mba(\mathcal{D}_1 \cap \mathcal{D}_2)$.
\end{proof}
Unfortunately, the mentioned interpolation is slightly improper.
It turns out that there are pre-Dynkin-systems $\mathcal{D}_1 \neq \mathcal{D}_2$
such that $\mba(\mathcal{D}_1) = \mba(\mathcal{D}_2)$.
{
\color{darkgray}
\begin{example}[Non-Injectivity of Credal Set Function]
\label{example:Non-Injectivity of Credal Set Function}
    Let $\Omega_3 = \{1,2,3\}$ and
    $\psi(\{1\}) = 1, \psi(\{2\}) = 0, \psi(\{3\}) = 0$ constitute
    a base probability space.
    Then, obviously $\mba(\PreDhull(\{ 1\})) = \mba(2^\Omega)$,
    but $\PreDhull(\{ 1\}) \neq 2^\Omega$.
\end{example}
}
The reason for this collision of credal sets is that not every pre-Dynkin-system $\mathcal{D} \subseteq 2^\Omega$
is a bipolar-closed set.
\begin{proposition}[Credal Set Function is Injective on Bipolar-Closed Sets]
    Let $\mba$ be the credal set and $\mdba$ be the dual credal set function.
    Let $\mathcal{D}_1, \mathcal{D}_2 \in \mathfrak{D}$ be pre-Dynkin-systems, which are bipolar-closed. If $\mathcal{D}_1 \neq \mathcal{D}_2$, then $\mba(\mathcal{D}_1) \neq \mba(\mathcal{D}_2)$.
\end{proposition}
\begin{proof}
    We prove the claim by contraposition:
    \begin{align*}
        \mba(\mathcal{D}_1) = \mba(\mathcal{D}_2) \Rightarrow \mdba(\mba(\mathcal{D}_1)) = \mdba(\mba(\mathcal{D}_2)) \Rightarrow \mathcal{D}_1 = \mathcal{D}_2.
    \end{align*}
\end{proof}
Hence, it is reasonable to focus on the set of bipolar-closed sets contained in $2^\Omega$. We define the \emph{set of interpolating pre-Dynkin-systems} $\mathfrak{C} \coloneqq \{ \mathcal{C} \subseteq 2^\Omega \colon \mathcal{C} = \mdba(\mba(\mathcal{C}))\}$. We know that $\mathfrak{C} \subseteq \mathfrak{D}$ (Proposition~\ref{prop:Dual Credal Set Function Maps to Pre-Dynkin-Systems}) and $\mathfrak{C}$ is even a lattice contained in $\mathfrak{D}$.
\begin{proposition}[Lattice of Bipolar-Closed Sets and Credal Sets]
    Let $\mba$ be the credal set function and $\mdba$ be the dual credal set function.
    The set of interpolating pre-Dynkin-systems $\mathfrak{C}$
    equipped with the $\subseteq$-ordering forms a lattice:
    \begin{align*}
        \bigvee_{i = 1}^n \mathcal{C}_i &= \mdba \left( \mba \left(\bigcup_{i  = 1}^n \mathcal{C}_i\right) \right)\\
        \bigwedge_{i = 1}^n \mathcal{C}_i &= \bigcap_{i  = 1}^n \mathcal{C}_i.
    \end{align*}
    In particular, this lattice $\mathfrak{C}$ is antitone isomorphic to the lattice of bipolar-closed sets in $2^{\Pba}$. It holds\footnote{We denote the composition of functions with $\circ$.}
    \begin{align*}
        \mba(\mathcal{C}_1 \vee \mathcal{C}_2) &= \mba(\mathcal{C}_1 \cup \mathcal{C}_2) = \mba(\mathcal{C}_1) \cap \mba(\mathcal{C}_2)\\
        \mba(\mathcal{C}_1 \wedge \mathcal{C}_2) &=  \mba(\mathcal{C}_1 \cap \mathcal{C}_2) = (\mba \circ \mdba) \left( \mba(\mathcal{C}_1) \cup \mba(\mathcal{C}_2)\right).
    \end{align*}
    
\end{proposition}
\begin{proof}
    By Theorem V.8.20 in \citep{birkhoff1940lattice} $\mathfrak{C}$ is a lattice and $\mba$ an antitone lattice isomorphism on $\mathfrak{C}$. The equations hold by simple manipulations
    \begin{align*}
        \mba(\mathcal{C}_1 \vee \mathcal{C}_2) = \mba(\mdba(\mba(\mathcal{C}_1 \cup \mathcal{C}_2))) = \mba(\mathcal{C}_1 \cup \mathcal{C}_2) = \mba(\mathcal{C}_1) \cap \mba(\mathcal{C}_2),
    \end{align*}
    and
    \begin{align*}
        \mba(\mathcal{C}_1 \wedge \mathcal{C}_2) = \mba(\mdba(\mba(\mathcal{C}_1 \cap \mathcal{C}_2)))
        = (\mba \circ \mdba) \left( \mba(\mathcal{C}_1) \cup \mba(\mathcal{C}_2)\right).
    \end{align*}
\end{proof}
Note that $\mathfrak{C}$ is not generally a sublattice of $\mathfrak{D}$. It is a lattice contained in the lattice $\mathfrak{D}$, but the closure operator for the supremum is distinct.
Interestingly, the order structure which both lattices, $\mathfrak{D}$ and $\mathfrak{C}$, induce on the set of credal sets via $\mba$ is identical. Every set $\mba(\mathcal{A})$ for arbitrary $\mathcal{A} \subseteq 2^\Omega$ is bipolar-closed (Corollary~\ref{corollary:rules for dual credal set function}). Hence, the lattice of bipolar-closed sets in $2^{ \Pba}$ is the domain of $\mba$ for elements in  $\mathfrak{D}$ and $\mathfrak{C}$.
In other words, the lattices $\mathfrak{D}$ and $\mathfrak{C}$
provide one and the same parametrized family of credal sets, thus one and the same parametrized family of
imprecise probabilities.

In comparison to other parametrized families of imprecise probability, such as 
distortion risk measures \citep{wirch2001distortion}, which heavily rely on convex analysis,
the duality used here is structurally weaker. Lattice isomorphisms
give a glimpse of structure to the involved dual spaces. Convex dualities
as exploited in \citep{frohlich2022risk} are far more informative, but apparently not able 
to handle the structural knob which we presented in this work: the set of sets which get assigned
precise probabilities.
Nevertheless, a natural question arises from this lattice duality: how does this lattice duality
relate to a convex duality?
We leave this question open to further research. A first attempt to an answer is discussed in
Appendix~\ref{sec:Credal Sets of Pre-Dynkin-system probabilities - Credal Sets of Distorted Probabilities}, where
we link the parametrized family of distorted probabilities to the pre-Dynkin-system
family of imprecise probabilities.

\section{A More General Perspective -- The Set of Gambles With Precise Expectation}
\label{sec:A More General Perspective - The Set of Gambles with Precise Expectation}
To this point, we have exclusively focused on probabilities and set systems of events to which we assign probabilities. In fact, there is a more general story to be told. In the literature on imprecise probability focus often lies on expectation-type functionals instead of probabilities and on sets of gambles (bounded functions from the base set $\Omega$ to the real numbers) instead of sets of events. One can easily see that the latter is more general and can recover the former. Indicator functions of events are gambles. An expectation-type functional evaluated on an indicator gamble of an event corresponds to a generalized probability of the event.
The converse direction, \ie recovering a unique expectation-type functional from an imprecise probability, however, is not always possible \citep[\S 2.7.3]{walley1991statistical}. In the following, we reiterate several questions which we asked in the preceding sections for probabilities and set systems.

\subsection{Partial Expectations Generalize Finitely Additive Probabilities on (Pre-)Dynkin-Systems}
We propose the following definition of partial expectation and show afterwards that it is a natural generalization of finitely additive probabilities defined on (pre-)Dynkin-systems.
\begin{definition}[Partial Expectation]
\label{def:partial expectation}
    Let $\{ L_i \}_{i \in I}$ be a non-empty family of linear subspaces of $\SpaceOfGambles$.
    We call $\pE \colon \bigcup_{i \in I} L_i \rightarrow \Reals$ a \emph{partial expectation} if and only if all of the following conditions are fulfilled:
    \begin{enumerate}[(a)]
        \item for any $i \in I$ and for all $f, g \in L_i$, then $ \pE(f + g) = \pE( f) + \pE(g)$, \quad (Partial Linearity),
        \item for any $i \in I$ and any $f \in L_i$, then $\pE(f) \ge \inf f$,\quad (Coherence).
    \end{enumerate}
\end{definition}
We remark that for this definition we leveraged the requirements for a linear prevision (Definition~\ref{def:coherent prevision}) on a linear space given in \citep[Theorem 2.8.4]{walley1991statistical}. In other words, a partial expectation is a functional which is defined on a union of linear subspaces and behaves like a ``classical'' (finitely additive) expectation on each of the subspaces, but not necessarily on all simultaneously. It is a linear prevision when restricted to one of the subspaces $L_i$ (\cf Definition~\ref{def:coherent prevision}).

There is a one-to-one correspondence of linear previsions and coherent additive probabilities.
For every coherent additive probability $\nu$ defined on an algebra $\mathcal{A}$, there is a unique linear prevision, which we equivalently denote $\nu$,
the set of all $\mathcal{A}$-measurable gambles, which agrees with the probability $\nu$ on the indicator gambles of the sets in $\mathcal{A}$ \citep[Theorem 3.2.2]{walley1991statistical}. The following Proposition exploits this correspondence. A finitely additive probability defined on a pre-Dynkin-system relates one-to-one to a partial expectation which is defined on the set of linear spaces induced by the simple gambles on the blocks of the pre-Dynkin-system. To this end, we introduce the following two notations: let $\mathcal{F} \subseteq 2^\Omega$ be an algebra. Then $\operatorname{S}(\Omega, \mathcal{F}) \subseteq \SpaceOfGambles$ denotes the linear subspace of simple gambles on $\mathcal{F}$, \ie scaled and added indicator gambles of a finite number of disjoint sets (\cf \citep[Definition 4.2.12]{rao1983theory}). 
Let $\mathcal{F}_\sigma \subseteq 2^\Omega$ be a $\sigma$-algebra. Then $\operatorname{B}(\Omega, \mathcal{F}_\sigma) \subseteq \SpaceOfGambles$ denotes the linear subspace of all bounded, real-valued, $\mathcal{F}_\sigma$-measurable gambles.
\begin{proposition}[Finitely Additive Probability on Pre-Dynkin-System and its Partial Expectation]
\label{prop:Finitely Additive Probability on Pre-Dynkin-System and its Partial Expectation}
    Let $\mathcal{D} \subseteq 2^\Omega$ be a pre-Dynkin-system. Let $\mu \colon \mathcal{D} \rightarrow [0,1]$ be a finitely additive probability defined on the pre-Dynkin-system with block structure $\{ A_i \}_{i \in I}$. Then $\mu$ is in one-to-one correspondence to a partial expectation $\pE \colon \bigcup_{i \in I} \operatorname{S}(\Omega, \mathcal{A}_i) \rightarrow \Reals$ defined on the union of linear spaces of simple gambles induced by all blocks $\mathcal{A}_i$ of $\mathcal{D}$.
\end{proposition}
\begin{proof}
    By Theorem~\ref{thm:pre-dynkin-systems are made out of algebras} we can decompose the pre-Dynkin-system $\mathcal{D}$ into a set of blocks $\{ \mathcal{A}_i\}_{i \in I}$. Since, $\mathcal{A}_i \subseteq 2^\Omega$ for all $i \in I$, each of the blocks induces a the linear subspace of simple gambles $\operatorname{S}(\Omega, \mathcal{A}_i) \subseteq \SpaceOfGambles$. Given a finitely additive measure $\mu \colon \mathcal{D} \rightarrow [0,1]$, we now define $\pE \colon \bigcup_{i \in I} \operatorname{S}(\Omega, \mathcal{A}_i) \rightarrow \Reals$ by
    \begin{align*}
        \pE(f) \coloneqq \int f d\mu|_{\mathcal{A}_i} \text{ if } f \in \operatorname{S}(\Omega, \mathcal{A}_i).
    \end{align*}
    For every $i \in I$, $\pE|_{\operatorname{S}(\Omega, \mathcal{A}_i)}$ is a linear prevision in one-to-one correspondence to the finitely additive probability $\mu|_{\mathcal{A}_i}$ \citep[Theorem 3.2.2]{walley1991statistical}. Hence, conditions (a) and (b) in Definition~\ref{def:partial expectation} are met \citep[Theorem 2.8.4]{walley1991statistical}. For $f \in \operatorname{S}(\Omega, \mathcal{A}_i) \cap \operatorname{S}(\Omega, \mathcal{A}_j)$ ($i \neq j, i,j \in I$) we know that $f \in \operatorname{S}(\Omega, \mathcal{A}_i \cap \mathcal{A}_j)$ (Lemma~\ref{lemma:intersection of simple gamble spaces}), hence,
    \begin{align*}
        \int f d\mu|_{\mathcal{A}_i} = \int f d\mu|_{\mathcal{A}_i \cap \mathcal{A}_j}  = \int f d\mu|_{\mathcal{A}_j}.
    \end{align*}
    Thus, $\pE$ is well-defined and there is no other partial expectation which agrees with $\mu$ on the indicator gambles of the sets in $\mathcal{D}$.
\end{proof}
The attentive reader might have noticed that we defined the partial expectation in Proposition~\ref{prop:Finitely Additive Probability on Pre-Dynkin-System and its Partial Expectation} on very specific linear subspaces of $\SpaceOfGambles$, namely the linear subspaces of simple gambles. In fact, the statement would still hold when enlarging the linear subspaces of simple gambles $\operatorname{S}(\Omega, \mathcal{A}_i)$ for every $i \in I$ to linear subspaces of functions which are ``convergence in measure''-approximated by gambles in $\operatorname{S}(\Omega, \mathcal{A}_i)$. For more details we refer the reader to \citep[Definition 4.4.5 and Corollary 4.4.9]{rao1983theory}.

But, it is not the case that we can extend the definition to all sets of bounded, $\mathcal{A}_i$-measurable functions, \ie bounded functions whose pre-images of sets in the smallest algebra which contains all open sets of the real numbers are contained in $\mathcal{A}_i$. For an algebra $\mathcal{A} \subseteq 2^\Omega$ the set of $\mathcal{A}$-measurable gambles is not necessarily a linear subspace of $\SpaceOfGambles$ \citep[page 129]{walley1991statistical}.

This is different for a $\sigma$-algebra $\mathcal{A}_\sigma \subseteq 2^\Omega$. The set of bounded, $\mathcal{A}_\sigma$-measurable functions $\operatorname{B}(\Omega, \mathcal{A}_\sigma)$ forms a linear subspace of $\SpaceOfGambles$ \citep[page 129]{walley1991statistical}. Here, measurability is defined as the pre-image of every Borel-measurable set in $\Reals$ is in $\mathcal{A}$. 
In this case, the set of linear spaces on which the partial expectation is defined is given by all bounded, measurable functions on the $\sigma$-blocks. 
\begin{proposition}[Finitely Additive Probability on Dynkin-Systems and its Partial Expectation]
\label{prop:Finitely Additive Probability on Dynkin-System and its Partial Expectation}
    Let $\mathcal{D}_\sigma \subseteq 2^\Omega$ be a Dynkin-system on the base set $\Omega$. Let $\mu \colon \mathcal{D}_\sigma \rightarrow [0,1]$ be a finitely additive probability defined on the Dynkin-system. Then $\mu$ is in one-to-one correspondence to a partial expectation $\pE \colon \bigcup_{i \in I} \operatorname{B}(\Omega, \mathcal{A}_i) \rightarrow \Reals$ defined on the union of linear spaces of measurable gambles induced by all $\sigma$-blocks $\mathcal{A}_i$ of $\mathcal{D}_\sigma$.
\end{proposition}
\begin{proof}
    By Theorem~\ref{thm:dynkin-systems are made out of sigma-algebras} we can decompose the Dynkin-system $\mathcal{D}_\sigma$ into a set of $\sigma$-blocks $\{ \mathcal{A}_i\}_{i \in I}$. Since, $\mathcal{A}_i \subseteq 2^\Omega$ for all $i \in I$, each of the $\sigma$-blocks induce a linear subspace of $\SpaceOfGambles$, which we denote as $\operatorname{B}(\Omega, \mathcal{A}_i)$. Given a finitely additive measure $\mu \colon \mathcal{D}_\sigma \rightarrow [0,1]$, we now define $\pE \colon \bigcup_{i \in I} \operatorname{B}(\Omega, \mathcal{A}_i) \rightarrow \Reals$ by
    \begin{align*}
        \pE(f) \coloneqq \int f d\mu|_{\mathcal{A}_i} \text{ if } f \in \operatorname{B}(\Omega, \mathcal{A}_i).
    \end{align*}
    For every $i \in I$, $\pE|_{\operatorname{B}(\Omega, \mathcal{A}_i)}$ is a linear prevision in one-to-one correspondence to the finitely additive probability $\mu|_{\mathcal{A}_i}$ \citep[Theorem 3.2.2]{walley1991statistical}. Hence, conditions (a) and (b) in Definition~\ref{def:partial expectation} are met \citep[Theorem 2.8.4]{walley1991statistical}. For $f \in \operatorname{B}(\Omega, \mathcal{A}_i) \cap \operatorname{B}(\Omega, \mathcal{A}_j)$ ($i \neq j, i,j \in I$) we know that $f \in \operatorname{B}(\Omega, \mathcal{A}_i \cap \mathcal{A}_j)$ (Lemma~\ref{lemma:intersection of measurable gamble spaces}), hence,
    \begin{align*}
        \int f d\mu|_{\mathcal{A}_i} = \int f d\mu|_{\mathcal{A}_i \cap \mathcal{A}_j}  = \int f d\mu|_{\mathcal{A}_j}.
    \end{align*}
    Thus, $\pE$ is well-defined and there is no other partial expectation which agrees with $\mu$ on the indicator gambles of the sets in $\mathcal{D}$.
\end{proof}
It remains to emphasize that there are partial expectations defined on families of linear subspaces which are not induced by finitely additive probabilities on Dynkin-systems. A simple example is given by a linear space which does not contain the constant gamble corresponding to the indicator gamble of the set $\Omega$. Hence, the definition of a partial expectation is indeed a generalization of the definition of a finitely additive probability on a pre-Dynkin-system.

Under the name ``partially specified probabilities'' \citet{Lehrer2006PartiallySpecifiedPD} introduced a closely related notion to our partial expectation. Lehrer, however, assumed that there is by definition an underlying probability distribution over the entire base set (or better said, a $\sigma$-algebra on the base set). Hence, his partially specified probabilities are by definition extendable (see Definition~\ref{def:general extendability for partial expectations}), a fact, which he implicitly exploited by re-defining the natural extension following \citep[Lemma 3.1.3 (e)]{walley1991statistical} of partially specified probabilities \citep[\S 3.2]{Lehrer2007CoherentRM}. Lehrer did not ask for the structure of the set of gambles with precise expectations, nor did he draw any connection to Walley's work, nor did he link his ``partially specified probabilities'' to finitely additive probabilities on pre-Dynkin-systems.

\subsection{System of Precision -- The Space of Gambles With Precise Expectations}
\label{sec:system of precision}
In Section~\ref{From Porbability measures to Dynkin systesms} we have shown that imprecise probabilities are precise on (pre-)Dynkin-systems. The natural analogue of this \emph{set structure of precision} is the \emph{space of gambles with precise expectation}, which actually forms a linear subspace.
\begin{theorem}(Imprecise Expectations Are Precise on a Linear Subspace of Precise Gambles)
\label{thm:Imprecise Expectations Are Precise on a Linear Subspace of Precise Gambles}
    Let $\SpaceOfGambles$ be the linear space of bounded, real-valued functions on $\Omega$. Let $L \colon \SpaceOfGambles \rightarrow \Reals$ and $U\colon \SpaceOfGambles \rightarrow \Reals$ be two functionals, for which all the following properties hold:
    \begin{enumerate}[(a)]
        \item Normalization: $L(\chi_\Omega) = U(\chi_\Omega) = 1$.
        \item Conjugacy: $U(f) = -L(-f)$ for $f \in \SpaceOfGambles$.
        \item Subadditivity of $U$: for $f,g \in \SpaceOfGambles$ we have $U(f + g) \le U(f) + U(g)$.
        \item Superadditivity of $L$: for $f,g \in \SpaceOfGambles$ we have $L(f + g) \ge L(f) + L(g)$.
        \item Positive Homogeneity: for $\alpha \in [0, \infty)$ and $f \in \SpaceOfGambles$ we have $L(\alpha f) = \alpha L(f)$ and $U(\alpha f) = \alpha U(f)$.
    \end{enumerate}
    Then $L$ and $U$ coincide on a linear space $\setofprecisegambles \subseteq \SpaceOfGambles$, the \emph{space of gambles with precise expectation}, which contains all constant gambles.
\end{theorem}
\begin{proof}
We define 
\begin{align}
\label{eq:setofprecisegambles}
    \setofprecisegambles \coloneqq \{ f \in \SpaceOfGambles \colon L(f) = U(f)\},
\end{align}
and show that $\setofprecisegambles$ forms a linear subspace of $\SpaceOfGambles$.
First, let $f,g \in \setofprecisegambles$, then
\begin{align*}
    L(f) + L(g) \overset{(d)}{\le}\le L(f + g) \overset{(\star)}{\le} U(f + g) \overset{(c)}{\le} U(f) + U(g) \overset{Eq.\ref{eq:setofprecisegambles}}{=} L(f) + L(g).
\end{align*}
For ($\star$) observe that $L(f) \le U(f)$ for all $f \in \SpaceOfGambles$, since
\begin{align*}
    L(f) + L(-f) \le L(0) \overset{(a),(e)}{=} 0 \overset{(a),(e)}{=} U(0) \le U(f) + U(-f),
\end{align*}
we have,
\begin{align*}
    L(f) + L(-f) \le U(f) + U(-f) \Leftrightarrow L(f) - U(f) \le U(f) - L(f) \Leftrightarrow  L(f) \le U(f).
\end{align*}
Second, let $f \in \setofprecisegambles$ and $\alpha \in \Reals$. If $\alpha \ge 0$, then
\begin{align*}
    L \left (\alpha f \right) \overset{(e)}{=} \alpha L \left ( f \right) \overset{Eq.\ref{eq:setofprecisegambles}}{=} \alpha U \left ( f \right) \overset{(e)}{=} U \left ( \alpha f \right).
\end{align*}
Otherwise,
\begin{align*}
    L \left (\alpha f \right) \overset{(b)}{=} -U \left (-\alpha f \right) \overset{(e)}{=} \alpha U \left ( f \right) \overset{Eq.\ref{eq:setofprecisegambles}}{=} \alpha L \left ( f \right) \overset{(b)}{=} -\alpha U \left (- f \right) \overset{(e)}{=}  U \left (\alpha f \right).
\end{align*}
Third, $\setofprecisegambles$ contains all constant gambles by (a), (b) and (e). Hence, we have shown that $\setofprecisegambles$ forms a linear subspace of $\SpaceOfGambles$ which contains all constant gambles.
\end{proof}
The choice of properties for the lower and upper expectation functional is not arbitrary. We tried to resemble the properties involved in the analogous statement for lower and upper probabilities (Theorem~\ref{thm:lower and upper probablity define Dynkin-system}). One can easily check that a lower and upper expectation with the given properties (a) -- (d) forms a lower and upper probability as required in Theorem~\ref{thm:lower and upper probablity define Dynkin-system} if the expectation is restricted to indicator gambles. However, we added property (e), positive homogeneity.

Without the property of positive homogeneity, the resulting set of gambles with precise expectations would not form a proper linear subspace, as then one can only guarantee closedness of $\setofprecisegambles$ under rational multiplication.
The condition of positive homogeneity ``fills up'' the gaps with all real-scaled functions. Instead of positive homogeneity one can as well demand a continuity assumption of the lower and upper functional $L$ and $U$, \eg \citep[Property (l) Theorem 2.6.1]{walley1991statistical}. We emphasize that coherent previsions (see Definition~\ref{def:coherent prevision}) fulfill all of the demanded properties \citep[Theorem 2.6.1]{walley1991statistical}.

Interestingly, the restriction $L|_\setofprecisegambles = U|_\setofprecisegambles$ is not necessarily a partial expectation. Otherwise it would form a coherent linear prevision (see Definition~\ref{def:coherent prevision}). This is different compared to Theorem~\ref{thm:lower and upper probablity define Dynkin-system}, where the lower and upper probability actually define a finitely additive probability on the set structure of precision, which is \emph{not} necessarily coherent. However, those two statements are not in contradiction. Any pair of lower and upper expectations, as we defined them here, induce a unique lower and upper probability. The resulting finitely additive probability on the set structure of precision gives rise to a partial expectation (Proposition~\ref{prop:Finitely Additive Probability on Pre-Dynkin-System and its Partial Expectation}) on a set of linear subspaces contained in the space of gambles with precise expecation $\setofprecisegambles$ of $L$ and $U$.

The converse direction, however, is not true. There is no unique lower and upper expectation functional with the given properties associated to a lower and upper probability fulfilling the axioms of Theorem~\ref{thm:lower and upper probablity define Dynkin-system} \citep[\S 2.7.3]{walley1991statistical}. 
Concluding, lower and upper expectation as defined here are not the ``perfect'' analogues of lower and upper probabilities.

This as well explains the mismatch between systems of precision for probabilities and expectations. A lower and upper expectation fulfills the properties of a lower and upper probability but not vice-versa. Hence, only weaker statements about the system of precision are possible for probabilities.
As a result, the analogue of the set structure of precision, a pre-Dynkin-system, is the space of gambles with precise expectations, a \emph{single} linear subspace.
In Definition~\ref{def:partial expectation}, however, we equated pre-Dynkin-systems with \emph{sets} of linear subspaces. In this case, a one-to-one correspondence between a finitely additive probability on a pre-Dynkin-system, and a generalized expectation, concretely a partial expectation, can be established. Hence, the analogy of pre-Dynkin-systems and linear subspaces of gambles depends on the correspondence of probability and expectation.

\subsection{Generalized Extendability is Equivalent to Coherence}
Partial expectations are, as we have shown, a natural generalization of finitely additive probabilities on pre-Dynkin-systems. Hence, it is not far-fetched to ask for definitions of coherence and extendability again, now in the more general context. It turns out that the same story can be re-told on a more general scale:
The definition of coherent probabilities (Definition~\ref{def:coherent probability}) is in fact just the reduction of the following definition of a coherent prevision to indicator gambles.
\begin{definition}[Coherent Prevision]\citep[Definition 2.5.1]{walley1991statistical}
\label{def:coherent prevision}
    Let $L \subseteq \SpaceOfGambles$ be an arbitrary subset of the linear space of bounded functions. A functional
    $\lP\colon L \rightarrow \Reals$ is a \emph{coherent lower prevision} if and only if 
    \begin{align*}
        \sup_{\omega \in \Omega} \sum_{i = 1}^j (f_i(\omega) - \lP(f_i)) - m (f_0(\omega) - \lP(f_0)) \ge 0,
    \end{align*}
    for non-negative $n,m \in \Naturals$ and $f_0, f_1, \ldots f_n \in L$. If $L = - L$, the conjugate \emph{coherent upper prevision}
    is given by $\uP(f) \coloneqq -\lP(-f)$ for all $f \in L$. If furthermore $\lP(f) = \uP(f)$ for all $f \in L$, we call $\aP \coloneqq \lP = \uP$ a \emph{linear  prevision}.
\end{definition}
This definition of coherent previsions is substantiated by consistency of gamblers regarding their betting behavior on gambles with uncertain outcome, \eg \citep[\S 2.3.1]{walley1991statistical}. Importantly, the rather opaque, but general definition of coherence can be simplified greatly for coherent previsions defined on linear subspaces of $\SpaceOfGambles$. Theorem 2.5.5 in \citep{walley1991statistical} shows that coherence for lower previsions on linear subspaces can be expressed as superadditivity, positive homogeneity, and accepting sure gains (see \citep[Definition 2.3.3]{walley1991statistical}).

Having introduced the notion of a coherent prevision, we now envisage the link between partial expectations and coherent previsions. We introduced extendability for finitely additive probabilities on pre-Dynkin-systems as a useful property. It guarantees that the probability can ``nicely'' be embedded into ``larger'' finitely additive probability which is defined on an encompassing algebra. Hence, the analogue for partial expectations is straightforward.
\begin{definition}[General Extendability]
\label{def:general extendability for partial expectations}
    A partial expectation $\pE\colon  \bigcup_{i \in I} L_i \rightarrow \Reals$ is \emph{extendable} if and only if
    there exists a partial expectation $\pE' \colon \SpaceOfGambles \rightarrow \Reals$ such that $\pE'|_{\bigcup_{i \in I} L_i } = \pE$.
\end{definition}
Interestingly, the extendability condition provided in Theorem~\ref{theorem:extendability condition - Horn-Tariski} has a (more general) cousin adapted to the setting of gambles instead of events.
\begin{proposition}[Extendability Condition for Previsions] (\cf \citep[Theorem 6.1]{maharam1972consistent})
    Let $\{ L_i \}_{i \in I}$ be a non-empty family of linear subspaces of $\SpaceOfGambles$.
    A \emph{partial expectation} $\pE \colon \bigcup_{i \in I} L_i \rightarrow \Reals$ is extendable if and only if for every finite collection of functions
    $f_1, \ldots , f_n \in \{ L_i \}_{i \in I}$,
    \begin{align*}
        \sum_{i = 1}^n f_i \ge 0 \implies \sum_{i = 1}^n \pE(f_i) \ge 0.
    \end{align*}
\end{proposition}
\begin{proof}
    It seems that Theorem 6.1 \citep{maharam1972consistent} is equivalent to our statement. However, there is a subtlety which we want to argue here is indeed irrelevant.
    Extendability of a partial expectation requires the existence of a positive, \emph{normed}, linear functional on $\SpaceOfGambles$, whose restriction on the according linear subspaces coincides with the partial expectation. Theorem 6.1 in \citep{maharam1972consistent} only guaratees that a positive, linear functionals exists. But, normedness of such functional is automatically given if $\chi_\Omega \in L_i$ for some $i \in I$. Otherwise, we extend the partial expectation $\pE$ to $\pE ' \colon \{ \alpha \chi_\Omega \colon \alpha \in \Reals\}\bigcup_{i \in I} L_i \rightarrow \Reals$ such that
    \begin{align*}
        \pE'(f) \coloneqq \begin{cases}
    \alpha \text{ if } f \in \{ \alpha \chi_\Omega \colon \alpha \in \Reals\}\\
    \pE(f) \text{ otherwise.}
    \end{cases}
    \end{align*} Then again, Theorem 6.1 \citep{maharam1972consistent} applies.
\end{proof}
Against the background that extendability and coherence define the same concept for finitely additive probabilities on pre-Dynkin-systems, the resulting equivalence of extendability and coherence for partial expectations is of little surprise.
\begin{proposition}[Extendability is Equivalent to Coherence]
    Let $\{ L_i \}_{i \in I}$ be a non-empty family of linear subspaces of $\SpaceOfGambles$. The partial expectation $\pE \colon \bigcup_{i \in I} L_i \rightarrow \Reals$ is extendable if and only if $\pE$ is a linear prevision, \ie is coherent.
\end{proposition}
\begin{proof}
    If $\pE$ is a linear prevision on $\bigcup_{i \in I} L_i$, then there is a linear prevision $\pE' \colon \SpaceOfGambles \rightarrow \Reals$ such that $\pE'|_{ \bigcup_{i \in I} L_i} = \pE$ \citep[Theorem 3.4.2]{walley1991statistical}.
    Conversely, if $\pE$ is an extendable partial expectation, then its extension is obviously a linear prevision, hence it is coherent. The restriction of a coherent linear prevision to any subset of gambles is coherent (and linear).
\end{proof}

\subsection{A Duality Theory for Previsions and Families of Linear Subspaces}
In Section~\ref{sec:the credal set and its relation to Dynkin-system structure} we step by step spelled out an order relationship between the set structure of precision and credal sets, a model for (coherent) imprecise probabilities.
Naturally the presented generalization begs the question whether a related relationship between credal sets and the spaces of gambles with precise expectations exists. We answer affirmatively. We redefine the credal and dual credal set function and shortly discuss its analogous properties. Again, we require a ``reference measure''. In this case, it is a fixed linear prevision  $\fixedaP$ on the space of all gambles $\SpaceOfGambles$, which is indeed in one-to-one correspondence to a finitely additive probability measure on $2^\Omega$.
\begin{definition}[Generalized Credal Set Function]
\label{def:generalized credal set function}
    Let $\setoflinearprevisions$ be the set of linear previsions on the Banach space $\SpaceOfGambles$. For a fixed linear prevision $\fixedaP \in \setoflinearprevisions$ we call
    \begin{align*}
        \mbageneral \colon 2^{\SpaceOfGambles} \rightarrow 2^\setoflinearprevisions, \quad \mbageneral(\mathcal{G}) \coloneqq \{ \aP \in \setoflinearprevisions\colon \aP(g) = \fixedaP(g), \forall g \in \mathcal{G}\},
    \end{align*}
    the \emph{generalized credal set function}.
\end{definition}
\begin{definition}[Generalized Dual Credal Set Function]
\label{def:generalized dual credal set function}
    Let $\setoflinearprevisions$ be the set of linear previsions on the Banach space $\SpaceOfGambles$. For a fixed linear prevision $\fixedaP \in \setoflinearprevisions$ we call
    \begin{align*}
        \mdbageneral \colon 2^\setoflinearprevisions \rightarrow 2^{\SpaceOfGambles} , \quad \mdbageneral(\mathcal{Q}) \coloneqq \{ g \in \SpaceOfGambles\colon \aP(g) = \fixedaP(g), \forall \aP \in \mathcal{Q}\},
    \end{align*}
    the \emph{generalized dual credal set function}.
\end{definition}
Why can we call those functions ``generalized''? Simply, because any system of sets is equivalently represented as its set of indicator gambles which span their own linear space of simple gambles, \ie linear combinations of indicator gambles.

The generalized credal set function maps, as the credal set function in Definition~\ref{def:Credal Set Function}, to weak$^\star$-closed, convex subsets of $\setoflinearprevisions$.
The generalized dual credal set function, however, reveals a first subtlety. It maps to linear subspaces of $\SpaceOfGambles$. The dual credal set function following Definition~\ref{def:Dual Credal Set Function} mapped to pre-Dynkin-systems. In Proposition~\ref{prop:Finitely Additive Probability on Pre-Dynkin-System and its Partial Expectation} and Proposition~\ref{prop:Finitely Additive Probability on Dynkin-System and its Partial Expectation} families of linear subspaces were the analogues of (pre-)Dynkin-systems. Here, a single linear subspace is the analogue of a pre-Dynkin-system. For a first step towards an explanation of this asymmetry see Section~\ref{sec:system of precision}. 
Finally, the pair of functions constitute a Galois connection.
\begin{proposition}[Properties of Generalized (Dual) Credal Set Function]
\label{prop:Properties of Generalized (Dual) Credal Set Function}
    Let $\mbageneral$ be a generalized credal set function and $\mdbageneral$ be a generalized dual credal set function. All the following properties hold:
    \begin{enumerate}[(a)]
        \item The generalized credal set function $\mbageneral$ maps to weak$^\star$-closed, convex sets.
        \item The generalized dual credal set function $\mdbageneral$ maps to a linear subspace.
        \item The generalized credal set function $\mbageneral$ and generalized dual credal set function $\mdba$ form a Galois connection.
    \end{enumerate}
\end{proposition}
\begin{proof}
\begin{enumerate}[(a)]
    \item We have fixed $\fixedaP$ to a linear prevision. Hence, it is coherent. For any $\mathcal{G} \subseteq \SpaceOfGambles$, $\mbageneral(\mathcal{G})$ is the set of all linear previsions which dominate $\fixedaP$ on $\mathcal{G}$. Theorem 3.6.1 in \citep{walley1991statistical} then states that this set is weak$^\star$-closed and convex.
    \item Let $\mathcal{Q} \subseteq \setoflinearprevisions$.
    \begin{description}
        \item [Additivity]
        Let $f,g \in \mdbageneral(\mathcal{Q})$. Then, for all $\aP \in \setoflinearprevisions$,
        \begin{align*}
            \aP(f + g) = \aP(f) + \aP(g) = \fixedaP(f) + \fixedaP(g) = \fixedaP(f + g),
        \end{align*}
        \ie $f + g \in  \mdbageneral(\mathcal{Q})$.

        \item [Homogeneity]
        Let $f \in \mdbageneral(\mathcal{Q})$ and $\alpha \in \Reals$. Then, for all $\aP \in \setoflinearprevisions$,
        \begin{align*}
            \aP(\alpha f) = \alpha \aP(f)= \alpha \fixedaP(f)= \fixedaP(\alpha f),
        \end{align*}
        \ie $\alpha f \in  \mdbageneral(\mathcal{Q})$.
        For homogeneity we need the easy fact that a linear prevision is not only positive homogeneous, but generally homogeneous. For this consider a linear prevision $\aP$ and any gamble $f \in \SpaceOfGambles$ with $\alpha < 0$, then $\aP(\alpha f) = - \aP(- \alpha f) =  \alpha \aP( f)$.
    \end{description}
    \item The two functions constitute a Galois connection (\cf Proposition~\ref{prop:galois connection by (dual) credal set function}), $\mathcal{G} \subseteq \mdbageneral(\mathcal{Q}) \Leftrightarrow \mathcal{Q} \subseteq \mbageneral(\mathcal{G})$. To this end, we show the left to right implication,
    \begin{align*}
        \aP \in \mathcal{Q} \Rightarrow \aP(g) = \fixedaP(g), \forall g \in \mathcal{G} \Rightarrow \aP \in \mbageneral(\mathcal{G}),
    \end{align*}
    and the right to left implication,
    \begin{align*}
        g \in \mathcal{G} \Rightarrow \aP(g) = \fixedaP(g), \forall \aP \in \mathcal{Q} \Rightarrow g \in \mdbageneral(\mathcal{Q}).
    \end{align*}
\end{enumerate}
This concludes the proof.
\end{proof}
Galois connections possess a series of helpful properties \citep[\S V.7 and V.8]{birkhoff1940lattice}. For instance, they give rise to a bipolar-closure operator. A non-empty subset $\mathcal{Q} \subseteq \setoflinearprevisions$ is bipolar-closed if and only if $\mathcal{Q} = \mbageneral(\mdbageneral(\mathcal{Q}))$. A non-empty subset $\mathcal{G} \subseteq \SpaceOfGambles$ is bipolar-closed if and only if $\mathcal{G} = \mdbageneral(\mbageneral(\mathcal{G}))$.
Furthermore, Proposition~\ref{prop:Properties of Generalized (Dual) Credal Set Function} provides necessary conditions for bipolar-closed sets. For instance, a bipolar-closed set $\mathcal{G} \subseteq \SpaceOfGambles$ is a linear subspace. But is every such linear subspace a  bipolar-closed set? No. 
{
\color{darkgray}
\begin{example}
    Let $\Omega_2 \coloneqq \{ 1,2\}$, then $\operatorname{B}(\Omega_2) = \{ \alpha_1 \chi_{\{ 1\}} + \alpha_2 \chi_{\{ 2\}} \colon \alpha_1, \alpha_2 \in \Reals \}$. Hence, linear functionals on $\operatorname{B}(\Omega_2)$ are defined via their behavior on the basis. Let $\fixedaP(\chi_{\{1\}}) = 1$. Then, $\{ \alpha_1 \chi_{\{ 1\}} \colon \alpha_1 \in \Reals\} \subseteq \operatorname{B}(\Omega_2)$ is a linear subspace, but it is, as one can easily check, not bipolar-closed, because the demand for normalization of any linear prevision $\aP$ which coincides with $\fixedaP$ on $\chi_{\{1 \}}$ requires $\aP(\chi_{\{ 2\}}) = 0$.
\end{example}
}
Again it seems to be more intricate than expected to characterize bipolar-closed sets. For pre-Dynkin-systems and finitely additive measures we already collected some first hints that sets of measure zero play an important role in the characterization of bipolar-closed sets. In the case of linear subspaces and linear previsions we observe a similar ``combinatorial restriction''. In order to improve understanding, let us replace $2^\Pba$ by $2^{\SpaceOfLinearFunctionalsOnGambles}$ in Definition~\ref{def:generalized credal set function}, Definition~\ref{def:generalized dual credal set function} and Proposition~\ref{prop:Properties of Generalized (Dual) Credal Set Function}\footnote{The proposition still holds.}, which is equivalent to stating that linear previsions are not necessarily normalized, nor positive. Then, by leveraging the Hahn-Banach-type Theorem 1.5.14 in \citep{rao1983theory}, one can easily see that linearity of a subset $\mathcal{G} \subseteq \SpaceOfGambles$ is not only a necessary, but as well a sufficient condition for bipolar-closedness for those modified ``credal set functions''. Thus, the restriction to actual linear previsions makes the characterization of bipolar-closed sets more complex. A compelling, more exhaustive answer still waits to be found.

Analogous to the discussion in Section~\ref{sec:interpolation from algebra to trivial dynkin-system}, it is possible to provide a lattice duality and interpolation scheme via the generalized (dual) credal set functions. Instead of the lattice of pre-Dynkin-systems $(\mathfrak{D}, \subseteq)$ the interpolation is directed by the lattice of linear subspaces $(\mathcal{L}, \subseteq)$ of $\SpaceOfGambles$. As commonly known, the lattice of linear subspaces has the two operations $L_1 \wedge L_2 \coloneqq L_1 \cap L_2$ and $L_1 \vee L_2 \coloneqq \lin (L_1 \cup L_2)$\footnote{We denote the linear span with $\lin$.}. Its minimal element is the trivial zero vector linear subspace $\{ 0\}$. Its maximal element is the entire space of all gambles $\SpaceOfGambles$.
Due to higher generality of the here presented (dual) credal set function, the interpolation provided is more fine-grained than for the previously given interpolation by pre-Dynkin-systems. The following set containment (trivially) holds:
\begin{align*}
    \{ \mbageneral(\mathcal{D}) \colon \mathcal{D} \in \mathfrak{D}\} \subseteq \{ \mbageneral(L) \colon L \in \mathcal{L} \},
\end{align*}
where $\mbageneral(\mathcal{D}) = \mbageneral(\{ \chi_D \colon D \in \mathcal{D}\})$. In other words, the lattice of pre-Dynkin-systems is ``contained'' in the lattice of of linear subspaces. However, as for pre-Dynkin-systems the interpolation via linear subspaces is improper. The reason for this is again that not every linear subspace is bipolar-closed. By restriction to linear, bipolar-closed subspaces one can clean up the setup. For details we refer to Section~\ref{sec:interpolation from algebra to trivial dynkin-system}. We do not make explicit the detailed reiteration of the same argument here.

In summary, we confirmed our findings of Section~\ref{sec:the credal set and its relation to Dynkin-system structure} extended to previsions. The generalized dual lattice setup underlines the structural consistency between the system of precision and its corresponding imprecise probability.

\section{Conclusion and Open Questions}
In this paper, we have explicated relations between the systems of precision and imprecise probabilities (respectively expectations).
First, we have shown that the system of precision forms a pre-Dynkin-system (respectively a linear subspace). This structural insight raises a series of follow-up questions: How does the system of precision of a coherent prevision relate to the set of desirable gambles of this prevision? How does the preference ordering change the set structure of precision for the corresponding beliefs? What is the role of coherence with respect to the system of precision?

Second, we defined finitely additive probabilities on pre-Dynkin-systems. The equivalence of extendability and coherence of such probabilities strengthens the link between quantum probability and imprecise probability. We speculate that further insights can be obtained by exploiting this relationship.
In addition, the generalization of finitely additive probabilities on pre-Dynkin-systems to partial expectations directly opens the door to machine learning applications. In robust machine learning the expected risk minimization framework is extended to more general expectation functionals. Partial expectation can, possibly after more computational investigations, deliver the desired robustness against dependencies in specific domains, such as privacy preservation, ``not-missing-at-random'' features, restricted data base access or multi-measurement data.

Finally, we developed a duality theory of systems of precision and imprecise probabilities (respectively expectations).
A Galois connection defines a parametrized family of imprecise probabilities
which follow an order structure provided by the lattice of pre-Dynkin-systems (respectively the lattice of linear subspaces).

In modern statistics, especially in machine learning, probabilistic statements are increasingly tailored to individuals. Individual probabilistic statements, however, require justification. One can interpret probabilities on pre-Dynkin-systems as probabilities which do not allow for such statements in the first place. One could perceive this fact as a weakness. We, in contrast, embrace its strength, when for ethical, legislative or other reasons individualistic ascriptions are harmful, unjustifiable, forbidden or not desirable. We provide a first, rough interpolation scheme via the lattice duality. It demonstrates the space of adjustability of probabilistic assumptions in real-world scenarios. The involved pre-Dynkin-systems are mathematical definitions of levels of group resolution. The question of how to choose such set-systems is related to the questions of intersectionality.

Several fundamental, technical questions remain open: how does the lattice duality imposed by pre-Dynkin-systems or linear spaces
relate to other dualities, such as convex duality, exploited in the field of imprecise
probability. Can one easily characterize the bipolar-closed sets? Why is there no clear analogy between pre-Dynkin-systems and linear subspaces?

We leave this collection of intriguing questions open to future work contributing to
an understanding of the system of precision and the imprecise probability model.

{
\small
\section{Acknowledgements}
The authors thank the International Max Planck Research School for
Intelligent System (IMPRS-IS) for supporting Rabanus Derr. 
Many thanks to Christian Fröhlich
for helpful discussions and feedback.

We thank all anonymous reviewers on a first version for the detailed and helpful feedback. In particular, with their help, the statements and proofs of Theorem~\ref{thm:Extendability Gives Coherence} and Theorem~\ref{thm:Extension Theorem - Finitely Additive Case} have been cleaned and shortened substantially.

\subsection{Funding}
This work was funded in part by the Deutsche Forschungsgemeinschaft
(DFG, German Research Foundation) under Germany’s Excellence Strategy –-
EXC number 2064/1 –- Project number 390727645. We acknowledge support from the Open Access Publication Fund of the University of Tübingen.

\subsection{Contributions of Authors}
Conceptualization, Rabanus Derr and Robert C. Williamson; Formal analysis, Rabanus Derr; Funding acquisition, Robert C. Williamson; Investigation, Rabanus Derr; Methodology, Rabanus Derr; Project administration, Robert C. Williamson; Supervision, Robert C. Williamson; Validation, Robert C. Williamson; Visualization, Rabanus Derr; Writing – original draft, Rabanus Derr; Writing – review \& editing, Robert C. Williamson.
}

\bibliographystyle{plainnat}
\bibliography{main}

\appendix

\section{Lemmas and Proofs}

\subsection{Compatibility Structure}
\label{appendix:compatibility structure}
We have shown in Theorem~\ref{thm:pre-dynkin-systems are made out of algebras} that every
pre-Dynkin-system can be decomposed into a union of algebras. A trivial follow-up question remains to be answered:
is every union of algebras a pre-Dynkin-system? The general answer is no.
{
\color{darkgray}
\begin{example}[Union of Algebras is not Necessarily Pre-Dynkin-System]
\label{example:Union of Algebras is not Necessarily Pre-Dynkin-System}
    Let $\Omega = \{ 1,2,3\}$. Then $\mathcal{A}_1 \coloneqq \{\emptyset, 1, 23, \Omega \}$ and $\mathcal{A}_2 \coloneqq \{\emptyset, 12, 3, \Omega \}$
    are algebras on $\Omega$. However, $\mathcal{A}_1 \cup \mathcal{A}_2$ do not form a pre-Dynkin-system, as \eg $1 \cup 3 = 13 \notin \mathcal{A}_1 \cup \mathcal{A}_2$.
\end{example}
}
However, we can give sufficient conditions under which the union over a set of algebras form a pre-Dynkin-system.
To this end, we have to introduce the so-called compatibility structure, which is made out of $\pi$-systems.
\begin{definition}[$\pi$-System]
\label{def:pi-system}
    Let $\Omega$ be an arbitrary base set. A $\pi$-system $\mathcal{I}$ is a non-empty subset of \ $2^\Omega$ such that for arbitrary
    $A_1, \ldots, A_n \in \mathcal{I}$, $\bigcap_{1 \le i \le n} A_i \in \mathcal{I}$.
\end{definition}
\begin{definition}[Compatibility Structure]
\label{def:compatibility structure}
    A set of $\pi$-systems $\{ \mathcal{I}_i \}_{ i \in I}$ is called a
     \emph{compatibility structure} of $\mathcal{A}$ if and only if the following two conditions hold:
    \begin{enumerate}[(a)]
        \item $\bigcup_{i \in I } \mathcal{I}_i = \mathcal{A}$
        \item For every $\pi$-system $\mathcal{I} \subseteq \mathcal{A}$
        there exist $i \in I$ such that $\mathcal{I} \subseteq \mathcal{I}_i$.
    \end{enumerate}
\end{definition}
A ``compatibility structure'' is a system of set systems. We call it compatibility structure, because every finite intersection
of elements in a $\pi$-system is compatible with every other finite intersection of elements from the same $\pi$-system.
The $\pi$-systems in the compatibility structure are not necessarily disjoint.
Given a set of $\pi$-systems, the union of these $\pi$-systems is not
necessarily a compatibility structure. In fact, there are possible
incompatibilities between the $\pi$-systems (\eg Example~\ref{example:Union of Algebras is not Necessarily Pre-Dynkin-System}).
It turns out that being a compatibility structure is a sufficient condition for a set of
algebras to form a pre-Dynkin-system.
\begin{theorem}[Union of Algebras is Pre-Dynkin-System if Compatibility Structure]
\label{theorem:Union of Algebras is Pre-Dynkin-system if Compatibility Structure}
    Let $\{\mathcal{A}_i \}_{i \in I}$ be a family of algebras.
    If $\{\mathcal{A}_i \}_{i \in I}$
    is a compatibility structure, then $\mathcal{A}_\cup \coloneqq \bigcup_{i \in I} \mathcal{A}_i$ is a pre-Dynkin-system.
\end{theorem}
\begin{proof}
    We show that
    $\mathcal{A}_\cup$ is a pre-Dynkin-system if it is the union
    of the compatibility structure. First,
    it contains the empty set. Second, if $A \in \mathcal{A}_\cup$
    there is an $i$ such that $A \in \mathcal{A}_i$. Then $A^c \in \mathcal{A}_i$, thus
    $A^c \in \mathcal{A}_\cup$. Third, let $\{A_j\}_{j \in [n]}$ be a
    subset of $\mathcal{A}_\cup$ with $n \in \Naturals$ such
    that $A_k \neq \emptyset$, $A_l \neq \emptyset$ and
    $A_k \cap A_l = \emptyset$ for all $k \neq l; k,l \in [n]$.
    Thus, $\{A_j\}_{j \in [n]} \cup \{ \emptyset\}$ is closed under finite intersection. Hence, it is a $\pi$-system. By definition of a compatibility
    structure there is an $i \in I$ such that $\{A_j\}_{j \in [n]} \cup \{\emptyset\} \subseteq \mathcal{A}_i$.
    Since $\mathcal{A}_i$ is an algebra, it is closed under disjoint union. Thus,
    $\bigcup_{j \in [n]} A_j \in \mathcal{A}_i \subseteq \mathcal{A}_\cup$.
\end{proof}

\subsection{Supremum of a Chain of Algebras is an Algebra}
The following lemma is used to prove that every pre-Dynkin-system
can be dissected into algebras (Theorem~\ref{thm:pre-dynkin-systems are made out of algebras}).
\begin{lemma}[Supremum of a Chain of Algebras is an Algebra]
\label{lemma:Supremum of a Chain of Algebras is an Algebra}
    Let $(\{\mathcal{A}_i\}_{i \in I}, \subseteq)$ be a non-empty chain of algebras.
    Then it has an upper bounding algebra $\mathcal{A}_{\sup}$.
\end{lemma}
\begin{proof}
    Let $\mathcal{A}_{\sup} = \bigcup_{i \in I} \mathcal{A}_{i}$.
    Trivially, $\mathcal{A}_i \subseteq \mathcal{A}_{\sup}$ for
    every $i \in I$.
    It remains to show that $\mathcal{A}_{\sup}$ is an algebra.
    It is non-empty by construction. Consider any
    finite subset $\{A_1, \ldots, A_n \} \subseteq \mathcal{A}_{\sup}$.
    Without loss of generality there exist $j \in I$ such that
    $\{A_1, \ldots, A_n \} \subseteq \mathcal{A}_j$. Thus,
    \begin{align*}
        \bigcup_{i = 1}^n A_i \in \mathcal{A}_j \subseteq \mathcal{A}_{\sup}.
    \end{align*}
    Furthermore, for every $A \in \mathcal{A}_{\sup}$ there exist $j \in I$
    such that $A \in \mathcal{A}_j$. Thus,
    \begin{align*}
        A^c \in \mathcal{A}_j \subseteq \mathcal{A}_{\sup},
    \end{align*}
    which concludes the proof.
\end{proof}

\subsection{Sufficient Conditions for Bipolar-Closed Sets}
\label{sufficient conditions for bipolar-closed sets}
In order to give sufficient conditions for bipolar-closed sets
we introduce the following notion.
\begin{definition}[Weak Atom with Respect to Pre-Dynkin-System]
\label{def:weak atom wrt pre-Dynkin-system}
    A \emph{weak atom} with respect to a pre-Dynkin-system $\mathcal{D}$ on $\Omega $
    is a set $A \subseteq \Omega$ such that every $B \in \mathcal{D}$ with
    $B \subseteq A$ is either $B = \emptyset$ or $B = A$.
\end{definition}
Weak atoms always exist, for instance the empty set is a weak atom with respect to any pre-Dynkin-system $\mathcal{D}$ on any non-empty set $\Omega$.
\begin{lemma}[Decomposition into Pre-Dynkin Set and Atom]
\label{lemma:decomposition into dynkin set and atom}
    Let $\mathcal{D} \subseteq 2^\Omega$ be a pre-Dynkin-system.
    Every element $B \in 2^\Omega \setminus \mathcal{D}$ can be expressed as a disjoint union
    of an element $B_\mathcal{D} \in \mathcal{D}$ and a weak atom $A \in 2^\Omega \setminus \mathcal{D}$
    with respect to $\mathcal{D}$.
\end{lemma}
\begin{proof}
    We define $B_\mathcal{D} \subseteq B$ as a maximal set $B_\mathcal{D} \in \mathcal{D}$ such that for all
    $H \in \mathcal{D}$ with $B_\mathcal{D} \subseteq H \subseteq B$ we have $H = B_\mathcal{D}$. By Zorn's lemma such an element always exist. To see this, consider the poset of all subsets of $B$ which are in $\mathcal{D}$ ordered by set inclusion.
    Then, we decompose $B \setminus B_\mathcal{D} = A$. 
    Clearly, $A \neq \emptyset$. Furthermore, $A \notin \mathcal{D}$ and there is no $H \in \mathcal{D}$
    such that $H \subseteq A$ except of the empty set, otherwise $B_\mathcal{D}$ would not have been maximal
    in our sense. Both statements follow from the closedness of $\mathcal{D}$ under disjoint union. Thus, $A$ is a weak atom with respect to $\mathcal{D}$.
\end{proof}
With these definitions at hand we can show the following crucial lemma.
\begin{lemma}[Weak Atoms Are Not in the Bipolar-Closed Sets -- Finite, Discrete Setting]
\label{lemma:Weak Atoms Are Not in the Bipolar-Closed Sets - Finite, Discrete Setting}
    Let $\Omega = [n]$. We fix a finitely additive probability
    $\psi$ on $2^\Omega$ such that $\psi(F)> 0$ for every $F \in 2^\Omega \setminus \{ \emptyset \}$.
    Let $\mathcal{D} \subseteq 2^\Omega$ be a pre-Dynkin-system. Then, for every weak atom $A$ with respect to $\mathcal{D}$ such that $A \notin \mathcal{D}$ we have
    \begin{align*}
        A \notin \mdba (\mba(\mathcal{D})),
    \end{align*}
    where $\mba$ and $\mdba$ are defined following Definition~\ref{def:Credal Set Function} and~\ref{def:Dual Credal Set Function}, respectively.
\end{lemma}
\begin{proof}
    We show the Lemma by constructing, for every weak atom $A \subseteq \Omega$ with respect to $\mathcal{D}$ such that $A \notin \mathcal{D}$, a probability measure on $2^\Omega$ which coincides with $\psi$ on the pre-Dynkin-system $\mathcal{D}$, but differs to $\psi$ on the event $A$.
    Observe, weak atoms $A \notin \mathcal{D}$
    exist as long as $\mathcal{D} \neq 2^\Omega$ (by Lemma~\ref{lemma:decomposition into dynkin set and atom}).
    Second, the weak atom $A$ contains an
    element $i \in A$ such that $\{ i\} \notin \mathcal{D}$.
    Otherwise, $A \in \mathcal{D}$.
    
    So, we fix an arbitrary weak atom $A$ with respect to $\mathcal{D}$ such that $A \notin \mathcal{D}$.
    Due to the finiteness of $[n]$ we can represent every measure $\nu \in \Pba$ as
    $\nu_1,\ldots, \nu_n$ with the constraints $0 \le \nu_k \le 1$ and $\sum_{k  = 1}^n \nu_k = 1$.
    With this in mind, we construct a probability $\nu \in \Pba$ such 
    that $\nu(B) = \psi(B)$ for all $B \in \mathcal{D}$ but $\nu(A) \neq \psi(A)$.
    
    To this end, we define $\nu$ on all elements in $[n]$. We choose
    \begin{align*}
        \nu_{i} = \psi_{i} + \epsilon
    \end{align*}
    for the $i \in A$ specified earlier,
    with $\epsilon \coloneqq \min_{k \in [n]} \psi_k$.
    By assumption, $\epsilon > 0$.
    
    By Theorem~\ref{thm:pre-dynkin-systems are made out of algebras} we know
    that we can decompose $\mathcal{D}$ into sub-algebras. In fact, we can decompose $\mathcal{D}$ into a finite set of sub-algebras as $\mathcal{D}$ is finite itself. Furthermore, set algebras
    on finite sets are build upon a finite set of atoms \citep[Remark 1.1.17.(2)]{rao1983theory}. Those atoms form a partition of the finite set, hence each set algebra on a finite set is in one-to-one correspondence to a partition (\cf \citep[Proposition 1.1.18]{rao1983theory}). Hence, $\mathcal{D}$ determines a
    family of partitions corresponding to its blocks. In detail, $\mathcal{D} = \bigcup_{o \in O} \mathcal{A}_o$ for set algebras $\mathcal{A}_o$ and finite $O$ (Theorem~\ref{thm:pre-dynkin-systems are made out of algebras}). Then, we define $\{ \mathcal{B}_o\}_{o \in O}$
    as the corresponding set of partitions. In particular, $\PreDhull(\mathcal{B}_o) = \mathcal{A}_o$. Thus, $\mathcal{D} = \bigcup_{o \in O} \PreDhull(\mathcal{B}_o) = \PreDhull(\bigcup_{o \in O} \mathcal{B}_o)$, because
    \begin{align*}
        \bigcup_{o \in O} \PreDhull(\mathcal{B}_o) \subseteq \PreDhull(\bigcup_{o \in O} \mathcal{B}_o)
        \subseteq \PreDhull(\mathcal{D})
        =\mathcal{D}
        = \bigcup_{o \in O} \mathcal{A}_o
        = \bigcup_{o \in O} \PreDhull(\mathcal{B}_o).
    \end{align*}
    
    For every partition $\mathcal{B}_o$ we have one $B_o \in \mathcal{B}_o$ such that $i \in B_o$.
    Obviously, $B_o \cap A \neq \emptyset$. But as well, $B_o \cap A^c \neq \emptyset$. Otherwise,
    $B_o \subseteq A$, thus $B_o = \emptyset$, because $A$ is a weak atom which is not contained in the
    pre-Dynkin-system $\mathcal{D}$. Thus, there is
    at least $j_o \in B_o$ with $j_o \notin A$.
    For one such $j_o$ we define the probability
    \begin{align*}
        \nu_{j_o} \coloneqq \psi_{j_o} - \frac{\epsilon}{|\{ j_o \colon o \in O\}|}.
    \end{align*}
    By choice of $\epsilon$ and $1 \le |\{ j_o \colon o \in O\}| \le \infty$ we have $\nu_{j_o} \ge 0$. Observe, we divide $\epsilon$ by the number of atoms on which we decrease the probability to guarantee normalization of $\nu$.
    Finally, for all $l \in [n] \setminus ( \{ j_o \colon o \in O\} \cup A)$ we set $\nu_l = \psi_l$.
    This assignment leads to the following conclusions:
    $\nu(B_o) = \psi(B_o)$ for every $o \in O$. More generally, for
    all $B \in \bigcup_{o \in O} \mathcal{B}_o$ we have $\nu(B) = \psi(B)$.
    
    The probability distribution $\nu$ is uniquely determined by the probability on the partitions,
    because $\mathcal{D} = \PreDhull(\bigcup_{o \in O} \mathcal{B}_o)$ (see above) and the
    $\pi$-$\lambda$-Theorem \citep[Lemma A.1.3]{williams1991probability}.
    So, $\nu \in \mba(\mathcal{D})$.
    But $\nu(A) = \nu(\{ i\}) + \nu(A \setminus \{ i\}) = \psi(\{ i\}) + \psi(A \setminus \{ i\}) + \epsilon \neq \psi(A)$.
    Thus, $A \notin \mdba (\mba(\mathcal{D}))$.
\end{proof}

\subsection{Intersection of Linear Subspaces}
\begin{lemma}[Intersection of Simple Gamble Spaces]
\label{lemma:intersection of simple gamble spaces}
    Let $\mathcal{A}, \mathcal{B} \subseteq 2^\Omega$ be two algebras. Then,
    \begin{align*}
        \operatorname{S}(\Omega, \mathcal{A}) \cap \operatorname{S}(\Omega, \mathcal{B}) = \operatorname{S}(\Omega, \mathcal{A} \cap \mathcal{B}).
    \end{align*}
\end{lemma}
\begin{proof}
    It is clear that $\operatorname{S}(\Omega, \mathcal{A} \cap \mathcal{B}) \subseteq \operatorname{S}(\Omega, \mathcal{A})$ and $\operatorname{S}(\Omega, \mathcal{A} \cap \mathcal{B}) \subseteq \operatorname{S}(\Omega, \mathcal{B})$, which makes the right to left set inclusion obvious.
    For the left to right inclusion observe that a function $f \in \operatorname{S}(\Omega, \mathcal{A}) \cap \operatorname{S}(\Omega, \mathcal{B})$ is a linear combination of indicator gambles of disjoint events $A_1, \ldots, A_n \in \mathcal{A}$, $\alpha_1, \ldots, \alpha_n \in \Reals$ and a linear combination of indicator gambles of disjoint events $B_1, \ldots, B_m \in \mathcal{B}$, $\beta_1, \ldots, \beta_m \in \Reals$. Without loss of generality we can demand that all $\alpha_1, \ldots, \alpha_n$ (respectively $\beta_1, \ldots, \beta_m$) have to be pairwise distinct, because we can replace indicator gambles scaled with the same factor by the appropriately scaled indicator function of the disjoint union of the events. Then, there is $i \in [n]$ for every $j \in [m]$ such that $\alpha_i = \beta_j$ and vice-versa. More importantly, there is $i \in [n]$ for every $j \in [m]$ such that $A_i = B_j$ and vice-versa. Hence, $\{ A_1, \ldots, A_n \} = \{ B_1, \ldots , B_m\} \in\mathcal{A} \cap \mathcal{B}$.
\end{proof}
\begin{lemma}[Intersection of Measurable Gamble Spaces]
\label{lemma:intersection of measurable gamble spaces}
    Let $\mathcal{A}_\sigma, \mathcal{B}_\sigma \subseteq 2^\Omega$ be two $\sigma$-algebras. Then,
    \begin{align*}
        \operatorname{B}(\Omega, \mathcal{A} _\sigma) \cap \operatorname{B}(\Omega, \mathcal{B}_\sigma) = \operatorname{B}(\Omega, \mathcal{A}_\sigma \cap \mathcal{B}_\sigma).
    \end{align*}
\end{lemma}
\begin{proof}
    It is clear that $\operatorname{B}(\Omega, \mathcal{A}_\sigma \cap \mathcal{B}_\sigma) \subseteq \operatorname{B}(\Omega, \mathcal{A}_\sigma)$ and $\operatorname{B}(\Omega, \mathcal{A}_\sigma \cap \mathcal{B}_\sigma) \subseteq \operatorname{B}(\Omega, \mathcal{B}_\sigma)$, which makes the right to left set inclusion obvious.
    For the left to right inclusion observe that for a function $f \in \operatorname{B}(\Omega, \mathcal{A}_\sigma) \cap \operatorname{B}(\Omega, \mathcal{B}_\sigma)$ the preimage of any Borel-measurable set in $\Reals$ is contained in $\mathcal{A}_\sigma$ and $\mathcal{B}_\sigma$.
\end{proof}

\section{Names of Pre-Dynkin-Systems and Dynkin-Systems}
\label{names of dynkin-systems}
See Table~\ref{tab:names for pre-dynkin-systems} for a list of names for pre-Dynkin-systems. Table~\ref{tab:names for dynkin-systems}
summarizes a list of names for Dynkin-systems.
\begin{table}[ht]
      \begin{tabular}{l}
         pre-Dynkin-system \citep{schurz2008finitistic}\\ \hline
         additive-class \citep[page 2]{rao1983theory}\\
         concrete logic \citep{ovchinnikov1999measures, de2007extending}\\
         partial field \citep{godowski1981varieties}\\
         quantum-mechanical algebra \citep{suppes1966probabilistic}\\
         semi-algebra \citep[page 13]{khrennikov2009interpretations}\\
         set-representable orthomodular poset \citep{ptak1998some}
    \end{tabular}
    \caption{A summary of names for pre-Dynkin-systems found in literature.}
    \label{tab:names for pre-dynkin-systems}
    \vspace{20pt}
    \begin{tabular}{l}
         Dynkin-system \citep{jun2011tighter}\\ \hline
         d-system \citep[page 193]{williams1991probability}\\
         $\lambda$-class \citep[page 7]{chow2003probability}\\
         quantum-mechanical $\sigma$-algebra \citep{suppes1966probabilistic}\\
         $\sigma$-class \citep{gudder1984extension}
    \end{tabular}
    \caption{A summary of names for Dynkin-systems found in literature.}
    \label{tab:names for dynkin-systems}
\end{table}

\section{Credal Sets of Pre-Dynkin-System Probabilities -- Credal Sets of Distorted Probabilities}
\label{sec:Credal Sets of Pre-Dynkin-system probabilities - Credal Sets of Distorted Probabilities}
To the best of the authors knowledge, there has not been any
attempt to parametrize a family of imprecise probabilities
via the set of induced precise probabilities.
In fact, a much better known class of imprecise probabilities
is parametrized via distortion functions \citep{wirch2001distortion}.
We use distorted probability functions as they regularly occur as
examples of imprecise probabilities \citep{walley1991statistical, wirch2001distortion, frohlich2022risk}. In particular,
there is a one-to-one correspondence of distorted probabilities as defined in the following
and so-called spectral risk measures, an important class of coherent upper previsions often used
in economics and finance \citep{frohlich2022risk}.
\begin{definition}[Credal Set of Distorted Probability]
\label{def:credal set for distorted probability}
    Let $\gamma \colon [0,1] \rightarrow [0,1]$ be a concave, increasing function
with $\gamma(0) = 0$ and $\gamma(1) = 1$. Let $\psi$ denotes a finitely additive probability on an algebra
$2^\Omega$ on $\Omega$. We overload the definition of a credal set
\begin{align*}
    \Mba(\psi, \gamma) \coloneqq \{ \nu \in \Pba\colon \nu(F) \le \gamma (\psi(F)), \ \forall F \in 2^\Omega\}.
\end{align*}
\end{definition}
How does the credal set of probabilities for distorted probabilities
relate to the credal set of probabilities for a
probability defined on a Dynkin-system?
\begin{lemma}[Distortion Lemma]
\label{lemma:distortion lemma}
    Let $\gamma \colon [0,1] \rightarrow [0,1]$ be a concave, increasing function
    with $\gamma(0) = 0$ and $\gamma(1) = 1$. If $\gamma$ is not the identity function, then $\gamma(x) > x$ for all $x \in (0,1)$.
\end{lemma}
\begin{proof}
    As $\gamma$ is concave, it is, in particular, quasi-concave.
    Thus, $x \mapsto \frac{\gamma(x)}{x}$ is decreasing \citep[Definition 10.1.1]{rubshtein2016foundations}.
    This gives the following inequalities
    \begin{align*}
        \frac{\gamma(x)}{x} \ge \frac{\gamma(x')}{x'} \ge \frac{\gamma(1)}{1} = 1,
    \end{align*}
    for $0 < x \le x'$.
    This implies, if there were $x \in (0,1)$ such that $\gamma(x) = x$, then $\gamma$
    would be the identity function. We excluded this by assumption, thus it follows $\gamma(x) > x$ for all $x \in (0,1)$.
\end{proof}
In the following Proposition we show that the set of events on which all measures of a credal set $\Mba(\psi, \gamma)$ coincide forms a pre-Dynkin-system. Actually, it is the \emph{system of certainty}, \ie the set of all events which get assigned either the value $0$ or the value $1$. We reuse the notation of the dual credal set function $\mdba$ which, as we noted earlier, maps a set of probabilities to the set structure on which those probabilities coincide.
\begin{proposition}[Events of Precise Probability for Distorted Probabilities]
\label{prop:events of precise probability for distorted probabilities}
    Let $\gamma \colon [0,1] \rightarrow [0,1]$ be a concave, increasing function
    with $\gamma(0) = 0$ and $\gamma(1) = 1$, which is not the identity function. Let $\psi$ denote a finitely additive probability on an algebra
    $2^\Omega$ on $\Omega$. 
    Let $\Mba(\psi, \gamma)$ be the credal set (Definition~\ref{def:credal set for distorted probability})
    and $\mdba$ the dual credal set function (Definition~\ref{def:Dual Credal Set Function}). Let
    $\mathcal{F}_{0} \coloneqq \{ F \in 2^\Omega \colon \psi(F) = 0\}$ denote the set of all measure zero sets and $\mathcal{F}_{01} \coloneqq \{ F \in 2^\Omega \colon \psi(F) = 0 \text{ or } \psi(F) = 1\}$ denote the set of all measure zero or one sets.
    Then
    \begin{align*}
        \mdba (\Mba(\psi, \gamma)) = \PreDhull(\mathcal{F}_{0}) = \mathcal{F}_{01}.
    \end{align*}
\end{proposition}
\begin{proof}
We show the equality of all sets via a circular set containment.
\begin{description}
    \item [ (a) ] $\mathcal{F}_{01} \subseteq \PreDhull(\mathcal{F}_{0})$\\
    If $F \in 2^\Omega$ such that $\psi(F) = 0$, then clearly $F \in \mathcal{F}_{0} \subseteq \PreDhull(\mathcal{F}_{0})$.
    If $F \in 2^\Omega$ such that $\psi(F) = 1$, then $\psi(F^c) = 0$. Thus, $F^c \in \mathcal{F}_{0}$. Because pre-Dynkin-systems
    are closed under complement $F \in \PreDhull(\mathcal{F}_{0})$.

    \item  [ (b) ] $\PreDhull(\mathcal{F}_{0}) \subseteq \mdba (\Mba(\psi, \gamma))$\\
    First, we write out
    \begin{align*}
        \mdba (\Mba(\psi, \gamma)) = \{ F \in 2^\Omega \colon \nu(F) = \psi(F), \ \forall \nu \in \{ \nu' \in \Pba\colon \nu'(F) \le \gamma ( \psi (F)), \ \forall F \in 2^\Omega\}\}
    \end{align*}
    Since $\gamma(0) = 0$, it is easy to see that $\mathcal{F}_{0} \subseteq \mdba (\Mba(\psi, \gamma))$. By Proposition~\ref{prop:Dual Credal Set Function Maps to Pre-Dynkin-Systems} we know that $\mdba (\Mba(\psi, \gamma))$
    is a pre-Dynkin-system. Thus, $\PreDhull(\mathcal{F}_{0}) \subseteq \mdba (\Mba(\psi, \gamma))$.

    \item [ (c) ] $\mdba (\Mba(\psi, \gamma)) \subseteq \mathcal{F}_{01}$\\
    We show this set inclusion via contraposition. If $F \in 2^\Omega$ has measure $\psi(F) \in (0,1)$, then $F \notin \mdba (\Mba(\psi, \gamma))$.
    For this we have to argue that there is a measure $\nu_F \in \Mba(\psi, \gamma)$ for every $F \in 2^\Omega$ with $\psi(F) \in (0,1)$ such that $\nu_F(F) \neq \psi(F)$.

    Observe that $\gamma \circ \psi$ defines a normalized, monotone, submodular set function on $2^\Omega$ \citep[page 17]{denneberg1994non}.
    Furthermore, any normalized, monotone, submodular set function induces a translation equivariant, monotone, positively homogeneous and subadditive functional $L_{\gamma \circ \psi}$ on all
    $f \in \SpaceOfGambles$ such that $L_{\gamma \circ \psi}(\chi_F) = (\gamma \circ \psi)(F)$ for all $F \in 2^\Omega$ \citep[page 260]{huber1981robust}\citep[page 130]{walley1991statistical}\citep[Proposition 5.1, Theorem 6.3]{denneberg1994non}. Hence, $L_{\gamma \circ \psi}$ is a coherent upper prevision \citep[page 65]{walley1991statistical}.
    Thus, Walley's extreme point theorem applies \citep[Theorem 3.6.2 (c)]{walley1991statistical}.\footnote{Even though this theorem is stated in terms of coherent lower previsions, it applies to coherent upper previsions, too. The weak$^\star$-compactness of the credal set $\Mba(\psi, \gamma)$, which is given by the coherence of $L_{\gamma \circ \psi}$ and \citep[Theorem 3.6.1]{walley1991statistical}, is crucial.}
    For any function $f\in \SpaceOfGambles$, in particular any $\chi_F$ with $F \in 2^\Omega$, there is a linear prevision $f \mapsto \langle f, \nu \rangle$ with $\nu \in \Pba$ on $\SpaceOfGambles$ dominated by $L_{\gamma \circ \psi}$ such that $\langle \chi_F, \nu \rangle = L_{\gamma \circ \psi}(\chi_F)$. More concretely, for any $F \in 2^\Omega$ there is a $\nu_F \in \Mba(\psi, \gamma)$ such that $\nu_F(F) = (\gamma \circ \psi)(F)$.
    If $\psi(F) \in (0,1)$ then Lemma~\ref{lemma:distortion lemma} applies and $\nu_F(F) = (\gamma \circ \psi)(F) = \gamma( \psi(F)) > \psi(F)$ gives the desired inequality. In conclusion, there is no $F \in 2^\Omega$ with measure $\psi(F) \in (0,1)$ such that $F \in \mdba (\Mba(\psi, \gamma))$. This implication finalizes the proof.
\end{description}
\end{proof}
We call the set $\mathcal{F}_{01}$ a \emph{system of certainty} for $\psi$. Hence,
the system of precision of a distorted probability
is the system of certainty.
Since the proposition clarifies the relation from distorted
imprecise probability to probability on Dynkin-system, the reverse
question immediately follows: when is $\Mba(\psi|_\mathcal{D}, \mathcal{D}) \subseteq \Mba(\psi, \gamma)$? In other words, given an extendable probability defined on a pre-Dynkin-system, what is a distortion function of an extension of this probability such that the credal set of the former is contained in the credal set of the latter. The simple example below gives an instantiation of this problem for which it is easy to find a solution. However, the problem is harder for more general cases.
{
\color{darkgray}
\begin{example}
    Let $\Omega = \{ 1,2\}$ and $\psi(1) = \psi(2)= 0.5$. Let  $\mathcal{D} = \{ \emptyset, \Omega\}$. The for the admittedly extreme distortion function $\gamma(0) = 0$, otherwise $\gamma(a) = 1$, for $a \in (0,1]$ the credal set $M(\psi, \gamma) = \Pba$ is the set of all probability measures. Hence, $M(\psi|_{\mathcal{D}}, \mathcal{D}) = \Pba \subseteq M(\psi, \gamma)$.
\end{example}
}

\section{Dynkin-Systems and Countably Additive Probability}
\label{appendix:Dynkin-Systems and Countably Additive Probability}
In the main text, we presented the majority of the results for the general case: pre-Dynkin-systems and finitely additive probabilities.
We did so, as we emphasize clearly here, not for the sake of mathematical generality.
There are probabilistic problems which demand for the use of finitely additive probabilities, \eg von Mises frequentistic notion of probability \citep{mises1964mathematical, schurz2008finitistic}.
As pre-Dynkin-systems and finitely additive probabilities walk hand in hand, so too do
Dynkin-systems and countably additive probabilities.
It is possible to strengthen some results when we assume Dynkin-systems and countably
additive probabilities. Furthermore, countably additive probabilities are more familiar
to students of probability theory. Finitely additive probabilities still
eke out an exotic living \citep{rao1983theory,kadane1999statistical, chichilnisky2010foundations}. Nevertheless, they are central in large parts of literature on imprecise probability \citep{walley1991statistical, augustin2014introduction}.

\subsection{Dynkin-Systems}
\label{appendix:Dynkin-Systems}
Many of the following results are analogous to the statements
for pre-Dynkin-systems. The subsequent theorem is analogous to
Theorem~\ref{thm:pre-dynkin-systems are made out of algebras}.
We remark, however, that we cannot use the same proof-technique as in Theorem~\ref{thm:pre-dynkin-systems are made out of algebras}
because the generalization of Lemma~\ref{lemma:Supremum of a Chain of Algebras is an Algebra} to $\sigma$-algebras
does not hold. Hence, we use a different technique to prove Theorem~\ref{thm:dynkin-systems are made out of sigma-algebras}
for Dynkin-systems and $\sigma$-algebras.
\begin{theorem}[Dynkin-Systems are made out of $\sigma$-Algebras]
\label{thm:dynkin-systems are made out of sigma-algebras}
    Let $\mathcal{D}_\sigma$ be a Dynkin-system on an arbitrary set $\Omega$. Then there is a unique
    family of maximal $\sigma$-algebras $\{\mathcal{A}_i \}_{i \in I}$
    such that $\mathcal{D}_\sigma = \bigcup_{i \in I} \mathcal{A}_i$. We
    call these $\sigma$-algebras the $\sigma$\emph{-blocks}
    of $\mathcal{D}_\sigma$.
\end{theorem}
\begin{proof}
    Since Dynkin-systems are pre-Dynkin-systems 
    Theorem~\ref{thm:pre-dynkin-systems are made out of algebras} guarantees that
    $\mathcal{D}_\sigma$ is constituted of a set of algebras $\{\mathcal{A}_i \}_{i \in I}$.
    Each algebra $\mathcal{A}_i$ is closed under finite intersection. Thus, it
    is a ``compatible collection'' following the terms of \citet{gudder1973generalized}.
    It follows that $\mathcal{A}_i$ is contained in a sub-$\sigma$-algebra of $\mathcal{D}_\sigma$
    by Theorem 2.1 in \citep{gudder1973generalized}. As any sub-$\sigma$-algebra is 
    an algebra and $\mathcal{A}_i$ is maximal, \ie there is no algebra contained in $\mathcal{D}$
    such that $\mathcal{A}_i$ is a strict sub-algebra of this algebra, $\mathcal{A}_i$ itself is a $\sigma$-algebra.
    Hence, $\mathcal{A}_i$ are the $\sigma$-blocks of $\mathcal{D}_\sigma$.
\end{proof}
Not every union of $\sigma$-algebras is a Dynkin-system (\cf Appendix~\ref{appendix:compatibility structure}). But if a union of $\sigma$-algebras forms a compatibility structure, then it is a Dynkin-system.
\begin{theorem}[Union of $\sigma$-Algebras is Dynkin-system if Compatibility Structure]
\label{theorem:Union of sigma-Algebras is Dynkin-system if Compatibility Structure}
    Let $\{\mathcal{A}_i \}_{i \in I}$ be a family of $\sigma$-algebras on an arbitrary $\Omega$.
    If $\{\mathcal{A}_i \}_{i \in I}$
    is a compatibility structure, then $\mathcal{A}_\cup \coloneqq \bigcup_{i \in I} \mathcal{A}_i$ is a Dynkin-system 
\end{theorem}
\begin{proof}
    We show that
    $\mathcal{A}_\cup$ is a Dynkin-system if it is the union
    of the compatibility structure. First,
    it contains the empty set. Second, if $A \in \mathcal{A}_\cup$
    there is an $i$ such that $A \in \mathcal{A}_i$. Then $A^c \in \mathcal{A}_i$, thus
    $A^c \in \mathcal{A}_\cup$. Third, let $\{A_j\}_{j \in J}$ be a
    subset of $\mathcal{A}_\cup$ with $J \subseteq \Naturals$ such
    that $A_k \neq \emptyset$, $A_l \neq \emptyset$ and
    $A_k \cap A_l = \emptyset$ for all $k \neq l; k,l \in J$.
    Thus, $\{A_j\}_{j \in J} \cup \{ \emptyset\}$ is closed under finite intersection. It is thus a $\pi$-system. By definition of a compatibility
    structure there is an $i \in I$ such that $\{A_j\}_{j \in J} \cup \{\emptyset\} \subseteq \mathcal{A}_i$.
    Since $\mathcal{A}_i$ is a $\sigma$-algebra, it is closed under countable disjoint union. Thus,
    $\bigcup_{j \in J} A_j \in \mathcal{A}_i \subseteq \mathcal{A}_\cup$.
\end{proof}

\subsection{Technical Setup}
In order to work on firm ground when introducing countably additive probabilities, we change our basic technical setup.
We summarize the used notations in Table~\ref{tab:notation countable}.
\begin{table}[t]
    \centering
    \begin{tabular}{r|l}
        $\Omega$ & Polish Space\\
        $\mathcal{D}_\sigma $ & Dynkin-system on $\Omega$ (Definition~\ref{def:dynkin-system})\\
        $\mu_\sigma$ & Countably additive probability defined on $\mathcal{D}_\sigma$ (Definition~\ref{def:probameasure on Dynkin-system})\\
        $\sigma(\mathcal{A})$ & $\sigma$-algebra hull of set system $\mathcal{A} \subseteq 2^\Omega$\\
        $\mathcal{F}_\sigma$ & Borel-$\sigma$-algebra on $\Omega$\\
        $\operatorname{ca}(\Omega, \mathcal{F}_\sigma)$ & Set of bounded, signed, countably additive measures on $\mathcal{F}_\sigma$\\
        $\Pca(\Omega, \mathcal{F}_\sigma)$ respectively $\Pca$ & Set of countably additive probability measures on $\mathcal{F}_\sigma$\\
        $\Mca(\mu_\sigma, \mathcal{D}_\sigma)$ & $\sigma$-Credal set of $\mu_\sigma$ on $\mathcal{D}_\sigma$ (Proposition~\ref{prop:credal sets of probabilities on dynkin-systems})\\
        $\underline{\mu}_{\mathcal{D}_\sigma}$, $\overline{\mu}_{\mathcal{D}_\sigma}$ & Lower respectively upper coherent $\sigma$-extension (Proposition~\ref{prop:coherent sigma-extension of probbaility})\\
    \end{tabular}
    \caption{Summary of used notations in Appendix~\ref{appendix:Dynkin-Systems and Countably Additive Probability}.}
    \label{tab:notation countable}
\end{table}
First, we assume now that $\Omega$ is a Polish space, that is, a separable completely metrizable
topological space \citep[page 20]{huber1981robust}.
We denote by $\mathcal{F}_\sigma$ the Borel-$\sigma$-algebra
with respect to the given topology.
With $\operatorname{ca}(\Omega, \mathcal{F}_\sigma)$ we denote the space of all finite, signed countably additive
measures on $(\Omega, \mathcal{F}_\sigma)$ (\cf \citep[page 20]{huber1981robust}).

We define $\Pca \subseteq \operatorname{ca}(\Omega, \mathcal{F}_\sigma)$ as the set of all 
probability measures on $(\Omega, \mathcal{F}_\sigma)$.
Every measure in $\Pca$ is regular \citep[page 20]{huber1981robust}. The set
$\Pca$ is equipped with the weak topology \citep[page 21]{huber1981robust}. This
makes $\Pca$ a Polish space \citep[page 29]{huber1981robust}. We emphasize that this is a different topology in comparison to the topology used above. The weak$^\star$ topology is weaker than the weak topology. And the weak topology is in this case metrizable.

Lastly, we denote the smallest $\sigma$-algebra which contains $\mathcal{A}$ with $\sigma(\mathcal{A})$.
We say that $\sigma(\mathcal{A})$ is the $\sigma$-algebra hull of $\mathcal{A}$.

\subsection{Dynkin Probability Spaces}
The definition of probability measures on pre-Dynkin-systems directly applies to Dynkin-systems, too.
Analogous to Kolmogorov's probability space, we can now leverage Definition~\ref{def:probameasure on Dynkin-system}
for a definition of a Dynkin probability space.
\begin{definition}[Dynkin Probability Space]\citep[page 296]{gudder1969quantum}
\label{def:Dynkin probability space}
    The triple $(\Omega, \mathcal{D}_\sigma, \mu_\sigma)$ is called a \emph{Dynkin probability space} if and only if (a) $\Omega$ is a non-empty base space,
    (b) $\mathcal{D}_\sigma$ is a Dynkin-system on $\Omega$ and (c) $\mu_\sigma$ is a countably additive probability measure on $\mathcal{D}_\sigma$ following Definition~\ref{def:probameasure on Dynkin-system}.
\end{definition}
Dynkin probability spaces\footnote{Dynkin probability spaces have been called ``quantum probability spaces'' in quantum probability theory \citep{gudder1969quantum}.}
generalize Kolmogorov's probability spaces. If $\mathcal{D}_\sigma$ were a $\sigma$-algebra,
then the Dynkin probability space would become a classical probability space.
Theorem~\ref{thm:dynkin-systems are made out of sigma-algebras} provides
another interesting link between Dynkin and classical probability spaces.
\begin{proposition}
\label{prop:dynkingivesmultikolmogorov}
    Every Dynkin probability space $(\Omega, \mathcal{D}_\sigma, \mu_\sigma)$ defines a collection
    of classical probability space  $\{(\Omega, \mathcal{A}_i, \mu_i)\}_{i \in I}$
    where $\mathcal{D}_\sigma = \bigcup_{i \in I} \mathcal{A}_i$ and
    $\mu_i \coloneqq \mu_\sigma|_{\mathcal{A}_i}$ are \emph{consistent}, \ie 
    $\mu_i(A) = \mu_{j}(A)$ for any $A \in \mathcal{A}_i \cap \mathcal{A}_j$ and $i,j \in I$.
\end{proposition}
Such ``Multi-Kolmogorov'' probability spaces have similarly been formalized in \citep[page 32]{khrennikov2009contextual} or
\citep[page 154]{vorob1962consistent}.
The multiplicity of probability spaces has been interpreted as 
a collection of contexts. Each context possesses its own classical probability space.\footnote{
In quantum physics, we obtain classical behavior in single context and quantum behavior across context.
In machine learning, each marginal scenario gets equipped an own probability space,
\ie an own context (\cf\citep{cuadras2002distributions, eban2014discrete}).}
Again, as for pre-Dynkin-systems and finitely additive probabilities, the question
of embedding Dynkin probability spaces into classical probability spaces arises.

\subsection{Conditions for Extendability for Countably Additive Probabilities}
\begin{definition}[$\sigma$-Extendability]
    Let $\mathcal{F}_\sigma$ be the Borel-$\sigma$-algebra on a Polish space $\Omega$ and $\mathcal{D}_\sigma$ a Dynkin-system
    contained in this $\sigma$-algebra. With $\Pca$ we denote the set of all countably additive probability
    measures on $(\Omega, \mathcal{F}_\sigma)$. We call a countably additive probability $\mu_\sigma$ on $\mathcal{D}_\sigma$
    \emph{$\sigma$-extendable}\footnote{We call this extendability ``$\sigma$'' to emphasize
the difference to Definition~\ref{def:extendability} and the
 properties of countably additive probability measures.},
    if and only if there is a countably additive probability measure $\nu \in \Pca$ such that $\nu|_{\mathcal{D}_\sigma} = \mu_\sigma$.
\end{definition}
Extendability of Dynkin probability spaces (as well as for finitely additive probabilities on pre-Dynkin-systems)
has already been part of discussions in quantum probability since 1969 \citep{gudder1969quantum}
up to more current times \citep{de2010measures}.
Most of the results, \eg \citep{gudder1984extension, de2007extending, de2010measures}, are only stated in terms of
finitely additive probability measure and/or apply only on a structurally restricted 
class of Dynkin-systems.
A theorem due to \citet{maharam1972consistent} in the literature on marginal probabilities can be tweaked to give a universal sufficient and necessary criterion
for $\sigma$-extendability.
\begin{theorem}[$\sigma$-Extendability Condition for Countably Additive Probability]
    Let $\mathcal{F}_\sigma$ be the Borel-$\sigma$-algebra on a Polish space $\Omega$ and $\mathcal{D}_\sigma$ a Dynkin-system
    contained in this $\sigma$-algebra such that $\sigma(\mathcal{D}_\sigma) = \mathcal{F}_\sigma$. Let $\mu_\sigma$ be a countably additive probability on $\mathcal{D}_\sigma$.
    For each $\sigma$-block $\mathcal{A}_i$ of $\mathcal{D}_\sigma$ we assume that $\mu_i \coloneqq \mu_\sigma|_{\mathcal{A}_i}$ is inner regular\footnote{A measure $\mu$ is inner regular if $\mu(G) = \sup \{  \mu(K) \colon K \subseteq G, K \text{ compact}\}$ \citep[page 808]{schechter1997handbook}.}.
    The probability $\mu_\sigma$ is $\sigma$-extendable if and only if
    \begin{align*}
        \sum_{j = 1}^n f_j(\omega) \ge 0, \ \forall \omega \in \Omega \implies \sum_{j =1}^n \int_{\Omega} f_j(\omega) d \mu_{i_j}(\omega) \ge 0
    \end{align*}
    for all finite families of measurable simple gambles $\{ f_j \}_{j \in [n]} \subseteq S(\Omega, \mathcal{D}_\sigma)$.
\end{theorem}
\begin{proof}
Proposition~\ref{prop:dynkingivesmultikolmogorov} separates the Dynkin probability space
into a set of Kolmogorov probability spaces. Then 
Theorem 8.1 in \citep{maharam1972consistent}
applies. As technical requirements, we emphasize that $\Omega$ is a Hausdorff space and
all restrictions $\mu_i$ are inner regular.
Thus, there exists a 
probability measure on the algebra generated by $\mathcal{D}_\sigma$.
Since the probability measure is countably additive, 
Caratheodory's Extension Theorem \citep[Theorem 1.7]{williams1991probability}
states that it can be uniquely extended to the Borel-$\sigma$-algebra $\mathcal{F}_\sigma = \sigma(\mathcal{D}_\sigma)$.
\end{proof}
Similar and related results can be found in \citep{vorob1962consistent, kellerer1964verteilungsfunktionen, horn1948measures}.
In particular, every so-called marginal scenario can be represented as a countably additive probability on a Dynkin-system (\cf \citep{vorob1962consistent} \citep[Example 4.2]{gudder1984extension}).
Marginal scenarios often arise in practical setups. They are defined as settings in which for a collection of random variables only specific partial joint distributions of the random variables are given. For those scenarios there are purely combinatorial conditions on the Dynkin-system sufficient
for the extendability \citep{vorob1962consistent, kellerer1964verteilungsfunktionen}.

    

\subsection{Credal Set of Countably Additive Probabilities on Dynkin-Systems}
The credal set of probabilities which we defined in Corollary~\ref{corollary:Coherent Extension of Probability}
contains finitely additive probabilities, some of which might not be countably additive.
For this reason we redefine the credal set for countably additive probabilities.
The credal set is the set of all countably additive probabilities which coincide with the
reference probability $\mu_\sigma$ on $\mathcal{D}_\sigma$.
\begin{proposition}[$\sigma$-Credal Set for Probabilities on Dynkin-Systems]
\label{prop:credal sets of probabilities on dynkin-systems}
    Let $\mathcal{F}_\sigma$ be the Borel-$\sigma$-algebra on a Polish space $\Omega$ and $\mathcal{D}_\sigma$ a Dynkin-system
    contained in this $\sigma$-algebra. Let $\mu_\sigma$ be a countably additive probability on $\mathcal{D}_\sigma$ and $\Pca$ the set of all countably additive probabilities on $\mathcal{F}_\sigma$.
    We call
    \begin{align*}
        \Mca(\mu_\sigma, \mathcal{D}_\sigma) \coloneqq \{ \nu \in \Pca\colon \nu(A) = \mu_\sigma(A), \ \forall A \in \mathcal{D}_\sigma\},
    \end{align*}
    the \emph{$\sigma$-credal set of $\mu_\sigma$ on $\mathcal{D}_\sigma$}. If $\Mca(\mu_\sigma, \mathcal{D}_\sigma) \neq \emptyset$, then it is weak closed and convex.
\end{proposition}
\begin{proof}
\begin{description}
    \item [Convexity]
    Let $\{ \nu_i \}_{i \in [n]} \subseteq \Mca(\mu_\sigma, \mathcal{D}_\sigma)$ and $\alpha_i \in [0,1]$ for all $i \in [n]$
    with $\sum_{i = 1}^n \alpha_i$, then $\sum_{i = 1}^n \alpha_i \nu_i \in \Mca(\mu_\sigma, \mathcal{D}_\sigma)$,
    because
    \begin{align*}
        \sum_{i = 1}^n \alpha_i \nu_i(A) = \sum_{i = 1}^n \alpha_i \mu_\sigma(A)= \mu_\sigma(A),\qquad  \forall A \in \mathcal{D}_\sigma.
    \end{align*}
    
    \item [Closedness]
    We assumed $\Omega$ to be a Polish space. It follows
    that $\Pca$ itself is a Polish space \citep[page 29]{huber1981robust}.
    Every Polish space is metrizable, consequently Lemma 21.2 in \citep{munkres2014topology} applies.
    Closedness of a set $\Q \subseteq \Pca$ can be identified
    via the convergence of sequences in $\Q$.

    Let $\nu_1, \nu_2 \ldots \in \Mca(\mu_\sigma, \mathcal{D}_\sigma)$ be a sequence of probability measures such
    that $\lim_{n\rightarrow \infty } \nu_n = \nu$, then 
    \begin{align*}
        \nu(A) = \lim_{n \rightarrow \infty}\nu_n (A) = \lim_{n \rightarrow \infty} \mu_\sigma (A) = \mu_\sigma(A),\qquad  \forall A \in \mathcal{D}_\sigma.
    \end{align*}
    It follows that $\nu \in \Mca(\mu_\sigma, \mathcal{D}_\sigma)$.
\end{description}
\end{proof}
Rather obviously the following corollary holds.
\begin{corollary}[$\sigma$-Extendability and $\sigma$-Credal Set]
    Let $\mathcal{F}_\sigma$ be the Borel-$\sigma$-algebra on a Polish space $\Omega$ and $\mathcal{D}_\sigma$ a Dynkin-system
    contained in this $\sigma$-algebra. The countably additive probability $\mu_\sigma$ on $\mathcal{D}_\sigma$ is $\sigma$-extendable
    if and only if $\Mca(\mu_\sigma, \mathcal{D}_\sigma) \neq \emptyset$.
\end{corollary}
Thus, we can redefine the lower and upper coherent extension if the Dynkin probability space is continuously extendable.
This extension is then derived from a set of countably additive probability measures.
\begin{proposition}[Coherent $\sigma$-Extension of Probability]
\label{prop:coherent sigma-extension of probbaility}
    Let $\mathcal{F}_\sigma$ be the Borel-$\sigma$-algebra on a Polish space $\Omega$ and $\mathcal{D}_\sigma$ a Dynkin-system
    contained in this $\sigma$-algebra. Assume the countably additive probability $\mu_\sigma$ on $\mathcal{D}_\sigma$ is $\sigma$-extendable. Then
    \begin{align*}
        \underline{\mu}_{\mathcal{D}_\sigma}(A) \coloneqq \inf_{\nu \in \Mca(\mu_\sigma, \mathcal{D}_\sigma)}  \nu(A),
        \quad
        \overline{\mu}_{\mathcal{D}_\sigma}(A) \coloneqq \sup_{\nu \in \Mca(\mu_\sigma, \mathcal{D}_\sigma)}  \nu(A),  \qquad \forall A \in \mathcal{F}_\sigma,
    \end{align*}
    define a coherent lower (respectively upper) probability on $\mathcal{F}_\sigma$ in the sense of \citep{walley1991statistical}.
\end{proposition}
\begin{proof}
    First, we notice the every countably additive probability is also finitely additive.
    Via Theorem~\ref{thm:Extendability Gives Coherence} and application of Theorem 3.3.4 (b) in \citep{walley1991statistical} we directly obtain the result.
\end{proof}
The coherent extension is one of two methods of measure extension presented
in this paper. Again, we ask how the coherent $\sigma$-extension relates to the inner
and outer measure construction.
\begin{corollary}[Extension Theorem -- Countably Additive Case]
\label{corollary:Extension Theorem - countably additive case}
     Let $\mathcal{F}_\sigma$ be the Borel-$\sigma$-algebra on a Polish space $\Omega$ and $\mathcal{D}_\sigma$ a Dynkin-system
    contained in this $\sigma$-algebra. Suppose the countably additive probability $\mu_\sigma$ on $\mathcal{D}_\sigma$ is $\sigma$-extendable. Then
    \begin{align*}
         \mu_*(A) \le \underline{\mu}_{\mathcal{D}_\sigma}(A) \le \overline{\mu}_{\mathcal{D}_\sigma}(A) \le \mu^*(A), \qquad \forall A \in \mathcal{F}_\sigma.
    \end{align*}
\end{corollary}
\begin{proof}
    Since $\Mca(\mu_\sigma, \mathcal{D}_\sigma) \subseteq \Mba(\mu_\sigma, \mathcal{D}_\sigma)$, clearly
    \begin{align*}
        \inf_{\nu \in \Mca(\mu_\sigma, \mathcal{D}_\sigma)}  \nu(A) \ge \inf_{\nu \in \Mba(\mu_\sigma, \mathcal{D}_\sigma)}  \nu(A), \qquad \forall A \in \mathcal{F}_\sigma,\\
        \sup_{\nu \in \Mca(\mu_\sigma, \mathcal{D}_\sigma)}  \nu(A) \le \sup_{\nu \in \Mba(\mu_\sigma, \mathcal{D}_\sigma)}  \nu(A), \qquad \forall A \in \mathcal{F}_\sigma.
    \end{align*}
    Hence, the analogous Theorem~\ref{thm:Extension Theorem - Finitely Additive Case} gives the result.
\end{proof}

\section{From Set Systems to Logical Structures and Back}
\label{appendix:From Set Systems to Logical Structures and Back}
In the year 1936, Stone demonstrated a celebrated representation result of logical algebras \citep{Stone1936theory}.
In particular, he showed that any Boolean algebra can be equivalently
represented by an algebra of sets and vice-versa.
This representation results opened the door to many generalizations
of the domain of probabilities (see \citep{narens2016introduction} for a nice summary). For instance, Boolean algebras can be replaced by weaker logical structures,
such as orthomodular lattices. Interestingly, modular ortholattices and orthomodular lattices have been investigated as representational
structures for quantum theoretic descriptions and measurements \citep{Birkhoff1936Logic, husimi1937studies}.
Furthermore, it turned out that an analogue to Stone's representation
holds for some orthomodular lattices. Some orthomodular lattices are isomorphic to pre-Dynkin-systems \citep{godowski1981varieties}. On the other hand, 
every (pre-)Dynkin-system is isomorphic to an ($\sigma$-)orthomodular poset, a order-theoretic generalization of an orthomodular lattice (\cf\citep{brabec1979compatibility}).
In summary, the investigation of probabilities on general logical structures
parallels our work presented here. We, however, have stuck 
to set structures as the domain of probabilities.

\end{document}